\documentclass[12pt]{article}
\usepackage{amssymb,amsmath}

\DeclareMathAlphabet{\eufrak}{U}{}{}{}  
\SetMathAlphabet\eufrak{normal}{U}{euf}{m}{n}
\SetMathAlphabet\eufrak{bold}{U}{euf}{b}{n}

\allowdisplaybreaks

\numberwithin{equation}{section}
\newenvironment{Proof}{\removelastskip\par\medskip
\noindent{\em Proof.} \rm}{\penalty-20\null\hfill$\square$\par\medbreak}

\def\real{{\mathord{{\rm I\kern-2.8pt R}}}}        
\def\inte{{\mathord{{\rm I\kern-2.8pt N}}}}
\def\PP{{\mathord{{\rm I\kern-2.8pt P}}}}

\def\real{{\mathord{\mathbb R}}}

\def\inte{{\mathord{\mathbb N}}}

\def\Cov{{\mathrm{{\rm Cov}}}}
\def\Var{{\mathrm{{\rm Var}}}}
\def\constante{c}

\newcommand{\ind}{\mathbf{1}}

\newcommand{\disp}{\displaystyle}

\newcommand{\lto}{\longrightarrow}

\newcommand{\om}{\omega}

\def\rit{\mathbb{R}}

\def\E{\mathop{\hbox{\rm I\kern-0.20em E}}\nolimits}
\def\Var{\mathop{\hbox{\rm Var}}\nolimits}
\def\Cov{\mathop{\hbox{\rm Cov}}\nolimits}

\def\med{\hskip 10pt}

\def\sig{\sigma(S^{n-1})} 

\newtheorem{prop}{Proposition}[section]
\newtheorem{lemma}[prop]{Lemma}

\newtheorem{corollary}[prop]{Corollary}
\newtheorem{theorem}[prop]{Theorem}

\def\Dom{\rm Dom \ \! }

\title{Dimension free and 
 infinite variance tail estimates on Poisson space 
} 

\author{Jean-Christophe Breton \and 
 Christian Houdr\'e
\and
Nicolas Privault}
\date{}

\begin{document}
\maketitle

\begin{abstract} 
 Concentration inequalities are obtained on Poisson space, 
 for random functionals 
 with finite or infinite variance. In particular, 
 dimension free tail estimates and exponential integrability results 
 are given for the Euclidean norm of vectors of independent functionals. 
 In the finite variance case these results are applied 
 to infinitely divisible random variables such as quadratic Wiener functionals, 
 including L\'evy's stochastic area and the square norm of Brownian paths. 
 In the infinite variance case, various tail estimates such as 
 stable ones are also presented. 
\end{abstract} 
 
\normalsize

\vspace{0.5cm}

\noindent {\bf Key words:} Concentration, 
 infinite divisibility, 
 stable laws, Poisson space, Ornstein-Uhlenbeck semi-group, 
 quadratic Wiener functionals, large deviations. 
\\ 
{\em Mathematics Subject Classification:} 
60F99, 60E07, 60G57. 
 
\small

\normalsize

\baselineskip0.7cm


\section{Introduction and notation} 
\label{3.0} 
\setcounter{equation}{0} 
 Let $\Omega^X$ denote the set of Radon measures 
$$\Omega^X = \left\{ 
 \omega = \sum_{i=1}^N \epsilon_{t_i} \ : \ 
 (t_i)_{i=1}^{i=N} \subset X, \ t_i\not=t_j, \ \forall i\not= j, \ 
 N\in \inte\cup \{ \infty \}\right\},$$ 
 where $X$ is a $\sigma$-compact metric space with distance $d_X$, 
 and $\epsilon_t$ denotes the Dirac measure at $t\in X$. 
 Let $\nu$ be a diffuse Radon measure on $X$, and let $P$ be 
 the Poisson measure with intensity $\nu$ on $\Omega^X$. 
 Let the linear, closable, finite difference operator 
$$ 
D:L^2(\Omega^X, P) 
 \longrightarrow L^2(\Omega^X\times X,P\otimes \nu)
$$ 
 be defined via 
$$ 
D_x F(\omega ) = F(\omega \cup \{ x \}) -F(\omega ), \ \ \ 
 \ \ \ dP\times \nu (d\omega , dx)\mbox{-a.e.}, 
$$ 
 where as a convention we identify $\omega \in \Omega^X$ with its support, 
 cf. e.g. \cite{nualart}, \cite{picard}, \cite{prirose}. 

 In \cite{ane}, \cite{wuls2}, \cite{hp}, Poisson tail estimates 
 are obtained under the hypothesis 
$$ 
 DF\leq K, \quad P\otimes \nu\mbox{-a.e.}, \quad \mbox{and} \quad 
 \|DF\|_{L^{\infty}(\Omega^X , L^2(X , \nu ))}\leq \tilde{\alpha} < \infty , 
$$ 
 for some $K\geq 0$. 
 While (modified) logarithmic Sobolev inequalities and the Herbst method are used 
 in \cite{ane} and \cite{wuls2}, 
 the methods of \cite{hp} rely on covariance representations (\cite{bgh}, \cite{h3}). 
 Recently the results of \cite{h3} have further led in \cite{HM} to 
 estimates for Lipschitz functions of stable random vectors. 
 Even more recently, dimension free concentration is obtained in \cite{HRB} 
 for the Euclidean norm as well as for various classes of functions 
 of independent infinitely divisible vectors having finite exponential moments. 
 
 In the present paper we first obtain new deviation inequalities on Poisson space 
 via the covariance method. Then, by replacing the bounds on $DF$ and 
 on $\|DF\|_{L^{\infty}(\Omega^X , L^2(X , \nu ))}$ by growth conditions, 
 deviation results for Poisson functionals with infinite variance are given. 

 Let us briefly describe the content of the paper. 
 In Section~\ref{sec:deviation-OU} we deal 
 with L\'evy measures with finite variance, 
 using the covariance representation method involving 
 the Ornstein-Uhlenbeck semi-group. 
 This leads to general deviation results for Poisson functionals 
 having finite exponential moments. 
 In Section~\ref{sec:appli-vect} we obtain dimension free deviation estimates 
 and exponential integrability properties for random vectors of 
 such Poisson functionals. 
 Since an infinitely divisible random vector can be represented as 
 a vector of Poisson stochastic integrals, these results are then applied 
 to derive deviation inequalities 
 for Lipschitz functions of infinitely divisible vectors. 
 In Section~\ref{sec:appli-quadratic}, 
 we study the particular case of quadratic Wiener functionals, 
 including the square norm of Brownian path, the sample variance of Brownian motion 
 and  L\'evy's stochastic area. 
 For such i.i.d. vectors, this also gives dimension free inequalities in Euclidean 
 norm, and large deviation estimates in $\ell^p$-norm, $p\in [1,\infty ]$, 
 recovering tail estimates of \cite{AG} for non-decoupled Gaussian chaos of degree $2$. 
 In Section~\ref{sec:no-variance} we adapt the method of \cite{HM} 
 to prove other tail estimates under weaker hypothesis on the gradient. 
 For example, if $\nu$ is the L\'evy measure of an $\alpha$-stable vector, 
 the bounds on $D$ can be replaced by the growth conditions 
\begin{equation} 
\label{growth} 
\sup_{x\in B_X(0,R)} \vert D_x F \vert \leq C' R \quad 
\mbox{and} \quad 
\Vert DF \Vert_{L^\infty(\Omega^X , L^2(B_X (0,R)))}^2 
 \leq C 
 R^{2-\alpha}, \quad R\geq R_0, 
\end{equation} 
 where $B_X (0,R) = \{ x\in X \ : \ d_X(0,x) \leq R \}$ is the 
 ball of radius $R$ in $X$. Here, $0$ denotes a fixed 
 arbitrary center in $X$, whose choice has no influence on the growth 
 conditions \eqref{growth}. 
 This leads to an estimate of stable type for the 
 deviation of $F$ from one of its medians. 

 Let us now introduce some notation which will be used throughout the 
 paper. 
 The multiple Poisson stochastic 
 integral $I_n (f_n)$ is defined as 
$$I_n (f_n) (\omega ) = 
 \int_{\Delta_n} 
 f_n(y_1,\ldots ,y_n) (\omega (dy_1)-\nu (dy_1)) 
 \cdots (\omega (dy_n)-\nu (dy_n)), 
$$ 
 for every square-integrable symmetric function 
 $f_n\in L^2 (X ,\nu )^{\circ n}$, 
 where 
$$\Delta_n = {\{(x_1,\ldots ,x_n)\in X^n \ : \ x_i\not= x_j , \ \forall i\not= j\}} 
. 
$$ 
 Recall the isometry formula 
$$E [I_n(f_n)I_m(g_m)] = n! 
 {\bf 1}_{\{ n=m \}} 
 \langle f_n , g_m\rangle_{L^2 (X,\nu )^{\circ n}},$$ 
 see \cite{nualartvives}, 
 and recall also that every square-integrable 
 random variable $F\in L^2(\Omega^X, P )$ admits 
 the Wiener-Poisson decomposition 
$$F=\sum_{n=0}^\infty I_n (f_n) 
. 
$$ 
 The operator $D$ defined above is such that 
$$ 
 D_x I_n (f_n) (\omega ) = n 
  I_{n-1}(f_n(*,x)) (\omega ) , 
 \ \ \ P (d\omega ) \otimes \nu (dx)\mbox{-a.e.}, 
 \ \ n\in \inte 
, 
$$ 
 and in particular, 
$$
D_x I_1 (f) (\omega ) = 
 f(x), \quad \nu (dx)\mbox{-a.e.}
$$ 
 We denote by $\Dom (D)$ the domain of $D$, i.e. the space of functionals $F\in L^2(\Omega^X , P)$ 
 such that $DF\in L^2(\Omega^X\times X,P\otimes \nu)$. 
 Recall also that the Ornstein-Uhlenbeck semi-group 
 $(P_t)_{t\in \real_+}$ is defined via 
$$
P_t I_n(f_n) = e^{-nt} I_n(f_n), \quad f_n\in 
 L^2 (X ,\nu )^{\circ n}, 
 \quad n\in \inte. 
$$ 
 In the sequel we also use the integral representation 
 of the Ornstein-Uhlenbeck semi-group $(P_t)_{t\in\real_+}$ 
 in terms of a probability kernel $p_t(\omega ,d\tilde{\omega}, d\hat{\omega} )$, 
 cf. e.g. \cite{surgailis}: 
\begin{equation} 
\label{eq:rep-OU}
 P_t F(\omega) = 
 \int_{\Omega^X \times \Omega^X } F(\tilde\om \cup \hat{\omega} )
 p_t(\omega ,d\tilde{\omega},d\hat{\omega} ) 
. 
\end{equation} 
 When $X=\real^n$, $\vert \cdot \vert_p$ denotes the 
 $\ell^p$-norm on $\real^n$, $p\geq 1$. 
 Assuming that 
$$\int_{\real^n} 1 \wedge \vert y \vert_2^2 \nu (dy) < \infty, 
$$  
 any $n$-dimensional infinitely divisible (ID) 
 random vector $F=(F_1,\dots, F_n)$ 
 without Gaussian component and with L\'evy measure $\nu$ can be represented as 
 the vector of single Poisson stochastic integrals 
\begin{equation} 
\label{eq:vaX2}
F=\left(\int_{ \{\vert y \vert_2 \leq 1\}} y_k\:  (\om(dy )-\nu (dy ))+\int_{\{\vert y \vert_2 > 1\}} y_k  \: \omega (dy )+b_k\right)_{1\leq k\leq n}
\end{equation} 
 for some $b\in \real^n$. 
 Indeed, the characteristic function of $F$ is given by 
$$ 
\varphi_F (u) = E[e^{i\langle F,u\rangle}] 
 =\exp\left( i\langle b,u \rangle 
 +\int_{\real^n} (e^{i\langle y ,u\rangle }-1-i\langle y ,u \rangle 
 \ind_{\{\vert y \vert_2 \leq 1\}})\nu (dy )\right) 
,  
$$ 
 $u\in\real^n$. 
\section{Deviation results from the Ornstein-Uhlenbeck semi-group} 
\label{sec:deviation-OU}
 As in \cite{hp}, 
 we need the following covariance identity on Poisson space, 
 which is obtained from the Ornstein-Uhlenbeck semi-group. 
\begin{lemma} 
\label{prop:id-cov}
 Let $F,G\in \Dom (D)$, then 
\begin{equation}
\label{eq:id-cov}
\Cov (F,G) = E\left[ 
 \int_0^\infty e^{-s}\int_X  D_y F P_s D_y G \nu(dy)ds 
 \right]. 
\end{equation} 
\end{lemma} 
\begin{Proof} 
 By orthogonality of multiple integrals of different 
 orders and continuity of $P_s$, $s\in \real_+$, on $L^2(\Omega^X,P)$, 
 it suffices to prove the identity for $F=I_n(f_n)$ and 
 $G=I_n(g_n)$: 
\begin{eqnarray*} 
E[I_n(f_n)I_n(g_n)] 
& = & n!\langle f_n, g_n\rangle_{L^2 (X, \nu )^{\circ n}}=n!\int_{X^n}f_n g_n \: d\nu^{\otimes n}\\
&=&n!\int_X\int_{X^{(n-1)}} f_n( x ,y) g_n( x , y) \: \nu^{\otimes (n-1)} (dx) \: \nu(dy)\\
&=&n \int_X E[I_{n-1} ( f_n(\cdot,y) ) I_{n-1} ( g_n(\cdot, y) ) ]\: \nu(dy)\\
&=&\frac 1 n E\left[ 
 \int_X D_yI_n(f_n) D_yI_n(g_n)\: \nu(dy) 
 \right] 
\\ 
&=& E \left[ 
 \int_0^\infty 
 e^{-ns} 
 \int_X D_yI_n(f_n) D_yI_n(g_n) \nu(dy) 
 ds 
 \right] 
\\ 
& = & E\left[ 
 \int_0^\infty e^{-s} \int_X D_y I_n(f_n)P_s 
 D_yI_n(g_n) \nu(dy) ds \right]. 
\end{eqnarray*} 
\end{Proof}
 Using the covariance identity \eqref{eq:id-cov}
 and the representation \eqref{eq:rep-OU} 
 we first state a general deviation result which slightly 
 improves the one presented in \cite{hp}. 
 In particular it will be 
 applied, in Section~\ref{sec:appli-vect}, 
 to obtain deviation inequalities on product spaces for 
 vectors of random functionals. 
 In this proposition and the following ones, the supremum on 
 $\Omega^X$ can be taken as an essential supremum with respect 
 to $P$. 
\begin{prop} 
\label{thd3} 
 Let $F \in \Dom (D)$ 
 be such that $e^{sF}\in \Dom (D)$, $0\leq s \leq t_0$, for some $t_0>0$. 
 Then 
$$ 
P(F-E[F]\geq x) \leq 
 \exp \left( 
\min_{0<t<t_0}\int_0^t h(s)\: ds -tx 
\right), 
 \qquad x>0, 
$$ 
 where 
\begin{equation} 
\label{eq:h1}
 h(s) = \sup_{(\om,\om')\in\Omega^X\times \Omega^X }\left| \int_{X}(e^{s D_yF(\om)}-1)\: D_{y}F(\om') \nu (dy)\right|, \quad s \in [0,t_0 ) . 
\end{equation}
 If moreover $h$ is nondecreasing and finite on $[0,t_0 )$ then 
\begin{equation} 
\label{*dev} 
P(F-E[F]\geq x) \leq 
 \exp \left( 
 -\int_0^x h^{-1}(s) ds\right), 
 \qquad 0<x<h(t_0^-), 
\end{equation} 
 where $h^{-1}$ is the left-continuous inverse of $h$: 
$$h^{-1} (x) = \inf \{ t >0 \ : \ h(t) \geq x\}, \qquad 0<x<h(t_0^-)
. 
$$ 
\end{prop} 
\begin{Proof} 
We start by deriving the following inequality 
for a centered random variable $F$: 
\begin{equation}
\label{eq:var1}
E[Fe^{sF}] 
 \leq 
 h(s) 
 E[e^{sF}] 
, \quad 0\leq s \leq t_0. 
\end{equation} 
 This follows from \eqref{eq:id-cov}. 
 Indeed, using the integral representation \eqref{eq:rep-OU} 
 of the Ornstein-Uhlenbeck semi-group $(P_t)_t$ for $P_v D_{y} F (\omega )$, we have 
\begin{eqnarray*}
\lefteqn{ 
 E[Fe^{sF}]
 = E\left[ \int_0^\infty e^{-v} 
 \int_X   D_{y}e^{sF} P_v D_{y}F \nu (dy) dv 
 \right] 
} 
\\ 
 & = &
\int_{\Omega^X} 
 \int_0^\infty e^{-v} 
 \int_X (e^{sD_{y}F (\omega ) }-1) e^{sF(\omega )} \int_{\Omega^X\times \Omega^X} D_{y}F( 
 \tilde{\omega} \cup \hat{\omega} 
 ) p_v (\om,d 
 \tilde{\omega}, d\hat{\omega} ) \nu (dy) dv 
 P( d\omega ) 
\\  
& \leq &
\int_{\Omega^X}  
 \int_0^\infty e^{-v} e^{sF(\omega )}
 \int_{\Omega^X\times \Omega^X} \left|\int_X (e^{sD_{y}F(\omega )}-1)  D_{y}F(\tilde{\omega}\cup \hat{\omega} ) \nu (dy) \right| p_v(\om,d\tilde{\omega},d\hat{\omega} ) dv 
 P( d \omega ) 
\\ 
& \leq&
\sup_{(\om,\om')\in\Omega^X\times \Omega^X}\left|\int_X (e^{sD_{y}F(\omega )}-1)  D_{y}F(\om') \nu (dy) \right|E\left[ e^{sF}
 \int_0^\infty e^{-v} dv 
 \right] 
\\ 
& = &
\sup_{(\om,\om')\in\Omega^X\times \Omega^X}\left|\int_X (e^{sD_{y}F(\omega )}-1)  D_{y}F(\om') \nu (dy) \right|E\left[ e^{sF}\right]
\end{eqnarray*} 
 which yields \eqref{eq:var1}. 
 In the general case, 
 we let $L(s)=E\left[e^{s(F-E[F])}\right]$ and obtain: 
$$
\frac{L'(s)}{L(s)} \leq h(s), 
 \quad 0 \leq s \leq t_0 
, 
$$ 
 which using Chebychev's inequality gives: 
\begin{equation}
\label{eq:minh}
P(F-E [F]\geq x)\leq \exp\left( 
-tx+\int_0^t h(s) ds\right) .
\end{equation} 
 Using the relation $\frac d{dt} \left(\int_0^t h(s)\: ds -tx\right)=h(t)-x$, 
 we can then optimize as follows:
\begin{eqnarray} 
\nonumber 
\min_{0<t<t_0}\int_0^t h(s)\: ds -tx& =& \int_0^{h^{-1}(x)} h(s)\: ds -xh^{-1}(x) 
\\
\nonumber 
&=&\int_0^x s\: dh^{-1}(s) -xh^{-1}(x) 
\\ 
\label{gth} 
&=& 
-\int_0^x h^{-1}(s)\: ds, 
\end{eqnarray}
 hence 
$$
P(F-E[F]\geq x) \leq \exp \left( 
 -\int_0^x h^{-1}(s) ds\right), 
 \ \ \ 0<x<h(t_0^-). 
$$ 
\end{Proof} 
 In the sequel we derive several corollaries from Proposition~\ref{thd3} 
 and discuss possible choices for the function $h$, 
 in particular for vectors of random functionals. 
 Note that since 
$$
h(t) \leq \int_X\|D_yF\|_\infty \left\|e^{t|D_yF|}-1\right\|_{\infty}\nu(dy), 
$$
 Proposition~\ref{thd3} recovers Proposition~3.3 in \cite{hp}, which 
 is obtained via a covariance identity relying on the Clark formula. 
 In the next proposition and following \cite{HRB}, 
 we obtain a better result by applying 
 Proposition~\ref{thd3} with more careful bounds. 
\begin{prop} 
\label{thd2.0} 
 Let $F : \Omega^X \to \real$ and let $K:X\to \real_+$ be a non-negative function such that 
\begin{equation} 
\label{j1} 
 D_y F (\omega ) 
 \leq K(y), 
 \qquad 
 y\in X,  \omega \in \Omega^X 
. 
\end{equation} 
 Then 
$$ 
P(F-E[F]\geq x) \leq 
 \exp \left( 
\min_{t >0}\int_0^t h(s)\: ds -tx 
\right), 
 \qquad x>0, 
$$ 
 where 
\begin {equation} 
\label{eq:h2.0}
 h(t) = \sup_{\om\in\Omega^X} \int_X \frac{e^{t K(y)}-1}{K(y)} 
 \vert D_y F(\om) \vert^2 \nu (dy), \quad t >0. 
\end{equation} 
 If moreover $h$ is finite on $[0,t_0)$ then 
\begin{equation} 
\label{j2} 
P(F-E[F]\geq x) \leq 
 \exp \left( 
 -\int_0^x h^{-1}(s) ds\right), 
 \qquad 0<x<h(t_0^-). 
\end{equation} 
 If $K(y)=0$, $y\in X$, we have: 
$$ 
P(F-E[F]\geq x) \leq \exp \left( 
 -\frac{x^2}{2\tilde{\alpha}^2} 
 \right), 
 \ \ \ x>0, 
$$ 
 with 
$$\tilde{\alpha}^2 = 
 \sup_{\om\in\Omega^X} \int_X \vert D_y F(\om) \vert^2 \nu (dy). 
$$ 
\end{prop} 
\begin{Proof}    
 Since when $K$ is $\real_+$-valued 
 the condition $D_yF_n(\omega ) \leq K(y)$, $\omega \in \Omega^X$, $y\in X$, is 
 satisfied we may apply Proposition~\ref{thd3} to $F_n=\max(-n,\min(F,n))$, $n\geq 1$, 
 and get 
\begin{eqnarray*}
h(t)&=&\sup_{(\om,\om')\in\Omega^X\times \Omega^X }\left|
 \int_X \frac{e^{tD_y F_n (\om)}-1}{D_y F_n (\om)}D_y F_n (\om)D_y F_n (\om')\: \nu(dy)\right|\\ 
&\leq&\sup_{(\om,\om')\in\Omega^X\times \Omega^X} 
 \int_X \frac{e^{t K(y)}-1}{K(y)}|D_y F_n (\om)|\:|D_y F_n (\om')|\: \nu(dy)\\
&\leq &\frac 12\sup_{(\om,\om')\in\Omega^X\times \Omega^X}\int_X \frac{e^{tK(y)}-1}{K(y)} 
 ( \vert D_y F_n (\om) \vert^2+ \vert D_y F_n (\om')\vert^2 )\: \nu (dy)\\
&\leq &\sup_{\om\in\Omega^X}\int_X \frac{e^{tK(y)}-1}{K(y)} 
 \vert D_y F_n (\om) \vert^2\: \nu (dy) 
\\ 
& \leq & \sup_{\om\in\Omega^X}\int_X \frac{e^{tK(y)}-1}{K(y)} 
 \vert D_y F (\om) \vert^2\: \nu (dy), 
\end{eqnarray*} 
 which allows to conclude. 
\end{Proof} 
 Note that if $K:X\to \real$ in \eqref{j1} is not necessarily positive 
 and $F, e^{sF}\in \Dom (D)$, $0\leq s \leq t_0$, for some $t_0>0$, 
 then applying Proposition~\ref{thd3} and the above argument directly to $F$ 
 yields: 
$$ 
P(F-E[F]\geq x) \leq 
 \exp \left( 
\min_{0<t<t_0}\int_0^t h(s)\: ds -tx 
\right), 
 \qquad x>0, 
$$ 
 and \eqref{j2} also holds provided $h$ is finite on $[0,t_0)$. 
\bigskip 

\noindent Part of the next corollary recovers a result of \cite{wuls2} (see also \cite{hp}). 
 This result is used in Corollary~\ref{corol:NP1} below 
 as well as in the infinite variance case in Section~\ref{sec:no-variance}. 
\begin{corollary} 
\label{prop:dev-NP}
 Let $F\in L^2(\Omega^X,P)$ be such that $DF\leq K$, $P\otimes \nu$-a.e., 
 for some $K\in \real$, and $\|DF\|_{L^{\infty}(\Omega^X, L^2(X , \nu ))}\leq \tilde{\alpha}$. 
 We have for $K>0$: 
\begin{equation} 
\label{eq:Dev-Nico}
 P(F-E[F]\geq x)\leq e^{x/K} 
 \left( 
 1+\frac{xK}{\tilde{\alpha}^2} 
 \right)^{ 
 - 
 \frac{x}{K} 
 - 
 \frac{\tilde{\alpha}^2}{K^2} 
 }, \qquad x > 0, 
\end{equation} 
 and for $K=0$: 
\begin{equation} 
\label{eq:Dev-Nico1}
 P(F-E[F]\geq x)\leq\exp\left(-\frac{x^2}{ 
 2 \tilde{\alpha}^2 }\right), 
 \qquad x>0 
. 
\end{equation} 
\end{corollary} 
\begin{Proof} 
 If $K\geq 0$, let us first assume that $F$ is a bounded random variable. 
 The function $h$ in \eqref{eq:h2.0} is such that 
$$ 
 h(t) \leq \frac{e^{t K}-1}{K} \|DF\|_{L^\infty(\Omega^X , L^2(X , \nu ))}^2 
 \leq \tilde{\alpha}^2\frac{e^{t K}-1}{K}, \qquad t>0. 
$$ 
 Applying \eqref{eq:minh} with $\tilde{\alpha}^2 ( e^{t K}-1 ) / K$ gives 
$$ 
P(F-E[F]\geq x) 
 \leq \exp\left( -tx+\frac{\tilde{\alpha}^2}{K^2}(e^{tK}-tK-1)\right).
$$ 
 Optimizing in $t$ with $t= K^{-1} \log ( 1 + K x / \tilde{\alpha}^2 )$ 
 (or using directly \eqref{*dev} with the inverse 
 $K^{-1} \log \left(1+K t / \tilde{\alpha}^2 \right)$) we have 
$$ 
P(F-E[F]\geq x)\leq \exp\left( 
 \frac xK-\left(\frac x K +\frac{\tilde{\alpha}^2}{K^2}\right)\log \left(1+\frac{xK}{\tilde{\alpha}^2}\right)\right),
$$ 
 which yields \eqref{eq:Dev-Nico0}, \eqref{eq:Dev-Nico1} and \eqref{eq:Dev-Nico}, 
 depending on the value of $K$. 
 For unbounded $F$, apply the above to $F_n=\max(-n,\min(F,n))$ 
 with $\vert DF_n \vert \leq \vert DF \vert$, $n\geq 1$. 
 Then \eqref{eq:Dev-Nico} follows since, 
 as  $n$ goes to infinity, $F_n$ converges to $F$ in $L^2(\Omega^X )$, 
 $DF_n$  converges to $DF$ in $L^2(\Omega^X , L^2(X,\nu ))$, and $DF_n\leq K$, $n\geq 1$. 
 The same argument applies if $K=0$. 
\end{Proof} 
 In case $K<0$ and $e^{tF}\in \Dom (D)$ for all $t>0$, Proposition~\ref{thd2.0} yields 
 in a similar way: 
\begin{equation} 
\label{eq:Dev-Nico0}
 P(F-E[F]\geq x)\leq e^{x/K} 
 \left( 
 1+\frac{xK}{\tilde{\alpha}^2} 
 \right)^{ 
 - 
 \frac{x}{K} 
 - 
 \frac{\tilde{\alpha}^2}{K^2} 
 }, \qquad 0 < x < -\frac{\tilde{\alpha}^2}{K}. 
\end{equation} 

\noindent If $F$ is an infinitely divisible 
 random variable in $\real^n$, without Gaussian component 
 and with L\'evy measure $\nu$, 
 the representation \eqref{eq:vaX2} shows that 
 for any Lipschitz$(c)$ function $f:\real^n \to \real$, 
\begin{eqnarray} 
\nonumber 
 \vert D_x 
 f(F) 
 (\omega ) \vert 
 & = & 
 \vert 
 f (F(\omega \cup \{x\})) 
 - 
 f (F(\omega ) ) 
 \vert 
\\ 
\nonumber 
& \leq & 
 c \Vert 
 F(\omega \cup \{x \})
 - 
 F(\omega )
 \Vert 
\\ 
\label{hjl} 
& = & 
 c \Vert 
 x \Vert, 
\end{eqnarray} 
 where $\Vert \cdot \Vert$ is any norm in $\real^n$. 
 Hence when $X=\real^n$ and $\nu$ has bounded support, 
 Corollary~\ref{prop:dev-NP} also recovers 
 Corollary~1 of \cite{h3} with 
$$ 
 K = \inf \{ r>0 \ : \ \nu (\{x \in X \ : \ 
 \Vert x \Vert > r \}) =0\}, 
$$ 
 and $\tilde{\alpha}^2 = \int_{\real^n} \Vert y \Vert^2 \nu (dy )$, 
 i.e. 
$$ 
P(f(F)-E[f(F)]\geq x)\leq e^{\frac{x}{cK}} 
 \left( 
 1+\frac{xK}{c \tilde{\alpha}^2} 
 \right)^{ 
 - 
 \frac{x}{cK} 
 - 
 \frac{\tilde{\alpha}^2}{K^2} 
 }, \quad x>0. 
$$ 
 On a product $X=\{ 1, \ldots, n\} \times Y$, where $Y$ is a  
 $\vert \cdot \vert_Y$-normed linear space, we have the identification 
$$\Omega^X \simeq \Omega^Y\times \cdots \times \Omega^Y, 
 \quad 
 \quad 
 \omega = (\omega_1,\ldots ,\omega_n) 
 \in \Omega^X, 
$$ 
 and 
$$D_{(i,y)} F (\omega) = 
 \sum_{j=1}^n 
 1_{\{i=j\}} 
 ( 
 F (\omega_1,\ldots ,\omega_{i-1} , \omega_i\cup \{y\},\omega_{i+1} ,\ldots ,\omega_n) 
 - 
 F (\omega_1,\ldots ,\omega_n) 
 ) 
, 
$$ 
 $i=1,\ldots, n$, $y\in Y$. 
 Proposition~\ref{thd3} can be directly applied with 
 $d\nu (i,y) = d\nu_i(y)$, $i=1,\ldots ,n$, $y\in Y$, and  
\begin{equation} 
\label{eq:h-direct1}
h(t) = \left\|\sum_{i=1}^n \int_Y 
 \|D_{(i,y)}F\|_\infty(e^{t|D_{(i,y)}F|}-1)\nu_i(dy)\right\|_\infty , 
\end{equation}
 or if $\nu = \nu_1 = \cdots = \nu_n$, with  
$$
h(t)=\tilde{\beta} \int_Y \vert y\vert_Y (e^{t\beta \vert y\vert_Y}-1)\nu(dy), 
$$
where 
$$
\tilde{\beta} =\sup_{\om, y\not=0} \sum_{i=1}^n \frac{|D_{(i,y)}F(\om)|}{\vert y\vert_Y},\med \mbox{ and } \med \beta=\sup_{i,\om,y\not=0}\frac{|D_{(i,y)}F(\om)|}{\vert y\vert_Y}.
$$ 
 In fact, a stronger result can be obtained as a corollary of Proposition~\ref{thd2.0}. 
\begin{corollary} 
\label{thd2} 
 Let $X=\{1,\ldots ,n\}\times Y$, where $Y$ is a 
 $\vert \cdot \vert_Y$-normed linear space 
 and $d\nu (i,y) = d\nu_i(y)$, $i=1,\ldots ,n$, $y\in Y$. 
 Let $F:\Omega^X \to \real$ 
 and let $\beta_i \geq 0$, $i=1,\ldots ,n$, be such that 
$$ 
 D_{(i,y)}F (\omega ) 
 \leq \beta_i \vert y \vert_Y, 
 \quad 
 i=1,\ldots , n , 
 \quad 
 y\in Y, 
 \omega \in \Omega^X 
. 
$$ 
 Then 
$$ 
P(F-E[F]\geq x) \leq 
 \exp \left( 
\min_{t>0}\int_0^t h(s)\: ds -tx 
\right), 
 \qquad x>0, 
$$ 
 where 
\begin {equation}
\label{eq:h2}
 h(t) = \sup_{\om\in\Omega^X} \sum_{i=1}^n \int_{Y}\frac{e^{t \beta_i \vert y \vert_Y }-1}{\beta_i \vert y\vert_Y}(D_{(i,y)}F(\om))^2 \nu_i (dy), \quad t >0. 
\end{equation}
 If moreover $h$ is finite on $[0,t_0)$ then 
\begin{equation*}
P(F-E[F]\geq x) \leq 
 \exp \left( 
 -\int_0^x h^{-1}(s) ds\right), 
 \qquad 0<x<h(t_0^-) 
. 
\end{equation*} 
 If $\beta_i=0$, $i=1,\ldots ,n$, {\it i.e.} for 
 decreasing functionals, we have: 
$$ 
P(F-E[F]\geq x) \leq \exp \left( 
 -\frac{x^2}{2\tilde{\alpha}^2} 
 \right), 
 \ \ \ x>0, 
$$ 
 with 
$$\tilde{\alpha}^2 = 
 \sup_{\om\in\Omega^X} \sum_{i=1}^n \int_{Y} 
 (D_{(i,y)}F(\om))^2 \nu_i (dy). 
$$ 
\end{corollary} 
\begin{Proof} 
 Apply Proposition~\ref{thd2.0} with 
 $K(i,y) = \beta_i \vert y \vert_Y$, $1\leq i \leq n$, $y\in Y$. 
\end{Proof} 
 As a consequence of \eqref{eq:h2}, 
 and if $\nu : = \nu_1 = \cdots = \nu_n$, one can take: 
\begin{equation} 
\label{em1} 
 h(t) = \frac{\tilde{\alpha}^2}{\beta} 
 \int_Y 
 \vert y \vert_Y 
 (e^{t \beta \vert y \vert_Y}-1)
 \nu (dy), \quad t\in [0,t_0], 
\end{equation} 
 with
$$
\tilde{\alpha}^2 = \sup_{\omega\in \Omega^X , \ y \not=0} 
 \sum_{i=1}^n  \frac{\vert D_{(i,y)}F (\omega ) 
 \vert^2}{\vert y \vert_Y^2},
 \quad 
 \mbox{and} 
 \quad 
\beta = \sup_{i, \omega ,y \not=0} \frac{\vert D_{(i,y)}F (\omega ) 
 \vert}{\vert y \vert_Y} 
. 
$$ 
 Taking 
\begin{equation} 
\label{also} 
h(t)  = \sum_{i=1}^n \beta_i \int_Y \vert y \vert_Y (e^{t\beta_i \vert y \vert_Y} - 1) 
 \nu_i (dy), \qquad t\in [0,t_0], 
\end{equation} 
 allows to recover the bound implied by \eqref{eq:h-direct1} in this case. 

\noindent For example, if $n=1$ and 
$$
F_1 (\omega ) = \int_Y u_1(y ) (\omega (dy ) - \nu (dy )) , 
 \ldots 
 , 
  F_m (\omega ) = \int_Y u_m(y ) (\omega (dy ) - \nu (dy )) 
$$ 
 are $m$ (not necessarily independent) single Poisson stochastic 
 integrals and $F=g(F_1,\ldots ,F_m)$, we have 
$$ 
\beta \leq \sup_{x_1,\ldots ,x_m ,y \not=0} 
 \frac{\vert g(x_1+u_1(y),\ldots , x_m + u_m (y)) 
 - g(x_1,\ldots ,x_m )\vert}{\vert y \vert_Y}.
$$
\bigskip 

\noindent 
 The following statement is obtained from 
 Corollary~\ref{prop:dev-NP} on a product space, 
 in the same way as Corollary~\ref{thd2} is obtained 
 from Proposition~\ref{thd2.0}. 
\begin{corollary}
\label{corol:NP1} 
 Let $X=\{1,\ldots ,n\}\times Y$ with 
 $d\nu (i,y) = d\nu_i(y)$, $i=1,\ldots ,n$, $y\in Y$. 
 Let $F$ be such that 
 $ D_{(i,y)}F(\om)\leq K$, $P\otimes \nu_i$-a.e., 
 $i=1,\ldots , n $, for some $K\in \real$ and 
$$\left\| 
 \sum_{i=1}^n \|D_{(i,\cdot )}F \|_{L^2(Y;\nu_i)}^2 \right\|_{L^\infty (\Omega^X,P)} 
 \leq \tilde{\alpha}^2.$$ 
 We have for $K>0$: 
\begin{equation} 
\label{eq:Dev-Nico-bis} 
 P(F-E[F]\geq x)\leq e^{x/K} 
 \left( 
 1+\frac{xK}{\tilde{\alpha}^2} 
 \right)^{ 
 - 
 \frac{x}{K} 
 - 
 \frac{\tilde{\alpha}^2}{K^2} 
 }, \qquad x > 0, 
\end{equation} 
 and for $K=0$: 
\begin{equation} 
\label{eq:Dev-Nico-bis1}
 P(F-E[F]\geq x)\leq\exp\left(-\frac{x^2}{ 
 2 \tilde{\alpha}^2 }\right), 
 \qquad x>0 
. 
\end{equation} 
 Moreover if $K<0$ and 
 $e^{tF}\in \Dom (D)$ for all $t>0$, then 
\begin{equation} 
\label{eq:Dev-Nico-bis0} 
 P(F-E[F]\geq x)\leq e^{x/K} 
 \left( 
 1+\frac{xK}{\tilde{\alpha}^2} 
 \right)^{ 
 - 
 \frac{x}{K} 
 - 
 \frac{\tilde{\alpha}^2}{K^2} 
 }, \qquad 0 < x < -\frac{\tilde{\alpha}^2}{K}. 
\end{equation} 
\end{corollary}
\section{Application to random vectors}
\label{sec:appli-vect}
 We start by applying Corollary~\ref{thd2} and Corollary~\ref{corol:NP1} 
 to random vectors $(F_1,\ldots ,F_n)$ on the 
 product space $\Omega^X \simeq \Omega^Y\times \cdots \times \Omega^Y$ 
 where $X=\{ 1, \ldots, n\} \times Y$ and $Y$ is a  
 $\vert \cdot \vert_Y$-normed linear space. 
 Corollary~\ref{thd2} yields 
\begin{eqnarray*} 
 P(g(F_1,\ldots , F_n) -E[g(F_1,\ldots , F_n)] 
 \geq x) 
 & \leq & 
\exp \left( 
 -\int_0^x h^{-1}(s)\: ds\right), 
\end{eqnarray*} 
 $0<x<h(t_0^-)$,  where $g:\real^n \to \real$, provided the function 
$$ 
 h(t) = 
 \left\| 
 \sum_{i=1}^n 
 \int_Y 
\Big( D_{(i,y)} g(F_1(\om),\ldots , F_n(\om))\Big)^2 \frac{e^{t \beta_i\vert y\vert_Y}-1}{\beta_i \vert y \vert_Y}\:
 \nu_i (dy) 
 \right\|_{L^\infty (\Omega^X , P)} 
, \quad t\in [0,t_0]
, 
$$ 
 is finite on $[0,t_0)$, with $\beta_i$ as in Corollary~\ref{thd2}. 
 Several particular cases are now presented. 
\subsection*{Random vectors with independent components} 
 If $F_1,\ldots ,F_n$ are $n$ independent random variables 
 defined on 
 $\Omega^X = \Omega^Y\times \cdots \times \Omega^Y$ 
 with $F_i = F_i(\omega_i)$, 
 $i=1,\ldots ,n$, and $g:\real^n \to \real$, 
 an $\ell^1$-Lipschitz$(c)$ function, we have 
\begin{eqnarray*} 
\lefteqn{ 
 \vert D_{(i,y)} 
 g(F_1,\ldots , F_n) 
 (\omega ) \vert 
} 
\\ 
 & = & 
 \vert 
 g(F_1(\omega_1 ) ,\ldots , F_i(\omega_i \cup \{y\}), \ldots , 
 F_n(\omega_n ) 
 ) 
 - 
 g(F_1(\omega_1 ) ,\ldots , F_n(\omega_n ) ) 
 \vert 
\\ 
& \leq & 
 c \vert 
 F_i(\omega_i \cup \{y\})
 - 
 F_i(\omega_i )
 \vert 
\\ 
& \leq & 
 c \vert 
 D_y F_i(\omega )
 \vert. 
\end{eqnarray*} 
 Now we can take in Corollary~\ref{thd2}: 
$$ 
 h(t) \leq c^2 
 \sup_{\om\in\Omega^X} 
 \sum_{i=1}^n 
 \int_{Y} 
 \frac{e^{ct \beta_i \vert y \vert_Y}-1}{c\beta_i \vert y\vert_Y } 
 (D_yF_i(\om))^2 \nu_i (dy), \quad t\in [0,t_0]
$$
 with 
$$\beta_i=\sup_{y\in Y, \ \omega_i \in\Omega^Y} \frac{|D_yF_i (\omega_i )|}{\vert y \vert_Y}. 
$$ 
 Moreover when $\nu = \nu_1=\cdots = \nu_n$, 
 we can take in \eqref{em1}:
$$
\tilde{\alpha}^2 =  
 c^2 
 \sup_{\omega ,y \not=0} 
 \sum_{i=1}^n \frac{\vert D_{y}F_i (\omega ) 
 \vert^2}{\vert y \vert_Y^2}, 
 \quad 
 \mbox{and} 
 \quad 
 \beta = 
 c 
 \sup_{i, \omega ,y \not=0} \frac{\vert D_yF_i (\omega ) 
 \vert}{\vert y \vert_Y}. 
$$ 
\subsection*{Independent vectors of Poisson stochastic integrals} 
 Assume that $Y$ is a normed linear space and that 
 $\int_{Y} 1\wedge \vert y \vert_Y^2 \nu_i(dy) < \infty$, 
 $i=1,\ldots , n$. 
 If $G=g(F_1,\dots, F_n)$ where $g:\real^n\to \real$ and $F_1,\dots, F_n$ 
 are independent Poisson stochastic integrals of the form \eqref{eq:vaX2}: 
$$ 
F_i(\om_i) = 
\int_{ \{\vert y\vert_Y \leq 1\}} y \:  (\om_i (dy)-\nu_i (dy))+\int_{\{\vert y\vert_Y > 1\}} y 
 \om_i (dy)+b_i 
, 
\qquad 1\leq i \leq n, 
$$ 
 we have 
\begin{equation}
\label{eq:D_(i,y)g} 
 D_{(i,y)} g(F_1,\ldots , F_n) 
 = g(F_1,\ldots , F_i+y,\ldots , F_n) - g(F_1,\ldots , F_n). 
\end{equation} 
 From Corollary~\ref{thd2} we have, denoting by 
 $(e_1,\ldots , e_n)$ the canonical basis on $\real^n$: 
\begin{eqnarray} 
\nonumber 
\lefteqn{ 
h(t)=\sup_{\om\in\Omega^X} \sum_{i=1}^n \int_{Y}\frac{e^{t \beta_i \vert y \vert_Y}-1}{\beta_i \vert y\vert_Y }
} 
\\
\nonumber  
&& (g(F_1(\omega_1),\dots, F_i(\omega_i)+y,\dots, F_n(\omega_i))-g(F_1(\omega_1),\dots, F_i(\omega_i),\dots, F_n(\omega_n))^2\nu_i (dy)\\ 
\label{rf} 
&\leq & \sup_{x\in\real^n} \sum_{i=1}^n \int_{Y}\frac{e^{t \beta_i \vert y \vert_Y}-1}{\beta_i \vert y\vert_Y} (g(x+ye_i)-g(x))^2\nu_i (dy), 
\end{eqnarray} 
 which recovers Theorem 1 of \cite{HRB} and (3) therein as a particular case. 
 We may also take 
\begin{equation}
\label{eq:alpha-beta}
\tilde{\alpha}^2 = \sup_{x ,y \not=0} 
 \sum_{i=1}^n \frac{
 \vert g( x + y e_i ) - g ( x ) \vert^2}{\vert y \vert_Y^2},
\quad \mbox{and} \quad 
\beta = \sup_{i, x ,y \not=0} \frac{
 \vert g(x + y e_i ) - g (x ) \vert 
 }{\vert y \vert_Y}, 
\end{equation}
 in \eqref{em1}. 
 If $g:\real^n\to \real$ is $\ell^1$-Lipschitz$(c)$, then 
 $\beta = c$, and so \eqref{also} gives: 
\begin{equation} 
\label{eqr} 
 h(t) = c \sum_{i=1}^n 
 \int_Y 
 \vert y \vert_Y 
 (e^{t c \vert y \vert_Y }-1)
 \nu_i (dy), \quad t\in [0,t_0]. 
\end{equation} 

\noindent 
 For $g(x)=\sup(x_1,\dots, x_n)$, with $Y=\rit$, $\beta=1$ and 
\begin{equation} 
\label{hj} 
 D_{(i,y)}g(F_1,\dots, F_n) 
 = 
\left\{ 
\begin{array}{ll} 
 0, & y\leq \sup(F_1,\dots, F_n) -F_i, 
\\ 
F_i+y-\sup(F_1,\dots, F_n), 
 & y > \sup(F_1,\dots, F_n) -F_i
, 
\end{array} 
\right. 
\end{equation} 
 $i=1,\ldots ,n$, $y\in\real$. 
 Hence \eqref{eq:h2} leads to 
\begin{eqnarray} 
\nonumber 
\lefteqn{ 
h(t)=\sup_{\om\in\Omega^X} \sum_{i=1}^n \int_{\rit}\frac{e^{t\vert y \vert}-1}{\vert y\vert}
} 
\\
\nonumber  
&& (g(F_1(\omega_1),\dots, F_i(\omega_i)+y,\dots, F_n(\omega_i))-g(F_1(\omega_1),\dots, F_i(\omega_i),\dots, F_n(\omega_n))^2\nu_i (dy)\\  
\nonumber
&=& \sup_{\om\in\Omega^X} \sum_{i=1}^n \int_{\sup(F_1,\ldots, F_n)-F_i}^\infty\frac{e^{t \vert y \vert}-1}{\vert y\vert} (F_i+y-\sup(F_1,\dots, F_n))^2\nu_i (dy),\\
\label{eq:h-sup}
&\leq& \sum_{i=1}^n \int_0^\infty y(e^{t y}-1)\nu_i (dy).
\end{eqnarray} 
 Note that in \eqref{eq:alpha-beta} the constants $\tilde{\alpha}^2$ and $\beta$ can be computed 
 in terms of the Lipschitz constant of $g$ with respect to the $\ell^1$-norm. 
 This however does not lead to dimension free estimates. 
 Next, we show, using \eqref{rf}, that dimension free 
 estimates can be obtained when $g$ is the Euclidean norm on $\real^n$. 
 The other results of \cite{HRB} can similarly be generalized to the present 
 framework. 
\subsection*{Dimension free inequalities for random vectors} 
 Dimension free inequalities for $\ell^2$-Lipschitz functions 
 of independent infinitely divisible random vectors 
 with finite exponential moments 
 have been obtained in Corollary~4 of \cite{HRB}.
 In the next proposition we extend this result to Poisson 
 random functionals. 
\begin{prop} 
\label{df1} 
 Let $f:\real^n\to\real$ be $\ell^2$-Lipschitz$(c)$, and let 
 $F=(F_1,\ldots ,F_n)$ be a vector of independent random functionals. 
 Let 
$$
 \beta_i = \sup_{y\in Y, ~ \omega\in \Omega^X} \frac{\vert D_y F_i(\omega )\vert}{\vert y\vert_Y}, 
 \quad  i=1,\ldots , n 
, 
$$ 
 and assume that 
\begin{eqnarray} 
\label{rgh} 
 h(t) 
& = & 
 8 \max_{i=1,\ldots ,n} 
 \beta_i 
 \int_{Y} 
 \vert y \vert_Y 
 (e^{t \beta_i \vert y \vert_Y}-1) 
 \nu_i (dy)
\\ 
\nonumber 
& & + 
 \frac{2n}{(E[\vert F-E[F]\vert_2])^2} 
 \max_{i=1,\ldots ,n} 
 \beta_i^3 
 \int_{Y} 
 \vert y \vert_Y^3 
 (e^{t \beta_i \vert y \vert_Y}-1) 
 \nu_i (dy)
\end{eqnarray}
 is finite in $t\in [0,t_0)$. 
 Then 
\begin{equation} 
\label{et1} 
P\left( f (F_1,\ldots , F_n) \geq E[ f (F_1,\ldots , F_n) ] 
 + c \sqrt{2\sum_{i=1}^n \Var F_i}+ c x \right) 
 \leq \exp \left( 
 -\int_0^x h^{-1}(s) ds\right), 
\end{equation}
 $0<x<h(t_0^-)$. 
\end{prop} 
\begin{Proof} 
 Define $\phi :\real^n\to\real$ by 
$$\phi (x ) = \sqrt{E[\vert x-G \vert_2^2]}, 
$$ 
 where $\vert x\vert_2$ is the Euclidean norm of $x\in \real^n$ 
 and $G$ is an independent copy of $F$. 
 As in the proof of Corollary~4 in \cite{HRB}, we have 
$$\vert \phi(x + u e_i ) - \phi (x ) \vert^2 
 \leq 
 \frac{8u^2E[(x_i-G_i)^2]}{E[\vert x-G\vert^2]} 
 + \frac{2u^4}{\sum_{k=1}^n \Var G_k}
, 
 \quad x\in \real^n, 
 \quad u\in \real. 
$$ 
 Hence for $\phi(F)$, Corollary~\ref{thd2} applies with 
\begin{eqnarray*} 
 h_\phi (t) & = &  
 \sup_{\om\in\Omega^X} \sum_{i=1}^n \int_{Y}\frac{e^{t \beta_i \vert y \vert_Y}-1}{\beta_i \vert y\vert_Y}(D_{(i,y)}\phi (F(\om)))^2 \nu_i (dy) 
\\ 
& \leq & 
 \sup_{\om\in\Omega^X} \sum_{i=1}^n \int_{Y}\frac{e^{t \beta_i \vert y \vert_Y}-1}{\beta_i \vert y\vert_Y} 
\\ 
& & (\phi (F_1,\ldots ,F_i+D_yF_i,\ldots ,F_n ) (\omega ) 
 - 
 \phi (F_1,\ldots , F_n)(\omega ) )^2 \nu_i (dy) 
\\ 
& \leq & 
 \sup_{\om\in\Omega^X} 
 \sum_{i=1}^n 
 \int_{Y}\frac{e^{t \beta_i \vert y \vert_Y}-1}{\beta_i \vert y\vert_Y} 
 \left( 
 8 \vert D_yF_i(\omega ) \vert^2 
 \frac{E_G[(F_i(\omega ) -G_i)^2]}{E_G[\vert F(\omega ) -G\vert^2]} 
 + \frac{2\vert D_yF_i(\omega ) \vert^4}{\sum_{k=1}^n \Var G_k } 
 \right) 
 \nu_i (dy)
\\ 
& \leq & 
 8 \sup_{\om\in\Omega^X} 
 \max_{i=1,\ldots ,n} 
 \int_{Y}\frac{e^{t \beta_i \vert y \vert_Y}-1}{\beta_i \vert y\vert_Y} 
 \vert D_yF_i(\omega ) \vert^2 
 \nu_i (dy)
\\ 
& & + 
 \frac{2}{\sum_{k=1}^n \Var G_k} 
 \sup_{\om\in\Omega^X} 
 \sum_{i=1}^n 
 \int_{Y}\frac{e^{t \beta_i \vert y \vert_Y}-1}{\beta_i \vert y\vert_Y} 
 \vert D_yF_i(\omega ) \vert^4 
 \nu_i (dy)
\\ 
& \leq & 
 8 \max_{i=1,\ldots ,n} 
 \beta_i 
 \int_{Y} 
 \vert y \vert_Y 
 (e^{t \beta_i \vert y \vert_Y}-1) 
 \nu_i (dy)
\\ 
& & + 
 \frac{2}{\sum_{k=1}^n \Var G_k} 
 \sum_{i=1}^n 
 \beta_i^3 
 \int_{Y} 
 \vert y \vert_Y^3 
 (e^{t \beta_i \vert y \vert_Y}-1) 
 \nu_i (dy) 
, \qquad t\in [0,t_0]. 
\end{eqnarray*} 
 Finally, Corollary~\ref{thd2} with the bounds 
$$E[\phi (F)] \leq \sqrt{2\sum_{k=1}^n \Var G_k} 
, 
$$ 
 and $\vert f(x)-E[f(F)]\vert \leq c \phi (x)$, 
 yields \eqref{et1}. 
\end{Proof} 
 The function $h$ in \eqref{rgh} is bounded independently of the 
 dimension $n$ if $(F_1,\ldots ,F_n)$ are i.i.d., since 
$$
 n\min_{1\leq k \leq n} (E[\vert F_k \vert ])^2 
 \leq (E[\vert F \vert_2 ])^2 
 \leq 
 n \max_{1\leq k \leq n} E[\vert F_k \vert^2] 
 . 
$$ 
 For the Euclidean 
 norm of independent infinitely divisible random vectors 
 with finite exponential moments, better results 
 have been obtained in Corollary~3 of \cite{HRB}.
 In the next proposition we extend this result to Poisson 
 random functionals. 
\begin{prop} 
\label{df} 
 Let $F=(F_1,\ldots ,F_n)$ be a vector of independent random functionals, 
 and let 
$$
 \beta_i = \sup_{y\in Y, ~ \omega\in \Omega^X} \frac{\vert D_y F_i(\omega )\vert}{\vert y\vert_Y}, 
 \quad  i=1,\ldots , n 
, 
$$ 
 and assume that 
\begin{eqnarray*} 
 h(t) 
& = & 
 8 \max_{i=1,\ldots ,n} 
 \beta_i 
 \int_{Y} 
 \vert y \vert_Y 
 (e^{t \beta_i \vert y \vert_Y}-1) 
 \nu_i (dy)
\\ 
& & + 
 \frac{2n}{(E[\vert F\vert_2])^2} 
 \max_{i=1,\ldots ,n} 
 \beta_i^3 
 \int_{Y} 
 \vert y \vert_Y^3 
 (e^{t \beta_i \vert y \vert_Y}-1) 
 \nu_i (dy)
\end{eqnarray*}
 is finite in $t\in [0,t_0)$. 
 Then 
\begin{equation} 
\label{et} 
P(\vert (F_1,\ldots , F_n) \vert_2 \geq 2E[\vert (F_1,\ldots , F_n) \vert_2 ] +x) 
 \leq \exp \left( 
 -\int_0^x h^{-1}(s) ds\right), 
\end{equation}
 $0<x<h(t_0^-)$. 
\end{prop} 
\begin{Proof} 
 Let $f(x) = (\vert x \vert_2 - E[\vert F \vert_2])^+$, $x\in \real^n$. 
 From \cite{HRB} we have the inequality 
$$
\vert f(x+ue_i)-f(x)\vert^2 
 \leq 8 \vert u \vert^2 
 \frac{\vert x_i \vert^2}{\vert x \vert_2^2} 
 + \frac{2\vert u \vert^4}{(E[\vert F \vert_2])^2}, 
 \quad x\in \real^n, 
 \quad u\in \real. 
$$ 
 Hence for $f(F)$, repeating the bounds in the proof of Proposition~\ref{df1} 
 we get 
\begin{eqnarray*} 
 h(t) 
 & \leq & 
 8 \max_{i=1,\ldots ,n} 
 \beta_i 
 \int_{Y} 
 \vert y \vert_Y 
 (e^{t \beta_i \vert y \vert_Y}-1) 
 \nu_i (dy)
\\ 
& & + 
 \frac{2}{(E[\vert F\vert_2])^2} 
 \sum_{i=1}^n 
 \beta_i^3 
 \int_{Y} 
 \vert y \vert_Y^3 
 (e^{t \beta_i \vert y \vert_Y}-1) 
 \nu_i (dy)
, \quad t\in [0,t_0]. 
\end{eqnarray*} 
 Finally, using $\vert x \vert_2 - E[\vert F \vert_2] 
 \leq (\vert x \vert_2 - E[\vert F \vert_2] )^+$ 
 and $E[(\vert F \vert_2 - E[\vert F \vert_2])^+] 
 \leq E[\vert F \vert_2]$ gives \eqref{et} for $g(F)=\vert F\vert_2 $. 
\end{Proof} 
 Similarly to Proposition~\ref{df1}, the deviation result of \eqref{et} 
 is dimension free if $(F_1,\ldots ,F_n)$ are i.i.d. 

\noindent Next, we obtain a dimension free deviation for the 
 Euclidean norm of a vector of $n$ i.i.d. random functionals with bounded support.
 The non-identically distributed case is done similarly, 
 while for single integrals it is in \cite{HRB}. 
\begin{corollary} 
\label{corol:norm-dev-borné} 
 Let $\nu = \nu_1=\cdots = \nu_n$ have bounded support 
 in $B_Y(0,R)$, let $\beta = \beta_1=\cdots = \beta_n$, 
 and let $F=(F_1,\ldots ,F_n)$ be an i.i.d. vector. 
 Then, for all $x>0$, 
\begin{equation}
\label{eq:dev-norme-bornée}
P(\vert F \vert_2 \geq x + 2E[\vert F \vert_2 ] ) 
 \leq 
 \exp\left( 
 \frac x{\beta R}-\left(\frac x {\beta R} +\frac{\tilde{\alpha}_R^2}{\beta^2R^2}\right)\log \left(1+\frac{x\beta R}{\tilde{\alpha}_R^2}\right)\right), 
\end{equation}
 where 
$$ 
 \tilde{\alpha}_R^2 
 = 
 \left( 
 8\beta^2 
 + 
 \frac{2 \beta^5  
 R^2 
 }{(E[\vert F_1\vert])^2} 
 \right) 
 \int_{Y} 
 \vert y \vert_Y^2 
 \nu (dy)
. 
$$ 
\end{corollary} 
\begin{Proof} 
 Apply Proposition~\ref{df} with  
\begin{eqnarray*} 
 h(t) & \leq & 
 \left( 
 \frac{8 \beta}{R} 
 + 
 \frac{2n \beta^4 
 R 
}{(E[\vert F\vert_2])^2} 
 \right) 
 (e^{t \beta R }-1) 
 \int_{Y} 
 \vert y \vert_Y^2 
 \nu (dy)
 \leq  \tilde{\alpha}_R^2
 \frac{e^{t \beta R }-1}{\beta R}, 
\end{eqnarray*} 
 and compute explicitly the right hand side of \eqref{et}.
\end{Proof} 
 The following result yields an exponential integrability property, 
 independent of $n$ for the $\ell^2$-norm of infinitely divisible random vector whose 
 L\'evy measures have bounded supports. 
 The non identically distributed case is similar. 
 For independent infinitely divisible random variables 
 an analog result is obtained in \cite{HRB}. 
\begin{corollary} 
\label{corol:intégrabilité1}
 Let $F=(F_1,\dots, F_n)$ be as in Corollary~\ref{corol:norm-dev-borné} then 
 for all $\lambda$, with $ 0 < \lambda  < \beta^2 R^2/(e \tilde{\alpha}_R^2)$, we have: 
\begin{equation} 
\label{eq:int-exp1}
E \left[ 
 \exp\left( 
 \frac {\vert F\vert_2 }{\beta R} \log_+ \frac{\lambda\vert F\vert_2 } {\beta R} \right) 
 \right] 
 <\infty, 
\end{equation} 
 with $\log_+ x = \max ( \log x, 0)$, $x>0$. 
\end{corollary} 
\begin{Proof} Let $\lambda<\beta^2 R^2/(e\tilde{\alpha}_R^2)$. 
 We have, using \eqref{eq:dev-norme-bornée}: 
\begin{eqnarray*} 
\lefteqn{ 
E \left[ 
 \exp \left( 
 \frac  {\vert F\vert_2 }{\beta R} \log_+ \frac{\lambda\vert F\vert_2 } {\beta R} 
 \right) \right] 
=\int_0^\infty P\left(\exp\left( 
  \frac  {\vert F\vert_2 }{\beta R}\log_+ \frac{\lambda\vert F\vert_2 }{\beta R} \right)\geq t\right) \: dt
} 
\\
& = & \int_{-\infty}^\infty 
 e^y 
 P\left( 
 \frac{\vert F\vert_2 }{\beta R}\log_+ \frac{\lambda\vert F\vert_2 }{\beta R} \geq y 
 \right) dy 
\\
& \leq & 1 
 + 
 \int_0^\infty 
 e^y 
 P\left( 
 \frac{\vert F\vert_2 }{\beta R}\log \frac{\lambda\vert F\vert_2 }{\beta R} \geq y 
 \right) dy 
\\
&\leq &1+\int_{\beta R/\lambda}^\infty P\left(\frac{\vert F\vert_2 }{\beta R}\log\frac{\lambda\vert F\vert_2 } {\beta R}\geq \frac{x}{\beta R}\log\frac {\lambda x}{\beta R}\right)\frac {1+ \log\frac{\lambda x}{\beta R}} {\beta R} e^{\frac x{\beta R} \log\frac{\lambda x}{\beta R}} \: dx 
\\
&\leq&1+\frac{2}{\beta R}\int_{\beta R/\lambda}^\infty P(\vert F\vert_2 \geq x) 
 e^{\frac x{\beta R} \log\frac {\lambda x}{\beta R}} 
 \log \frac{\lambda x}{\beta R} 
 dx\\
&\leq& 1+\frac 2{\beta R}\int_{\beta R/\lambda}^\infty e^{\frac {x-2E[\vert F\vert_2 ]}{\beta R}} 
 e^{\frac x{\beta R}\log\frac{\lambda x}{\beta R}} 
 \log \frac{\lambda x}{\beta R} 
\\ 
 & & \times 
 \exp 
 \left(-\left(\frac {x-2E[\vert F\vert_2 ]}{\beta R}+\frac{\tilde{\alpha}_R^2}{\beta^2R^2}\right)\log\Big(1+\frac{(x-2E[\vert F\vert_2 ])\beta R}{\tilde{\alpha}_R^2}\Big)\right) 
 dx\\
&\leq& 1+\frac 2{\beta R}\int_{\beta R/\lambda-2E\vert F\vert_2 }^\infty e^{\frac {u}{\beta R}} 
 e^{\frac {u+2E[\vert F\vert_2 ]}{\beta R}\log \frac{\lambda (u+2E[\vert F\vert_2 ])}{\beta R}} 
\\ 
 & & \times 
 \exp\left(-\left(\frac {u}{\beta R}+\frac{\tilde{\alpha}_R^2}{\beta^2R^2}\right) 
 \log \left(1+\frac{u\beta R}{\tilde{\alpha}_R^2}\right)\right) 
 \log \frac{\lambda (u+2E[\vert F\vert_2 ])}{\beta R} du 
. 
\end{eqnarray*}
 It then suffices to study the dominant term in the above integral: 
\begin{align*} 
&\int_{\beta R/\lambda}^\infty \exp 
 \left( 
 \frac u{\beta R}-\frac u {\beta R}\log \Big( 
 1+\frac{u\beta R}{\tilde{\alpha}_R^2} 
 \Big) 
 \right) 
 e^{\frac {u}{\beta R} \log (1+\frac{\lambda u}{\beta R})} 
  \: du\\
&= \int_{\beta R/\lambda}^\infty \exp 
\left( 
 \frac u{\beta R}-\frac u{\beta R}\log 
 \frac{1+u\beta R/\tilde{\alpha}_R^2}{1+\lambda(u+2E[\vert F\vert_2 ])/(\beta R)} 
 \right) 
 du. 
\end{align*}
 Since 
$$\lim_{u\to \infty} 
 \log 
 \frac{1+u\beta R/\tilde{\alpha}_R^2}{1+\lambda(u+2E[\vert F\vert_2 ])/(\beta R)} 
 = 
 \log 
 \frac{\beta^2 R^2}{\lambda \tilde{\alpha}_R^2} > 1 
, 
$$ 
 for $\lambda<\beta^2 R^2/(e\tilde{\alpha}_R^2)$, 
 the convergence of the integral follows from 
$$ 
\int_{\beta R/\lambda}^\infty \exp 
 \left( 
 \frac u{\beta R}-\frac u {\beta R} 
 \log \frac{\beta^2 R^2}{\lambda \tilde{\alpha}_R^2} 
 \right) 
 du<\infty.
$$ 
\end{Proof}
 Since $\tilde{\alpha}_R^2$ given in Corollary~\ref{corol:norm-dev-borné} does not depend on the dimension, the condition on $\lambda$ in the above corollary 
 is also dimension free.

\subsection*{Random vectors with non-independent components} 
 First, we obtain the following from Corollary~\ref{corol:NP1}:
\begin{corollary}
\label{corol:NP2} 
 Let $f:\real^n\to\real$ be $\ell^2$-Lipschitz$(c)$ and let 
 $F=(F_1,\ldots ,F_n)$ such that 
 $\sum_{j=1}^n 
 |D_{(i,y)} F_j (\om)|^2 \leq K^2$, $P\otimes\nu_i(d\omega , dy)$-a.e., $i=1,\ldots ,n$, 
 for some $K\geq 0$ and 
$$ 
 \left\| 
 \sum_{i,j=1}^{n} \|D_{(i,\cdot )} F_j (\omega ) \|_{L^2(Y;\nu_i)}^2 
 \right\|_{L^\infty (\Omega^X , P)} 
 \leq \tilde{\alpha}^2 
.$$ 
 Then 
\begin{equation} 
\label{eq:Dev-Nico-ter}
P(f(F)-E[f(F)]\geq x)\leq e^{\frac{x}{cK}} 
 \left(1+\frac{xK}{\tilde{\alpha}^2}\right)^{-\frac{x}{cK} - \frac{\tilde{\alpha}^2}{K^2} }, \quad x>0. 
\end{equation}
\end{corollary} 
\begin{Proof} 
 Note that since $f$ is $\ell^2$-Lipschitz$(c)$ we have for $G=f(F)$:
\begin{eqnarray*}
( D_{(i,y)}G(\om) )^2 
 & =& 
 (f(F_1(\om\cup\{(i,y)\}),\dots, F_i(\om\cup\{(i,y)\}), \dots , F_n(\om\cup\{(i,y)\})) 
\\ 
 & & -f(F_1(\om),\dots, F_i(\om), \dots F_n(\om)))^2 
\\
&\leq& c^2 \sum_{j=1}^n 
 |F_j(\om\cup\{(i,y)\})-F_j (\om)|^2 
\\ 
& = & 
 c^2 \sum_{j=1}^n 
 |D_{(i,y)} F_j (\om)|^2 
.
\end{eqnarray*}
 So that $\vert D_{(i,y)}G (\omega ) \vert \leq cK$ 
 and $\sum_{i=1}^n \left\|D_{(i,\cdot )}G (\omega ) 
 \right\|_{L^2(Y;\nu_i)}^2 \leq c^2\tilde{\alpha}^2$, $P(d\omega )$-a.s., 
 and Corollary~\ref{corol:NP1} applies to $G$.
\end{Proof} 
\noindent 
 From Corollary~\ref{corol:NP2} 
 we can derive an exponential integrability 
 result for the Euclidean norm of a vector of 
 arbitrary functionals on $X= \{1,\ldots , n\} \times Y$, 
 provided $\nu_1=\cdots =\nu_n$ has support in $B_Y (0,R)$. 
 This completes the sharper result stated in Corollary~\ref{corol:intégrabilité1} 
 in the case of independent components. 
 However, in the infinitely divisible case, it is slightly 
 less sharp as 
 Corollary~3 of \cite{Rosinski95}. 
\begin{corollary} 
\label{cor:intégrabilité2}
 Let $\nu = \nu_1=\cdots = \nu_n$ have bounded support 
 in $B_Y(0,R)$, and let 
 $F=(F_1,\ldots ,F_n)$ be 
 a vector of $n$ (non-necessarily independent) random functionals. 
 Assume that 
$$\sum_{j=1}^{n} \frac{\vert D_{(i,y)}F_j (\omega ) \vert^2}{\vert y \vert_Y^2} 
 \leq \tilde{\alpha}^2 < \infty 
, \quad P\otimes \nu_i (d\omega,dy)-a.e., \ i=1,\ldots , n, 
$$ 
 and 
$$ 
 \left\| 
 \sum_{i,j=1}^{n} \|D_{(i,y)}F_j (\omega ) \|_{L^2 (Y; \nu_i)}^2 
 \right\|_{L^\infty (\Omega^X , P)} 
 \leq \tilde{\tilde{\alpha}}^2 < \infty 
. 
$$ 
 Then \eqref{eq:int-exp1} holds for $0<\lambda<\tilde{\alpha}^2R^2/\tilde{\tilde{\alpha}}^2$: 
\begin{equation}
\label{eq:int-exp2}
E\left[ 
 \exp\left( 
 \frac{\vert F\vert_2 }{\tilde{\alpha} R}\log_+ 
 \frac{\lambda^2\vert F\vert_2^2}{\tilde{\alpha} R}\right) \right] 
 <\infty. 
\end{equation}
\end{corollary} 
\begin{Proof}
 First, note that 
 $|D_{(i,y)}F_j|\leq \tilde{\alpha} \vert y\vert_Y\leq \tilde{\alpha} R$, since  $D_{(i,y)}F_j$ 
 is zero for $\vert y\vert_Y>R$, $\nu$ being supported on $B_Y(0,R)$. We can thus apply 
 Corollary~\ref{corol:NP2} and get 
\begin{equation} 
\label{eq:Dev-Nico-ter2} 
 P(f(F)-E[f(F)]\geq x)\leq e^{\frac{x}{\tilde{\alpha} R}} 
 \left(1+\frac{x\tilde{\alpha} R}{\tilde{\tilde{\alpha}}^2}\right)^{-\frac{x}{\tilde{\alpha} R} - \frac{\tilde{\tilde{\alpha}}^2}{\tilde{\alpha}^2 R^2} }, \quad x>0, 
\end{equation} 
 which is \eqref{eq:Dev-Nico-ter} with $K=\tilde{\alpha} R$ and $c=1$. 
 Finally \eqref{eq:int-exp2} follows from \eqref{eq:Dev-Nico-ter2} 
 as \eqref{eq:int-exp1} follows from \eqref{eq:dev-norme-bornée} 
 in Corollary~\ref{corol:intégrabilité1}. 
\end{Proof}
 In the previous corollary, $\tilde{\tilde{\alpha}}$ is dimension dependent, unlike 
 Corollary~\ref{corol:intégrabilité1}, so that the exponential integrability 
 is not dimension free in the dependent case. 
 As an application of Corollary~\ref{corol:NP2} we obtain 
 an upper large deviation bound in the dependent case, 
 for random functionals with bounded support. 
\begin{corollary} 
\label{corol:LD2} 
 Let $\nu := \nu_1=\cdots=\nu_n$ have bounded support in $B_Y (0,R)$ and let 
 $F=(F_1,\dots, F_n)$ be a vector of $n$ (non-necessarily independent) 
 random functionals. 
 Assume that 
$$\sum_{j=1}^{n} \frac{\vert D_{(i,y)}F_j (\omega ) \vert^2}{\vert y \vert_Y^2} 
 \leq \tilde{\alpha}^2 < \infty 
, \quad P\otimes \nu_i (d\omega,dy)-a.e., \ i=1,\ldots , n. 
$$ 
 and 
$$\left\| 
 \sum_{i,j=1}^{n} \|D_{(i,y)}F_j (\omega ) \|_{L^2 (Y; \nu_i)}^2 
 \right\|_{L^\infty (\Omega^X , P)} 
 < \infty 
. 
$$ 
 Then for any $\ell^2$-Lipschitz$(c)$ function $f:\real^n\to \real$, we have 
\begin{equation*} 
\limsup_{x\to\infty}\frac{\log P(|f(F)|\geq x)}{x\log x}\leq-\frac c {\tilde{\alpha} R}.
\end{equation*}
\end{corollary}
 When restricted to single Poisson integrals, the previous result recovers 
 the upper estimate of Corollary~4 in \cite{Rosinski95} 
 in which a deviation result is obtained 
 for the norm of infinitely divisible vector with L\'evy measure having a bounded support. 
 See also \cite{jurek} for related results in the framework of 
 large deviations for Poisson stochastic integrals. 
\section{Quadratic Wiener functionals} 
\label{sec:appli-quadratic} 
 The results of the previous section apply in particular to quadratic Wiener 
 functionals since they have infinitely divisible laws, cf.~\cite{matsumoto}, 
 and can be represented as Poisson stochastic integrals with 
 finite variance. 
 Note that exact estimates for the tail probabilities of (quadratic) Wiener 
 functionals have been obtained in \cite{imkeller2}, see also 
 \cite{borell2}, \cite{ledouxchaos}, \cite{mayer}, 
 \cite{mckean2}. 
 Here we present dimension free results for norms of vectors of independent quadratic 
 functionals. 
 In this section we take $X=\real$. 
\subsection*{Second order Wiener integrals} 
 It is well-known (see e.g. \cite{matsumoto}) 
 that every centered quadratic Wiener functional can be 
 determined by a symmetric Hilbert-Schmidt operator $A:L^2(\real_+)\to L^2 (\real_+)$ with 
 eigenvalues $(a_k)_{k\in\inte }$ and a complete orthonormal basis 
 of eigenvectors $(h_k)_{k\in\inte}$  in $L^2(\real_+)$. 
 In particular it can be expressed as 
 a second order Wiener integral $J_2(f_2)$ with respect 
 to a standard Brownian motion $(B_t)_{t\in \real_+}$, with 
$$ 
J_2(f_2) 
 = \frac{1}{2} 
 \sum_{k=0}^\infty a_k \left( \left( \int_0^\infty h_k (t) dB_t\right)^2 
 - 1\right), 
$$ 
 where the series converges in $L^2(\Omega^X)$, and $f_2$ has 
 the decomposition 
$$ 
 f_2 = \frac{1}{2} \sum_{k=0}^\infty a_k h_k\otimes h_k, 
$$ 
 converging in $L^2(\real_+^2)$. 
 Note that $J_2(f_2)$ is distinct from the double Poisson stochastic integral 
 $I_2(f_2)$. 
 The variance of $J_2(f_2)$ is 
$$\Var [J_2(f_2)] = \Vert f_2 \Vert_{L^2(\real_+^2)}^2 
 = \frac{1}{4} \sum_{k=0}^\infty a_k^2. 
$$ 
 In the sequel we consider a vector $(J_2^1(f_2^1),\ldots ,J_2^n(f_2^n))$ 
 of mutually independent second order Wiener integrals of 
 $f_2^1,\ldots ,f_2^n\in L^2(\real_+^2)$ with respect to possibly 
 different Brownian motions. 
 Denote also by $(a_k^i)_{k\in\inte}$ the eigenvalues associated to 
 $J_2^i(f_2^i)$, $i=1,\ldots ,n$. 
 For each $i=1,\ldots ,n$, 
 $J_2^i(f_2^i)$ is infinitely divisible, integrable, and centered 
 with L\'evy measure 
\begin{equation}
\label{eq:mesure-I2}
\nu_i (dy) 
 = {\bf 1}_{\{y>0\}} 
 \sum_{k; a_k^i>0} 
 \frac{1}{2\vert y \vert } e^{-y /a_k^i} dy
 + 
 {\bf 1}_{\{y < 0\}} 
 \sum_{k; a_k^i<0} 
 \frac{1}{2\vert y \vert } e^{- y /a_k^i} dy
, 
\end{equation}
 cf. Theorem~2 of \cite{matsumoto}. 
 Hence from \eqref{eq:vaX2}, $J_2^i(f_2^i)$ has the representation 
\begin{equation} 
\label{eq:vaX2.1}
J_2^i(f_2^i) (\omega_i ) = \int_{-\infty}^\infty y (\omega_i (dy)-\nu_i (dy)), 
\end{equation} 
 as a single Poisson stochastic integral. 
 Denote by 
\begin{align*}
& a_+ = \max_{1\leq i \leq n} \max_{k, a_k^i>0} a_k^i 
,\quad a_- = \max_{1\leq i \leq n} \max_{k, a_k^i<0} (-a_k^i ) 
, 
 \quad a = \max_{1\leq i \leq n} \max_{k \in\inte } |a_k^i| 
, 
\end{align*}
 the maxima of the spectral 
 radii associated to $J_2^1(f_2^1),\ldots ,J_2^n(f_2^n)$. 
\\ 
\noindent 
 In the next proposition we apply Corollary~\ref{thd2} to obtain 
 a deviation result for $\ell^1$-Lipschitz functions of quadratic 
 Wiener functionals. 
 Note that Corollary~4 of \cite{HRB} (or Proposition~\ref{df1} applied 
 to Poisson stochastic integrals) 
 would yield dimension free deviation results when $g$ is $\ell^2$-Lipschitz, 
 however with an additional range condition. 
\begin{prop} 
\label{prop:I2}
 Let $(J_2^1(f_2^1),\ldots ,J_2^n(f_2^n))$ 
 be a vector of independent second order Wiener integrals. 
 For any $\ell^1$-Lipschitz$(c)$ function $g$: 
$$ 
 P( g(J_2^1(f_2^1),\ldots ,J_2^n(f_2^n)) 
 -E[g(J_2^1(f_2^1),\ldots ,J_2^n(f_2^n)) ] \geq x) 
 \leq 
 \exp \left( 
 -\int_0^x h^{-1}(s) ds\right), 
$$ 
 $x>0$, where $h^{-1}$ is the inverse of the function 
$$ 
 h(t) = \frac{1}{2} \sum_{i=1}^n 
 \sum_{k=0}^\infty 
 \frac{ct(a_k^i)^2}{1-ct|a_k^i|}, 
 \quad t\in [0,(ca)^{-1}) 
. 
$$ 
 Moreover, 
\begin{eqnarray}
\nonumber 
\lefteqn{ 
\hskip -30pt
P( g(J_2^1(f_2^1),\ldots ,J_2^n(f_2^n)) 
 -E[g(J_2^1(f_2^1),\ldots ,J_2^n(f_2^n)) ] \geq x) 
}
\\ 
\label{eq:int2}
 & \leq & 
\exp \left( 
 - \frac{x}{ac} 
 + \frac{2\sum_{i=1}^n \Vert f_2^i\Vert_{L^2(\real_+^2)}^2}{a^2c} 
 \log \left( 1 + \frac{a x}{2\sum_{i=1}^n \Vert f_2^i\Vert_{L^2(\real_+^2)}^2} \right) 
 \right) 
\\ 
\nonumber 
 & \leq & 
\exp \left( 
 - \frac{1}{c} 
 \left(1-\frac{\log 3}{2}\right) \min 
 \left( \frac{x}{a} 
 , \frac{x^2}{4 \sum_{i=1}^n \Vert f_2^i\Vert_{L^2(\real_+^2)}^2} \right) 
\right) 
, 
\qquad 
 x>0. 
\end{eqnarray} 
\end{prop}
\begin{Proof} 
 From Corollary~\ref{thd2} and 
 \eqref{eq:vaX2.1}, 
 \eqref{eq:mesure-I2}, 
 \eqref{eqr}, we have 
\begin{eqnarray*} 
h(t) & \leq & c \sum_{i=1}^n 
 \int_{-\infty}^\infty 
 \vert y \vert (e^{t c \vert y \vert }-1)
 \nu_i (dy) 
\\ 
& = & \frac{c}2 \sum_{i=1}^n 
 \sum_{\begin{subarray}{c}k=0\\a_k^i>0\end{subarray}}^\infty
 \int_0^\infty 
 (e^{t c y }-1) e^{-y/a_k^i} dy 
 + 
 \frac{c}2 \sum_{i=1}^n \sum_{\begin{subarray}{c}k=0\\a_k^i<0\end{subarray}}^\infty
 \int_{-\infty}^0 
 (e^{-t c y }-1) e^{-y/a_k^i} dy 
\\ 
 &= & \frac{1}{2} \sum_{i=1}^n 
 \sum_{\begin{subarray}{c}k=0\\a_k^i>0\end{subarray}}^\infty   \frac{ct(a_k^i)^2}{1-cta_k^i}+ 
 \frac{1}{2} \sum_{i=1}^n \sum_{\begin{subarray}{c}k=0\\a_k^i<0\end{subarray}}^\infty   \frac{ct(a_k^i)^2}{1+cta_k^i}
\\
&= & \frac{1}{2} \sum_{i=1}^n \sum_{k=0}^\infty   \frac{ct(a_k^i)^2}{1-ct|a_k^i|}. 
\end{eqnarray*} 
 Then one can take 
$$ 
 h(t) = \frac{1}{2} \sum_{i=1}^n \sum_{k=0}^\infty 
 \frac{ct(a_k^i)^2}{1-cta}
 \leq 
 \frac{2ct}{(1-cta)} \sum_{i=1}^n 
 \Vert f_2^i\Vert_{L^2(\real_+^2)}^2 , 
 \quad t\in [0,(ca)^{-1}),  
$$ 
 and in this case, 
$$
h^{-1}(t) = \frac{t}{ 
 c a t
 + 2 c\sum_{i=1}^n \Vert f_2^i \Vert_{L^2(\real_+^2)}^2 
 },
$$
 from which \eqref{eq:int2} follows with explicit computations.
\end{Proof} 
\noindent 
 Alternatively, and since 
 $\lim_{t\to \infty} h^{-1}(t) = 1/(ca)$, we have for any $\varepsilon >0$ 
\begin{eqnarray} 
\nonumber 
\lefteqn{ 
\hspace{-3cm} P( g(J_2^1(f_2^1),\ldots ,J_2^n(f_2^n)) 
 -E[g(J_2^1(f_2^1),\ldots ,J_2^n(f_2^n)) ] \geq x) 
} 
\\ 
\label{ee2} 
&\leq & C_1(c,n,\varepsilon ) \exp \left( 
 -x((ca)^{-1}-\varepsilon) \right), 
\end{eqnarray} 
 $x>0$, for some constant $C_1(c,n,\varepsilon )$ depending on $c$, $n$ 
 and $\varepsilon$. 
 It follows that there is a constant $C_2 (c,n,\lambda )$ such that 
\begin{equation} 
\label{tcn} 
E[e^{\lambda \vert g (J_2^1(f_2^1),\ldots ,J_2^n(f_2^n)) \vert 
 }] < C_2 (c,n,\lambda ) 
 < \infty 
, 
\end{equation} 
 for all $\lambda < 1/(a c)$, and every $\ell^1$-Lipschitz$(c)$ function $g:\real^n\to \real$. 
 In fact, \eqref{eq:int2} implies that for all $\lambda < 
 - 1 - \frac{2}{a^2c} \sum_{i=1}^n \Vert f_2^i\Vert_{L^2(\real_+^2)}^2$, 
$$ 
E\left[\exp \left( 
 -\frac{1}{ac} \vert g (J_2^1(f_2^1),\ldots ,J_2^n(f_2^n)) \vert 
 + \lambda \log (1+ \vert g (J_2^1(f_2^1),\ldots ,J_2^n(f_2^n)) \vert ) \right) 
 \right] 
 < C_3 (c,n,\lambda ) 
, 
$$ 
 for some $C_3 (c,n,\lambda ) < \infty$. 
 For the supremum of $J_2^1(f_2^1),\ldots ,J_2^n(f_2^n)$, 
 which is a Lipschitz function with respect to the $\ell^\infty$-norm, 
 hence with respect to the $\ell^1$-norm, the previous corollary 
 can be strengthened by making use of \eqref{hj}. 
\begin{prop}
\label{prop:I2-sup}
 Let $(J_2^1(f_2^1),\ldots ,J_2^n(f_2^n))$ 
 be a vector of independent second order Wiener integrals. 
 Then, 
$$ 
 P( \sup (J_2^1(f_2^1),\ldots ,J_2^n(f_2^n)) 
 -E[\sup(J_2^1(f_2^1),\ldots ,J_2^n(f_2^n)) ] \geq x) 
\leq 
\exp \left( 
 -\int_0^x h^{-1}(s) ds\right), 
$$ 
 $x>0$, where $h^{-1}$ is the inverse of the function 
$$ 
 h(t) = \frac{1}{2} \sum_{i=1}^n
 \sum_{\begin{subarray}{c}k=0\\ a_k^i>0\end{subarray}}^\infty 
 \frac{t(a_k^i)^2}{1-ta_k^i}, 
 \quad t\in [0,1/a_+) 
. 
$$ 
 Moreover, 
\begin{eqnarray}
\nonumber 
\lefteqn{ 
\hskip -30pt
P( \sup(J_2^1(f_2^1),\ldots ,J_2^n(f_2^n)) 
 -E[\sup(J_2^1(f_2^1),\ldots ,J_2^n(f_2^n)) ] \geq x) 
}
\\ 
\label{eq:int2-sup}
 & \leq & 
\exp \left( 
 - \frac{x}{a_+} 
 + \frac{2\sum_{i=1}^n \Vert f_2^i\Vert_{L^2(\real_+^2)}^2}{(a_+)^2} 
 \log \left( 1 + \frac{a_+ x}{2\sum_{i=1}^n \Vert f_2^i\Vert_{L^2(\real_+^2)}^2} \right) 
 \right) 
\\ 
\nonumber 
 & \leq & 
\exp \left( 
 - \left(1-\frac{\log 3}{2}\right) \min 
 \left( \frac{x}{a_+} 
 , \frac{x^2}{4 \sum_{i=1}^n \Vert f_2^i\Vert_{L^2(\real_+^2)}^2} \right) 
\right) 
, 
\qquad 
 x>0. 
\end{eqnarray} 
\end{prop}
\begin{Proof} Follow the lines of the proof 
 of Proposition~\ref{prop:I2} starting from \eqref{eq:h-sup} instead of \eqref{eqr}. 
\end{Proof}
 Note that in dimension one and for second order Wiener integrals, 
 \eqref{eq:int2-sup} above implies the upper 
 deviation bound of \cite{mayer} (Example~5.1), 
 since 
$$a_+ = 2 \sup_{\Vert h \Vert_{L^2(\real_+)} = 1} 
  \langle f_2,h\otimes h\rangle_{L^2(\real_+^2)} 
. 
$$ 
 Counterparts of \eqref{ee2}, \eqref{tcn} for $\sup(J_2^1(f_2^1),\ldots ,J_2^n(f_2^n))$ 
 can be derived in the same way. 
 Our next result is a first lower bound. 
\begin{prop} 
\label{ll2} 
 Let $(J_2^1(f_2^1),\ldots ,J_2^n(f_2^n))$ 
 be a vector of (centered) mutually independent quadratic Wiener functionals. 
 For any $b\in (0,1)$, there exists $x_{b} >0$ such that 
$$ 
 P(\vert 
 (J_2^1(f_2^1),\ldots ,J_2^n(f_2^n)) 
 \vert_\infty \geq x) 
 \geq 
 \frac{1 - b}{2x} 
 \sum_{i=1}^n 
 a^i e^{-x/a^i } 
, \quad x> x_{b} 
, 
$$ 
 with $a^i= \max_{k \in\inte } |a_k^i|$, $1\leq i \leq n$.  
\end{prop} 
\begin{Proof} 
 Let $F_1,\ldots ,F_n$ be $n$ independent 
 random variables with respective 
 distribution $ID(m_1,0,\nu_1)$, ..., $ID(m_n,0,\nu_n)$. 
 We have for $x>0$: 
\begin{eqnarray*} 
P(\vert (F_1 ,\ldots , F_n ) \vert_\infty \geq x) 
 & \geq & P(\exists i\in \{1,\ldots ,n\} \ : \ 
 \vert F_i \vert \geq x) 
\\ 
 & = & 1 - P( \vert F_i \vert < x, \ 1\leq i \leq n ) 
\\ 
 & = & 1 - \prod_{i=1}^n 
 P( \vert F_i \vert < x ) 
. 
\end{eqnarray*} 
 Writing 
$$F_i = F^+_i + F^-_i - m_i$$ 
 with 
$$F^+_i (\omega_i ) = \int_0^\infty y \omega_i (dy ) , 
 \qquad 
 F^-_i (\omega_i ) 
 = 
 \int_{-\infty}^0 y \omega_i (dy ) , 
 \qquad 
 m_i = \int_{-\infty}^\infty y \nu_i (dy ) , 
$$ 
 we have 
\begin{eqnarray*} 
 P(F_i < x) 
 & = & \int_0^\infty 
 P(F^+_i < x+m_i +y)dP(F_i^-=-y) 
\\ 
 & \leq & \int_0^\infty 
 P(\omega_i ([x+m_i +y , \infty )) = 0 ) dP(F^-_i =-y) 
\\ 
 & = & \int_0^\infty 
 \exp (-\nu_i ([x+m_i +y , \infty ) )) dP(F^-_i =-y) 
. 
\end{eqnarray*} 
 Here, 
$$ 
 \nu_i ([x,\infty )) = 
 \int_{x}^\infty \sum_{k, a_k^i> 0} \frac{1}{2\vert y\vert} e^{-y /a_k^i} \; dy 
\sim_{x\to \infty} N_+\int_{x}^\infty \frac{e^{-y/a_+^i}}{2y} dy\\
\sim_{x\to \infty} a_+^iN_+^i\frac{e^{-x/a_+^i}}{2x} 
, 
$$ 
 where $f(x)\sim_{x\to\infty} g(x)$ means that $\lim_{x\to +\infty }f(x)/g(x) =1$, 
 and $a_+^i=\max_{k, a_k^i>0} |a_k^i|$, $N_+^i=\#\{k, \: a_k^i=a_+^i\}$. 
 Hence for all $b'\in (0,1)$ there exists $x_{i,b'}>0$ such that 
$$ 
 P(F_i < x) 
 \leq 
 \exp \left(-  
 a_+^iN_+^i\frac{e^{-x/a_+^i}}{2x} (1-b') 
 \right) 
, \qquad 
 x> x_{i,b'}. 
$$ 
 Similarly we have 
$$ 
\nu_i ((-\infty ,-x]) 
  = 
 \int_{-\infty}^{-x}\sum_{k, a_k^i<0} \frac{1}{2\vert y\vert} e^{-y /a_k^i} \; dy 
 \sim_{x\to \infty} 
 a_-^iN_-^i\frac{e^{-x/a_-^i}}{2x} 
, 
$$ 
 with 
 $a_-^i=\max_{k, a_k^i<0} |a_k^i|$ and $N_-^i=\#\{k, \: a_k^i=-a_-^i \}$. 
 Hence $x_{i,b'}$ can be chosen such that 
$$ 
 P(F_i > - x) 
 \leq 
 \exp \left(-  
 a_-^iN_-^i\frac{e^{-x/a_-^i}}{2x} 
 (1-b') \right) 
, 
 \qquad 
 x> x_{i,b'} 
, 
$$ 
 thus 
$$ 
 P(\vert F_i \vert < x) 
 \leq 
 \exp \left(-  
 a^i \frac{e^{-x/a^i}}{2x} 
 (1-b') \right) 
, \qquad 
 x> x_{i,b'} 
. 
$$ 
 For $x > \max \{ x_{1,b'}, \ldots , x_{n,b'}\}$. 
 It follows that 
$$ 
P(\vert (F_1 ,\ldots , F_n ) \vert_\infty \geq x) 
 \geq 
 1 - \exp \left(- 
 \sum_{i=1}^n
 a^i \frac{e^{-x/a^i}}{2x} 
 (1-b') \right) 
, 
$$ 
 and so for any $b\in (0,1)$, there exists $x_b>0$ such that 
$$ 
P(\vert (F_1 ,\ldots , F_n ) \vert_\infty \geq x) 
 \geq 
 \frac{1-b}{2x} \sum_{i=1}^n
 a^i e^{-x/a^i} 
, 
 \qquad 
 x> x_b. 
$$ 
\end{Proof} 
 Note that without the independence assumption on 
 $(J_2^1(f_2^1),\ldots ,J_2^n(f_2^n))$, a similar argument 
 leads to the estimate 
$$ 
 P(\vert 
 (J_2^1(f_2^1),\ldots ,J_2^n(f_2^n)) 
 \vert_\infty \geq x) 
 \geq 
 a (1-b) 
 \frac{e^{-x/a }}{2x} 
, 
$$ 
 for any $b\in (0,1)$ and $x$ large enough. 
 A version of Proposition \ref{ll2} 
 can also be stated for $\sup(J_2^1(f_2^1),\ldots ,J_2^n(f_2^n))$. 
 For $n=1$, and for second order Wiener integrals, 
 this also implies the lower deviation bound obtained in Example~5.1 in \cite{mayer}. 
\begin{prop} 
\label{ll2-sup} 
 Let $(J_2^1(f_2^1),\ldots ,J_2^n(f_2^n))$ 
 be a vector of independent quadratic Wiener functionals. 
 For any $b\in (0,1)$, there exists $x_{b} >0$ such that 
$$ 
 P(\sup(J_2^1(f_2^1),\ldots ,J_2^n(f_2^n)) \geq x) 
 \geq 
 \frac{1 - b}{2x} 
 \sum_{i=1}^n 
 a_+^i e^{-x/a_+^i } 
, \quad x> x_{b} 
, 
$$ 
 with $a_+^i= \max_{k \in\inte, a_k^i>0 } a_k^i$, $1\leq i \leq n$.  
\end{prop} 
\begin{Proof} 
 We follow the lines of proof of Proposition \ref{ll2}. 
 Let $F_1,\ldots ,F_n$ be $n$ independent 
 random variables with respective 
 distribution $ID(m_1,0,\nu_1)$, ..., $ID(m_n ,0,\nu_n )$. 
 We have for $x>0$: 
$$
P(\sup(F_1 ,\ldots , F_n)\geq x) 
 \geq 1 - \prod_{i=1}^n 
 P(F_i < x ) 
, 
$$ 
 which leads to 
$$ 
P(\sup(F_1 ,\ldots , F_n )\geq x) 
 \geq 
 1 - \exp \left(- 
 \sum_{i=1}^n
 a_+^i \frac{e^{-x/a_+^i}}{2x} 
 (1-b') \right) 
, 
$$ 
 for $x$ sufficiently large. 
 Hence, for any $b\in (0,1)$, there exists $x_b>0$ such that 
$$ 
P(\sup(F_1 ,\ldots , F_n ) \geq x) 
 \geq 
 \frac{1-b}{2x} \sum_{i=1}^n
 a_+^i e^{-x/a_+^i} 
, 
 \qquad 
 x> x_b. 
$$ 
\end{Proof} 
 Without the independence assumption on 
 $(J_2^1(f_2^1),\ldots ,J_2^n(f_2^n))$ we get 
\begin{equation} 
\label{le} 
 P(\sup 
 (J_2^1(f_2^1),\ldots ,J_2^n(f_2^n)) 
  \geq x) 
 \geq 
 a_+ 
 (1-b) 
 \frac{e^{-x/a_+ }}{2x} 
, 
\end{equation} 
 for any $b\in (0,1)$ and $x$ large enough. 
 In the next corollary we derive an exact tail estimate for the 
 $\ell_p$-norm of vectors of independent quadratic Wiener functionals, recovering, 
 in the special case of second order integrals, the result obtained 
 in \cite{AG} for non-decoupled Gaussian chaos, see also \cite[Cor. 3.9]{LT}. 
\begin{corollary} 
\label{c4.3} 
 Let $p\in [1,\infty ]$, and let $(J_2^1(f_2^1),\ldots ,J_2^n(f_2^n))$ 
 be a vector of independent quadratic Wiener functionals. 
 Then 
\begin{equation} 
\label{*.**} 
\lim_{x\to +\infty}\frac{\log P(\vert (J_2^1(f_2^1),\dots, J_2^n(f_2^n))\vert_p \geq x)}{x} 
 = -\frac 1 a.
\end{equation} 
\end{corollary}
\begin{Proof} 
 For any $b\in(0,1)$, from Proposition~\ref{ll2} and Proposition~\ref{prop:I2}, 
 there exists $x_b>0$ such that 
\begin{eqnarray*} 
\lefteqn{ 
 (1 - b ) 
 \frac{e^{-x/a }}{2x} 
  \sum_{i=1}^n 
 a^i
 \leq 
 P(\vert 
 (J_2^1(f_2^1),\ldots ,J_2^n(f_2^n)) 
 \vert_\infty \geq x) 
} 
\\
& \leq & 
 P(\vert 
 (J_2^1(f_2^1),\ldots ,J_2^n(f_2^n)) 
 \vert_p \geq x) 
\\
& \leq & 
\exp \left( 
 - \frac{x-M}{a} 
 + \frac{2\sum_{i=1}^n \Vert f_2^i\Vert_{L^2(\real_+^2)}^2}{a^2} 
 \log \left( 1 + \frac{a (x-M)}{2\sum_{i=1}^n\Vert f_2^i\Vert_{L^2(\real_+^2)}^2} \right) 
 \right), 
\end{eqnarray*}
 $x>\max (x_b,M)$, with $M=E[\vert 
 (J_2^1(f_2^1),\ldots ,J_2^n(f_2^n)) 
 \vert_1]$. 
\end{Proof} 
 Note that for $n=1$ and for second order Wiener integrals, 
 the above result coincides with Theorem~2.2 of \cite{borell2} (see also 
 \cite{ledouxchaos} and \cite{mckean2}), 
 since $a/2$ is also the strong operator norm of the linear map canonically associated 
 to $f_2$, i.e. 
$$a = 
 2 \sup_{\Vert h \Vert_{L^2(\real_+)} = 1} 
 \vert \langle f_2,h\otimes h\rangle_{L^2(\real_+^2)} \vert 
. 
$$ 
 A result of \cite{borell2} states that 
\begin{equation} 
\label{eq1} 
\lim_{x\to +\infty}\frac{\log P(\sup_{t\in \real_+ } 
 \vert J_m (f_m^t)\vert \geq x)}{x} 
 = -\frac{1}{2 \sup_{t\in \real_+} \Vert f_m^t\Vert_{L^2(\real_+^m)}}
, 
\end{equation} 
 provided $(J_m (f_m^t))_{t\in \real_+}$ is a process of $m$-th 
 order integrals with a.s. continuous sample paths (see also Remark~4.3 in 
 \cite{mayer}). 
 It is clear that for $n=1$, $m=2$ and $p=+\infty$, \eqref{*.**} 
 and \eqref{eq1} coincide. 
 However, \eqref{eq1} does not imply \eqref{*.**}, since as is 
 well known the process $(J_m(f_m^t))_{t\in \real_+}$ cannot 
 be jointly measurable and have independent components. 
 For the supremum of $J_2^1(f_2^1),\dots, J_2^n(f_2^n)$ we similarly have: 
\begin{corollary} 
\label{c4.3-sup} 
 Let $(J_2^1(f_2^1),\ldots ,J_2^n(f_2^n))$ 
 be a vector of independent quadratic Wiener functionals, then 
$$
\lim_{x\to +\infty}\frac{\log P(\sup (J_2^1(f_2^1),\dots, J_2^n(f_2^n))\geq x)}{x} 
= -\frac 1 {a_+}.
$$ 
\end{corollary}
\begin{Proof}
 Apply Proposition~\ref{prop:I2-sup} 
 and Proposition~\ref{ll2-sup}. 
\end{Proof} 
 A left deviation estimate 
 for $\sup(J_2^1(f_2),\ldots ,J_2^n(f_2))$ can be independently obtained from 
$$
P(\sup (J_2^1(f_2),\ldots ,J_2^n(f_2))\leq x)=\prod_{k=1}^n P(J_2^k(f_2)\leq x) 
, 
$$
 which can then be 
 estimated from Proposition~\ref{ll2-sup} or Proposition~9.17 of \cite{imkeller2}. 
 Counterparts of \eqref{eq:int2-sup} in Proposition~\ref{prop:I2-sup}, 
 as well as Proposition~\ref{ll2-sup}, Corollary \ref{c4.3-sup} and \eqref{le} 
 can also be stated for the left deviation of $\inf(J_2^1(f_2),\ldots ,J_2^n(f_2))$, 
 replacing $a_+^i$ by $a_-^i$, $i=1,\ldots ,n$, and 
 $a_+$ by $a_-$. 
 Since 
$$a_- = - 2 \inf_{\Vert h \Vert_{L^2(\real_+)} = 1} 
  \langle f_2,h\otimes h\rangle_{L^2(\real_+^2)} 
, 
$$ 
 this will imply 
 the one-dimensional left tails of \cite{mayer} (Example~5.1). 
 For an arbitrary norm $\Vert \cdot \Vert$ on $\real^n$, we have 
$$\Vert x \Vert = \left\| 
 \sum_{i=1}^{i=n} x_i e_i \right\| 
 \leq 
 \vert x \vert_1 
 \max_{1\leq i \leq n} \Vert e_i\Vert . 
$$ 
 Hence, 
$$\limsup_{x\to + \infty} 
 \frac{\log P(\Vert (J_2^1(f_2^1),\dots, J_2^n(f_2^n)) \Vert \geq x)}{x} 
 \leq -\frac{1}{a \max_{1\leq i \leq n } \Vert e_i\Vert} 
. 
$$ 
 Similarly, since $\Vert x \Vert \geq c(n) \vert x \vert_\infty$ for some $c(n)>0$, 
$$\liminf_{x\to + \infty} 
 \frac{\log P(\Vert (J_2^1(f_2^1),\dots, J_2^n(f_2^n)) \Vert \geq x)}{x} 
 \geq -\frac{1}{a c(n)} 
. 
$$ 
 For the Euclidean norm, we also have 
 the following dimension free deviation inequality 
 obtained from Proposition~\ref{df} for an i.i.d. 
 vector. The independent but non identically distributed case is 
 similar with more notation. 
\begin{prop} 
\label{prop:I2-bis} 
 Let $(J_2^1(f_2),\ldots ,J_2^n(f_2))$ be an i.i.d. vector of second order Wiener integrals, 
 and let $b\in(0,1)$. 
 Then, 
\begin{equation}
\label{eq:J2}
P( \vert (J_2^1(f_2),\ldots ,J_2^n(f_2)) \vert_2 
 -2E[\vert (J_2^1(f_2),\ldots ,J_2^n(f_2)) \vert_2 ] \geq x)
\leq 
e^{-(1-b)\frac x a +K_b}, \quad x>0, 
\end{equation}
 and 
\begin{equation}
\label{eq:J3}
P( \vert (J_2^1(f_2),\ldots ,J_2^n(f_2)) \vert_2 
 -2E[\vert (J_2^1(f_2),\ldots ,J_2^n(f_2)) \vert_2 ] \geq x)
\leq 
e^{-(1-b)\frac x a}, 
 \quad
 x\geq \frac{2a}b K_{b/2} 
, 
\end{equation}
 where
\begin{eqnarray*}
K_b&=&
-\frac{16 \Vert f_2\Vert_{L^2(\real_+^2)}^2}{a^2}\log b 
-8\Vert f_2\Vert_{L^2(\real_+^2)}^2\left(\frac{2}{a^2}+\frac{1}{(E[\vert J_2(f_2)\vert ])^2}\right) (1-b)\\
&& +\frac{4\Vert f_2\Vert_{L^2(\real_+^2)}^2}{(E[\vert J_2(f_2)\vert ])^2} \left(\frac{1-b^2}{b^2}\right) 
. 
\end{eqnarray*}
\end{prop}
\begin{Proof} 
 Applying Proposition~\ref{df} with $\beta=1$ gives 
\begin{eqnarray} 
\label{eq:minhh}
\lefteqn{ 
\hspace{-1cm} 
 P( \vert (J_2^1(f_2),\ldots ,J_2^n(f_2)) \vert_2 
 -2E[ \vert (J_2^1(f_2),\ldots ,J_2^n(f_2)) \vert_2 ] \geq x) 
} 
\\ 
\nonumber 
& \leq & 
 \exp\Big(-\int_0^x h^{-1} (s)ds\Big) 
 = 
 \min_{0<t<1/a} \exp\Big(-tx+\int_0^t h(s)ds\Big) 
, \quad x>0, 
\end{eqnarray} 
 where
\begin{eqnarray*} 
 h(t) & = & 
 4
 \sum_{\begin{subarray}{c}k=0\\a_k>0\end{subarray}}^\infty
 \int_0^\infty 
 (e^{t y }-1)
 e^{-y/a_k} dy 
 + 
 \frac{1}{(E[\vert J_2(f_2)\vert ])^2} 
 \sum_{\begin{subarray}{c}k=0\\a_k>0\end{subarray}}^\infty
 \int_0^\infty 
 y^2 
 (e^{t y }-1)
 e^{-y/a_k} dy 
\\ 
&&+4 
 \sum_{\begin{subarray}{c}k=0\\a_k<0\end{subarray}}^\infty
 \int_{-\infty}^0 
 (e^{-t y }-1) e^{-y/a_k} dy 
 + 
 \frac{1}{(E[\vert J_2(f_2)\vert ])^2} 
 \sum_{\begin{subarray}{c}k=0\\a_k<0\end{subarray}}^\infty
 \int_{-\infty}^0 
 y^2 
 (e^{-t y }-1)
 e^{-y/a_k} dy 
\\ 
 & = & 
 4
 \sum_{k=0}^\infty
 \frac{ta_k^2}{1-t|a_k|} 
 + 
 \frac{2}{(E[\vert J_2(f_2)\vert ])^2} 
 \sum_{k=0}^\infty
 |a_k|^3\left(\frac{1}{(1-t|a_k|)^3}-1\right) 
\\ 
 & \leq & 
 \frac{16 \Vert f_2\Vert_{L^2(\real_+^2)}^2}a
 \frac{ta}{1-ta} 
 + 
 \frac{8a\Vert f_2\Vert_{L^2(\real_+^2)}^2}{(E[\vert J_2(f_2)\vert ])^2} 
 \left(\frac{1}{(1-ta)^3}-1 \right) 
. 
\end{eqnarray*}
 Letting 
$$
A=\frac{16 \Vert f_2\Vert_{L^2(\real_+^2)}^2}{a^2} \quad \mbox{ and } \quad B=\frac{8\Vert f_2\Vert_{L^2(\real_+^2)}^2}{(E[\vert J_2(f_2)\vert ])^2} 
, 
$$
 we have 
\begin{eqnarray*}
\int_0^t h(s) ds & \leq & -A\left(\log(1-ta)+ta\right)+B\left(\frac 1{2(1-ta)^2}-ta-\frac 12\right)\\
&=&-A\log(1-ta)-(A+B)ta+ B\frac {ta}{(1-ta)^2} -\frac{B}{2} \frac{(ta)^2}{(1-ta)^2} 
. 
\end{eqnarray*}
 Taking $t=(1-b)/a$, the min in \eqref{eq:minhh} is bounded by 
$$
\min_{0<t<1/a}\exp \Big(-tx+\int_0^t h(s) ds \Big) \leq \exp\Big(-(1-b)\frac{x}{a}+K_b\Big) 
,  
$$
where 
$$
K_b=\int_0^{\frac{1-b}a} h(s) ds=-A\log b -(A+B) (1-b)+B \frac{1-b}{b^2}-B\frac{(1-b)^2}{2b^2}
, 
$$
and \eqref{eq:J2} follows. 
 Taking $t=(1-b/2)/a$ in \eqref{eq:minhh} yields 
$$
\min_{0<t<1/a}\exp\Big(-tx+\int_0^t h(s) ds\Big) \leq \exp\Big(-(1-b/2)\frac{x}{a}+K_{b/2}\Big)\leq \exp\Big(-(1-b)\frac{x}{a}\Big),
$$
 $x\geq \frac {2a} b K_{b/2}$, and \eqref{eq:J3} follows.
\end{Proof} 
 Note that the growth of $K_b$ is in $1/b^2$, as $b\to 0$. 
\subsection*{Square norm of Brownian paths on $[0,T]$} 
 An example of quadratic Wiener functional for which the coefficients $(a_k)_{k\in\inte}$ 
 can be explicitly computed is given by 
 the (compensated) integrated squared Brownian motion 
$$
\eufrak{h}_T = \int_0^T (B(t))^2dt - \frac{T^2}{2}
$$ 
 on the interval $[0,T]$. 
 In this case, from \cite{manabe} or \S3.1.1 of \cite{matsumoto}, 
 we have $a_k = \frac{4T^2}{(2k+1)^2\pi^2}$, $k\geq 0$, and 
 the above results apply with $a=\frac{4T^2}{\pi^2}$ and 
$$
\sum_{k=0}^\infty a_k^2  = \frac{T^2}{2}. 
$$ 
 Letting $(\eufrak{h}_T^1,\ldots ,\eufrak{h}_T^n)$ be a vector 
 of i.i.d. copies of $\eufrak{h}_T$, Proposition \ref{prop:I2} states 
 in this case that for any $\ell^1$-Lipschitz$(c)$ function $g:\real^n\to \real$: 
\begin{eqnarray*}
\lefteqn{ 
\hskip -45pt
P( g(\eufrak{h}_T^1,\ldots ,\eufrak{h}_T^n) 
 -E[g(\eufrak{h}_T^1,\ldots ,\eufrak{h}_T^n) ] \geq x) 
}
\\ 
 & \leq & 
\exp \left( 
 - \frac{\pi^2 x}{4cT^2} 
 + \frac{n\pi^4}{32T^2c} 
 \log \left( 1 + \frac{8x}{n\pi^2} \right) 
 \right),
\\ 
 & \leq & 
 \exp \left( 
 - \left( 1- \frac{\log 3}{2} \right) 
 \min \left( 
 \frac{\pi^2 x}{4cT^2} 
 , 
 \frac{x^2}{nc T^2} 
 \right) 
 \right),
 \qquad x>0, 
\end{eqnarray*} 
 and for all $\varepsilon >0$: 
\begin{equation*} 
P( g(\eufrak{h}_T^1,\ldots ,\eufrak{h}_T^n) 
 -E[g(\eufrak{h}_T^1,\ldots ,\eufrak{h}_T^n) ] \geq x) 
\leq C_1 (c,n,T,\varepsilon ) \exp \left( 
 -x(\pi^2/(4cT^2)-\varepsilon) \right), 
\end{equation*} 
 $x>0$. 
 It follows that 
$$
E[e^{\lambda \vert 
 g(\eufrak{h}_T^1,\ldots ,\eufrak{h}_T^n) \vert }] < C_2 (c,n,T,\lambda ) 
, 
$$ 
 for all $\lambda < \pi^2/(4cT^2)$ 
 and some constant $C_2 (c,n,T,\lambda ) < \infty$. 
 For any $b\in (0,1)$ we also get the lower bound 
$$ 
P( \vert (\eufrak{h}_T^1,\ldots ,\eufrak{h}_T^n)\vert_\infty 
 \geq x) 
\geq n(1-b)\frac{\pi^2}{8 x T^2} e^{-x\pi^2/(4T^2)}
, 
$$ 
 for $x$ sufficiently large, 
 and from Corollary \ref{c4.3}: 
$$ 
\lim_{x\to +\infty}\frac{\log P(\vert (\eufrak{h}_T^1,\ldots ,\eufrak{h}_T^n)\vert_p \geq x)}{x} 
 = -\frac{\pi^2}{4T^2}, 
 \qquad p\in [1,\infty ].
$$
 For the Euclidean norm and from Proposition \ref{prop:I2-bis}, 
 we again have dimension free deviation inequalities: for any $b\in (0,1)$, 
$$
P( \vert (\eufrak{h}_T^1,\ldots ,\eufrak{h}_T^n) \vert_2 
 -2E[\vert (\eufrak{h}_T^1,\ldots ,\eufrak{h}_T^n) \vert_2 ] \geq x) 
\leq \exp \left( -\frac{\pi^2 (1-b)}{4T^2}x 
 + K_b\right), \quad x\geq 0, 
$$ 
 and 
$$
P( \vert (\eufrak{h}_T^1,\ldots ,\eufrak{h}_T^n) \vert_2 
 -2E[\vert (\eufrak{h}_T^1,\ldots ,\eufrak{h}_T^n) \vert_2 ] \geq x) 
\leq 
 \exp \left( -\frac{\pi^2 (1-b)}{4T^2}x 
 \right), 
 \qquad
 x\geq \frac{8T^2}{\pi^2 b} K_{b/2} 
, 
$$ 
 where 
\begin{eqnarray*}
K_b&=&- \frac{\pi^4}{8T^2}\log b-T^2\left(\frac{\pi^4}{8T^4}+\frac{1}{(E[\vert \eufrak{h}_T\vert ])^2}\right)(1-b) +\frac{T^2}{2(E[\vert \eufrak{h}_T\vert ])^2}\frac{1-b^2}{b^2} 
. 
\end{eqnarray*}

\subsection*{Sample variance of Brownian motion on $[0,T]$} 
 A second example is given by the (compensated) sample variance of Brownian motion 
$$ 
\eufrak{v}_T = \int_0^T \left(B(t) - \frac{1}{T} \int_0^T B(s)ds\right)^2dt - \frac{T^2}{6}
$$ 
 on the interval $[0,T]$. 
 From \cite{manabe}, or \S3.3.1 of \cite{matsumoto}, 
 $a_k = \frac{T^2}{k^2\pi^2}$, $k\geq 1$, $a_0=0$, 
 $a=\frac{T^2}{\pi^2}$ and 
$$
\sum_{k=1}^\infty a_k^2  = \frac{T^2}{6}. 
$$ 
 Letting $(\eufrak{v}_T^1,\ldots ,\eufrak{v}_T^n)$ be a vector 
 of i.i.d. copies of $\eufrak{v}_T$, we have from Proposition \ref{prop:I2} 
 that for every $\ell^1$-Lipschitz$(c)$ function $g:\real^n\to \real$: 
\begin{eqnarray*}
\lefteqn{ 
\hskip -45pt
P( g(\eufrak{v}_T^1,\ldots ,\eufrak{v}_T^n) 
 -E[g(\eufrak{v}_T^1,\ldots ,\eufrak{v}_T^n) ] \geq x) 
}
\\ 
 & \leq & 
\exp \left( 
 - \frac{\pi^2 x}{cT^2} 
 + \frac{n\pi^4}{6T^2c} 
 \log \left( 1 + \frac{6x}{n\pi^2} \right) 
 \right),
\\ 
 & \leq & 
\exp \left( 
 - \left( 1- \frac{\log 3}{2}\right) 
 \min \left( \frac{\pi^2 x}{cT^2} 
 , 
 \frac{6x^2}{2ncT^2} 
 \right) 
 \right) 
, \qquad x>0, 
\end{eqnarray*} 
 and for all $\varepsilon >0$: 
\begin{equation*} 
P( g(\eufrak{v}_T^1,\ldots ,\eufrak{v}_T^n) 
 -E[g(\eufrak{v}_T^1,\ldots ,\eufrak{v}_T^n) ] \geq x) 
\leq C_1 (c,n,T,\varepsilon ) \exp \left( 
 -x(\pi^2/(cT^2)-\varepsilon) \right), 
\end{equation*} 
 $x>0$. 
 It follows that 
$$
E[e^{\lambda \vert 
 g(\eufrak{v}_T^1,\ldots ,\eufrak{v}_T^n) \vert }] < 
 C_2 (c,n,T,\lambda ) 
, 
$$ 
 for all $\lambda < \pi^2/(cT^2)$ 
 and some $C_2 (c,n,T,\lambda ) < \infty$. 
 For any $b\in (0,1)$ we also get the lower bound 
$$ 
P( \vert (\eufrak{v}_T^1,\ldots ,\eufrak{v}_T^n)\vert_\infty 
 \geq x) 
\geq n(1-b)\frac{\pi^2}{2 x T^2} e^{-x\pi^2/T^2}
, 
$$ 
 for $x$ sufficiently large, 
 and from Corollary \ref{c4.3}: 
$$ 
 \lim_{x\to +\infty}\frac{\log P(\vert (\eufrak{v}_T^1,\ldots ,\eufrak{v}_T^n)\vert_p \geq x)}{x} 
 = -\frac{\pi^2}{T^2} 
 , 
 \qquad 
 p\in [1,\infty ].
$$
 From Proposition \ref{prop:I2-bis}, for any $b\in (0,1)$: 
$$
P( \vert (\eufrak{v}_T^1,\ldots ,\eufrak{v}_T^n) \vert_2 
 -2E[\vert (\eufrak{v}_T^1,\ldots ,\eufrak{v}_T^n) \vert_2 ] \geq x) 
\leq
\exp \left( -\frac{\pi^2(1-b)}{T^2}x 
 + K_b\right),\quad  x\geq 0,
$$ 
 and 
$$
P( \vert (\eufrak{v}_T^1,\ldots ,\eufrak{v}_T^n) \vert_2 
 -2E[\vert (\eufrak{v}_T^1,\ldots ,\eufrak{v}_T^n) \vert_2 ] \geq x) 
\leq
 \exp \left( -\frac{\pi^2(1-b)}{T^2}x 
 \right),\quad  x\geq \frac{2T^2}{\pi^2 b} K_{b/2} 
,
$$ 
 where 
\begin{eqnarray*}
K_b&=&- \frac 2 3 \frac{\pi^4}{T^2}\log b-\frac{T^2}3\left(\frac{2\pi^4}{T^4}+\frac{1}{(E[\vert \eufrak{v}_T\vert ])^2}\right)(1-b) +\frac{T^2}{6(E[\vert \eufrak{v}_T\vert ])^2}\frac{1-b^2}{b^2} 
. 
\end{eqnarray*}
\subsection*{L\'evy's stochastic area} 
 Let $(B^1(t),B^2(t))_{t\in \real_+}$, be a two-dimensional Brownian motion. 
 L\'evy's stochastic area $S_T$ on $[0,T]$ is 
$$
S_T = \frac{1}{2}\int_0^T \big(B^1(t)dB^2(t)-B^2(t)dB^1(t)\big), 
$$ 
 cf. \cite{levy}. 
 For $S_T$, the expression of the coefficients $(a_k)_{k\in\inte }$ is 
 intricate (see \cite{matsumoto}), 
 hence, we can not directly specialize the results of Proposition \ref{prop:I2}, 
 \ref{prop:I2-bis} and Corollary \ref{c4.3} in that case. 
 However, since the L\'evy measure of $S_T$ has the analytic expression 
\begin{equation}
\label{eq:mesure-A-Levy}
\nu (dy) 
 = \frac{1}{2y\sinh \frac{\pi y}{T}} dy, 
\end{equation} 
 (cf. page 175 of \cite{levy}, \S3.2.1 of \cite{matsumoto} or Example~15.15 of \cite{Sato}), 
 we can derive 
 results similar to the ones obtained for general second order Wiener-It\^o integrals. 
\begin{prop} 
\label{p12} 
 Let $g:\real^n\rightarrow \real$ be $\ell^1$-Lipschitz$(c)$, and 
 let $(S_T^1,\ldots , S_T^n)$ be an i.i.d. vector of L\'evy's stochastic areas on $[0,T]$. 
 Then, 
\begin{eqnarray*} 
P( g(S_T^1,\ldots ,S_T^n ) -E[g(S_T^1,\ldots ,S_T^n)] \geq x) 
 & \leq & 
 \left(1+\frac{\pi x}{4 n c T} \right)^{4n} 
 \exp \left( 
 - \frac{\pi x }{cT} 
 \right) 
.
\end{eqnarray*} 
\end{prop}
\begin{Proof} 
 Using the representation of $S_T$ as the compensated Poisson stochastic integral 
\begin{equation}
\label{eq:A-Lévy}
 \int_{-\infty}^\infty 
 y (\omega (dy ) -\nu (dy)) 
\end{equation} 
 and \eqref{eqr} derived from Corollary~\ref{thd2}, we have 
\begin{eqnarray*} 
h(t) & = & nc \int_{-\infty}^\infty 
 \vert y \vert_Y (e^{t c \vert y \vert_Y }-1)
 \nu (dy) 
\\ 
& = & nc \int_{-\infty}^\infty 
 \vert y \vert_Y (e^{t c \vert y \vert_Y }-1)
 \frac{1}{2y\sinh \frac{\pi y}{T}} dy 
\\ 
 & = & 2 n c \int_0^\infty 
 \frac{e^{t c y }-1}{e^{\pi y/T}- e^{-\pi y/T}}
 dy 
\\ 
&\leq & 2nc\int_0^\infty \frac{e^{tcy}-1}{e^{\pi y/T}-1} dy\\
& \leq & \frac{2 n c^2 t T}{\pi} 
 \int_0^\infty 
 e^{- \frac{y}{2} (\frac{\pi}{T}-ct)} 
 dy 
\\ 
& \leq & \frac{4 n c^2 t T}{\pi 
 (\frac{\pi}{T}-ct)} 
, \qquad 0<t<\pi/(cT), 
\end{eqnarray*} 
 using the inequality 
\begin{equation}
\label{eq:exp}
 \frac{e^{u x}-1}{e^{v x}-1}\leq \frac u v e^{(u-v)x/2}, \qquad 
 x>0, 
\end{equation} 
for $0<u<v$. Hence 
\begin{eqnarray*} 
 - \int_0^x h^{-1} (s) ds & \leq & 
 - 
 \frac{\pi^2}{T} 
 \int_0^x 
 \frac{s}{4 n c^2 T+ c \pi s } ds 
 = 
 - 
 \frac{\pi}{cT} 
 x 
 + 
 4 n 
 \log \left(1+\frac{\pi x}{4 n c T} 
 \right) 
.
\end{eqnarray*} 
\end{Proof}
 The above result can also be obtained from Theorem~1 in 
 \cite{h3} in place of \eqref{eqr}. 
 Alternatively, we have 
 $\lim_{t\to \infty} h^{-1}(t) = \pi (cT)^{-1}$ 
 since $\lim_{t\to \pi (cT)^{-1} } h(t) = +\infty$. 
 Hence for all $\varepsilon >0$, we also derive as in \eqref{ee2} 
\begin{equation*} 
P( g(S_T^1,\ldots ,S_T^n ) -E[g(S_T^1,\ldots ,S_T^n)] \geq x) 
\leq C_1 (c,n,T,\varepsilon ) \exp \left( 
 -x(\pi(cT)^{-1}-\varepsilon) \right), 
\end{equation*} 
 $x>0$, for some constant $C_1 (c,n,T,\varepsilon )$ depending on $T$, $c$, $n$ 
 and $\varepsilon$.  
 This last inequality is not dimension free. 
 Nevertheless it yields 
\begin{equation} 
\label{nvr} 
 E\left[ 
 e^{\lambda \vert g (S_T^1,\ldots ,S_T^n ) \vert }\right] < 
 C_2 (c,n,T,\lambda ) < \infty , 
\end{equation} 
 for all $\lambda < \pi/(cT)$, and every $\ell^1$-Lipschitz$(c)$ 
 function $g:\real^n \to \real$. 
\begin{prop} 
\label{ll2.1} 
 Let $p\geq 1$. 
 For all $b\in (0,1)$ there exists $x_{b} >0$ such that 
$$ 
 P(\vert 
 (S_T^1,\ldots ,S_T^n) 
 \vert_\infty \geq x) 
 \geq 
 (1-b) 
 \frac{nTe^{-\pi x/T}}{2\pi x} 
, 
 \qquad 
 x> x_b. 
$$ 
\end{prop} 
\begin{Proof} 
 Given $n$ i.i.d. random variables $F_1,\ldots ,F_n$ 
 with distribution $ID(m,0,\nu)$, 
 we have from the proof of Proposition~\ref{ll2}: 
$$ 
P(\vert (F_1 ,\ldots , F_n ) \vert_\infty \geq x) 
 \geq 1 - (P( \vert F_1 \vert < x ))^n 
, 
$$ 
 while for $\nu$ given in \eqref{eq:mesure-A-Levy} we have the equivalence: 
\begin{eqnarray*}
\nu([x,\infty[)=\int_{x}^\infty \frac{dy}{2y \sinh\frac{\pi y}T}& 
 \sim_{x\to\infty} 
 &\int_{x}^\infty \frac{e^{-\pi y/T}}{y} dy\sim_{x\to\infty} 
 \frac{Te^{-\pi x/T}}{\pi x} 
, 
\end{eqnarray*} 
 hence for all $b'\in (0,1)$ there exists $x_{i,b'}$ such that 
$$ 
 P(\vert F_i \vert < x) 
 \leq 
 \exp \left(-  
 (1-b') \frac{Te^{-\pi x/T}}{2\pi x} 
 \right) 
, \qquad 
 x> x_{i,b'}, 
 \quad 1\leq i \leq n 
. 
$$ 
 Thus for $x$ large enough, 
$$ 
P(\vert (S_T^1 ,\ldots , S_T^n ) \vert_\infty \geq x) 
 \geq 
 1 - \exp \left(- 
 (1-b) \frac{nTe^{-\pi x/T}}{2\pi x} 
 \right) 
. 
$$ 
\end{Proof} 
 The corollary below is a direct consequence of Proposition~\ref{p12} and 
 Proposition~\ref{ll2.1}. It recovers a univariate result of \cite{berthuet} 
 and extends it to $\ell^p$-norms 
 of i.i.d. random vectors, independently of their dimension. 
 For non identical variables $S_{T_1}^1, \ldots , S_{T_n}^n$, 
 replace $T$ by $\max_{1\leq k \leq n } T_k$. 
\begin{corollary} 
\label{corol:lim-A-Levy}
 Let $p\in [1,\infty ]$. 
 Then, 
$$ 
 \lim_{x\to +\infty}\frac{\log P(\vert (S_T^1,\dots, S_T^n)\vert_p \geq x)}{x} 
 = -\frac \pi T.
$$ 
\end{corollary}
 Note that from the above results we have $a=T/\pi$ and 
 $\sum_{k=0}^\infty a_k^2 = T^2/4$. 
 Moreover as a consequence of Proposition~\ref{df}, we have: 
\begin{prop} 
\label{prop:A-Levy} 
 Let $\vert \cdot \vert_2$ denote the Euclidean norm on $\real^n$ 
 and let $b\in(0,1)$. 
 We have 
\begin{equation}
\label{eq:AL2}
P( \vert (S_T^1,\ldots ,S_T^n )\vert_2 -2E[\vert (S_T^1,\ldots ,S_T^n)\vert_2 ] \geq x) 
\leq 
\exp\Big(-\frac{(1-b)\pi}T x +K_b\Big), \quad x>0, 
\end{equation}
 and 
\begin{equation}
\label{eq:AL3}
P( \vert (S_T^1,\ldots ,S_T^n )\vert_2 -2E[\vert (S_T^1,\ldots ,S_T^n)\vert_2 ] \geq x) 
\leq 
\exp\Big(-\frac{(1-b)\pi}T x \Big), 
 \quad 
 x\geq \frac{2T}{\pi b} K_{b/2} 
, 
\end{equation}
 where 
\begin{eqnarray*}
K_b&=& -32 \log b-32(1-b)+\frac{16T^2}{\pi^2(E[\vert S_T^1\vert ])^2}\frac{(1-b)^2}{b^2} 
. 
\end{eqnarray*}
\end{prop}
\begin{Proof} 
 Since for $F$ given in \eqref{eq:A-Lévy}, $\beta=1$, from Proposition~\ref{df} 
 and \eqref{gth}, we have 
$$
P( \vert (S_T^1,\ldots ,S_T^n ) \vert_2 -2E[\vert (S_T^1,\ldots ,S_T^n)\vert_2 ] \geq x) 
\leq \exp \Big(-tx+\int_0^t h(s) ds \Big),
$$
 for all $x>0$ and $0<t<\pi/T$. 
 From \eqref{eq:exp} we have 
\begin{eqnarray*} 
 h(t) & = & 
 16 \int_0^\infty 
 e^{ty/2} 
 \frac{\sinh \frac{ty}{2}}{\sinh \frac{\pi y}{T} } 
 dy 
 + 
 \frac{4}{(E[\vert S_T^1\vert ])^2} 
 \int_0^\infty 
 y^2 
 e^{ty/2} 
 \frac{\sinh \frac{ty}{2}}{\sinh \frac{\pi y}{T}} 
 dy 
\\ 
 & \leq & 
 \frac{16 t T}{\pi} 
 \int_0^\infty 
 e^{- \frac{y}{2} (\frac{\pi}{T}-t)} 
 dy 
 + 
 \frac{4tT}{\pi(E[\vert S_T^1\vert ])^2} 
 \int_0^\infty 
 y^2 
 e^{- y (\frac{\pi}{T}-t)/2} 
 dy 
\\ 
 & \leq & 
 \frac{32T}{\pi} \frac{(t T/\pi)}{1-(tT/\pi)} 
 + 
 \frac{32T^3}{\pi^3(E[\vert S_T^1\vert ])^2}\frac{(tT/\pi)}{(1-(tT/\pi))^3} 
, 
\end{eqnarray*} 
 hence 
\begin{eqnarray} 
\label{eq:minhhh}
\lefteqn{ 
 P( \vert (S_T^1,\ldots ,S_T^n ) \vert_2 -2E[\vert (S_T^1,\ldots ,S_T^n)\vert_2 ] \geq x) 
} 
\\ 
\nonumber 
& \leq & 
\exp \left( 
 -tx-32\big((tT/\pi)+\log(1-(tT/\pi))\big)+\frac{16 T^2}{\pi^2(E[\vert S_T^1\vert ])^2}\frac{(tT/\pi)^2}{(1-(tT/\pi))^2} \right)  
. 
\end{eqnarray} 
 For all $b\in(0,1)$, \eqref{eq:AL2} follows by taking $t=(1-b)\pi/T$ in \eqref{eq:minhhh}, 
 and \eqref{eq:AL3} is a consequence of \eqref{eq:AL2} where $b$ is replaced with $b/2$. 
\end{Proof} 
\section{The infinite variance case}  
\label{sec:no-variance} 
 In \cite{HM}, deviation results have been 
 derived for Lipschitz functions of stable random vectors. 
 In this section, we extend these results to general Poisson functionals 
 under arbitrary intensity measures. 
 Deviations are now given with respect to a median rather 
 than to the mean (which may not exist). 
 For $A$ in ${\cal B}(X)$ (the Borel $\sigma$-field of $X$), 
 let $\nu_R (A) = \nu (A\cap B_X (0,R))$, 
 where $0$ denotes an arbitrary fixed point in $X$. 
 The proofs of the forthcoming results are inspired by that of Theorem~1 in \cite{HM}: 
 configurations are truncated and we will use 
 the following notation on the configuration space. For 
 a fixed $R>0$ and any $\omega\in \Omega^X$, let 
$$ 
\omega_R=\omega\cap B_X (0,R),\quad \omega^c_R=\omega\cap B_X (0,R)^c 
 = \{ 
 x\in  \omega \ : \  
 d_X(0,x) > R\} 
. 
$$
 Given a stochastic functional $F$ on the configuration space, 
 we also set 
$$
F_R(\omega ) =F(\omega_R)=F(\omega\cap B_X (0,R)), 
$$ 
 and denote by $\gamma$ a non-negative and non-increasing function such that 
$$ 
P(\{ \omega \in \Omega^X \ : \  
 \omega\cap B_X (0,R)^c\not =\emptyset \}) 
\leq \gamma (R), 
$$ 
 for all $R$ large enough. 
 The next Lemma will be used in the sequel. 
 It allows to control $m(F_R)-m(F)$ as in \cite{HM}.
\begin{lemma} 
\label{lemma:BHP} 
 Let $F$ be a stochastic functional on the configuration space such that 
 there exists a non-negative and non-decreasing 
 function $\tilde \beta$ (resp. non-increasing function $\tilde{\gamma}$) defined on $\real_+$, 
 such that for all $R$ greater than a given $R_0$: 
\begin{equation} 
\label{*1.1} 
 P\big(F_R-m(F_R)\geq \tilde \beta (R)\big) 
\leq \tilde{\gamma} (R) 
. 
\end{equation} 
 Then we have 
\begin{equation}
\label{eq:tech3} 
 m(F_R)-m(F)\leq \tilde \beta (R), 
\end{equation} 
 for all $R$ such that 
\begin{equation} 
\label{ccond} 
 R\geq 
 \max \left( 
 R_0 
 , 
 \inf_{0<\delta < 1/2} 
 \max \left( 
 \gamma ^{-1}(\delta) 
 , 
 \tilde{\gamma}^{-1} 
 \left( \frac{1}{2} - \delta \right) 
 \right) 
 \right) 
. 
\end{equation} 
\end{lemma} 
\begin{Proof} 
 The case $m(F)\geq m(F_R)$ being 
 trivial, we consider henceforth $m(F_R)\geq m(F)$. 
 Let $0 < \delta < 1/2$ and assume 
$$ 
R\geq \gamma ^{-1}(\delta). 
$$ 
 We have 
\begin{eqnarray*}
 0 < \frac{1}{2}-\delta & \leq & 
\frac{1}{2}-\gamma (R) 
\\ 
& \leq & P(F\leq m(F)) -P(\{ \omega \in \Omega^X \ : \ \omega^c_R\not =\emptyset \}) 
\\ 
 & = & 
 E[{\bf 1}_{\{F(\omega_R\cup \omega^c_R )\leq m(F) \}} 
 - {\bf 1}_{\{ \omega^c_R \not= \emptyset \}} 
 ] 
\\ 
 & \leq & 
E[{\bf 1}_{\{F(\omega_R )\leq m(F) \}} 
 ] 
\\ 
& = & 
 P(F_R \leq m(F))
\\
& = & P(F_R-m(F_R)\leq m(F)-m(F_R))\\
& = & P(-F_R-m(-F_R)\geq m(F_R)-m(F))
, 
\end{eqnarray*} 
 where we used the fact that $-m(F_R)$ is a median of $-F_R$. 
 Consider the decreasing function 
$$
 H_R(x)=P\big(-F_R-m(-F_R)\geq x\big), \quad x\in \real, 
$$
 and let 
$I_R(y)=\sup\{z\geq 0, H_R(z)\geq y\}$ denote 
 its inverse. 
 We have 
\begin{equation} 
\label{drv} 
m(F_R) -m(F)\leq I_R\big(P(F\leq m(F))-P( \{ \omega \in \Omega^X \ : \  
 \omega^c_R\not =\emptyset \} )\big) 
 \leq I_R \Big(\frac{1}{2} - \delta \Big). 
\end{equation} 
 Assume further that 
$$ 
 R\geq \tilde\gamma ^{-1} \Big(\frac{1}{2} - \delta \Big). 
$$ 
 From \eqref{*1.1} applied to $-F_R$, we have
\begin{eqnarray*} 
H_R(\tilde \beta (R)) 
 & = & 
 P(-F_R-m(-F_R)\geq \tilde \beta (R)) 
 \leq \tilde \gamma (R) 
\leq \frac{1}{2} - \delta, 
\end{eqnarray*}
 that is finally 
$I_R 
 \left( 
 \frac{1}{2} - \delta 
 \right) 
 \leq \tilde\beta (R)$, and from \eqref{drv}: 
$$ 
m(F_R)-m(F)\leq I_R \Big(\frac 12-\delta \Big) \leq \tilde \beta (R). 
$$ 
\end{Proof} 
 The next result provides a general deviation property for stochastic functionals 
 with infinite variance on Poisson space. 
\begin{theorem} 
\label{theo:BHP} 
 Let $F$ be a stochastic functional on the configuration space such that 
 there exists a non-negative and non-decreasing function $\beta$, defined on $\real_+$, 
 and a constant $C>0$ such that for all $R$ greater than a given $R_0$: 
\begin{description} 
\item{(i)} $\sup_{y\in B_X(0,R)} 
 \vert 
 D_y F (\omega ) 
 \vert  
 \leq \beta (R), \quad P(d\omega ) $-a.s., 
\item{(ii)} $\Vert DF \Vert_{L^\infty(\Omega^X, L^2(\nu_R))}^2 
 \leq C \beta^2 (R) \gamma (R)$. 
\end{description} 
 Then 
\begin{equation*}
P(F-m(F)\geq x)\leq (1 + Ce ) 
 \gamma \circ \beta ^{-1}(x/4), 
\end{equation*}
 for all 
$$
x\geq 2\: \beta \left( 
 \gamma^{-1}\left(\frac 1{2( 1 + Ce )}\right)\right).
$$ 
\end{theorem}
\begin{Proof} 
 Configurations are truncated 
 to deal on the one hand with 
 the functional restricted to the truncated configuration and on the 
 other hand with the rest of the configuration which 
 is controlled using the function $\gamma$. 
 We have 
\begin{eqnarray}
\nonumber
P(F-m(F)\geq x)&=&P(F-m(F)\geq x,\; \omega^c_R=\emptyset)+P(F-m(F)\geq x,\; \omega^c_R\not =\emptyset)\\
\label{eq:tech2}
& \leq &P(F_R-m(F)\geq x)+P(\{ \omega \in \Omega^X \ : \  \omega^c_R\not =\emptyset \} ).
\end{eqnarray}
 For the first term, in order to apply \eqref{eq:Dev-Nico} in Corollary~\ref{prop:dev-NP} 
 (which provides a deviation result 
 from the mean rather than from a median), let
$$ 
g(x) = (x-m(F_R))^+ \wedge r, 
 \quad x\in \real. 
$$ 
 Then $E[g(F_R)]\leq r P(F_R\geq m(F_R)) \leq r/2$. Moreover if 
 $F_R\geq m(F_R) +r$ then $g(F_R)\geq g(m(F_R) +r)\geq r$, 
 hence 
$$
\{F_R\geq m(F_R) +r\}\subset \{ g(F_R)\geq r\}, 
$$ 
 and 
\begin{equation} 
\label{and} 
P(F_R - m(F_R)\geq r ) \leq P(g(F_R) \geq r ) \leq P(g(F_R)-E[g(F_R)]\geq r/2) 
. 
\end{equation} 
 On the other hand, $g(F_R)$ satisfies 
$$ 
 D_y g(F_R)(\omega ) 
 \leq \vert g(F_R(\omega \cup \{ y \})) 
 - g(F_R(\omega )) 
 \vert 
 \leq 
 \vert 
 F_R(\omega \cup \{ y \})  
 - F_R(\omega ) 
 \vert 
 = 
 \vert 
 D_y F(\omega_R ) 
 \vert 
, 
$$ 
 since $g:\real\to \real$ is Lipschitz$(1)$. 
 Thus 
$$
\sup_{y \in B_X(0,R)} 
 D_yg(F_R)\leq \beta (R) 
 \quad 
\mbox{and} 
 \quad 
\|Dg(F_R)\|_{L^\infty(\Omega^X , L^2(\nu_R))}^2\leq C \gamma (R)\beta (R)^2 
, 
$$ 
 and from \eqref{and} and Corollary~\ref{prop:dev-NP} we get 
\begin{equation*} 
P(F_R - m(F_R)\geq x ) 
 \leq 
 e^{x/(2\beta (R))} 
 \left(1+\frac{x}{2 C \gamma (R)\beta (R)}\right)^{-x/(2\beta (R))} 
, 
\end{equation*} 
and taking $x=2\beta (R)$ we have: 
\begin{equation}
\label{eq:*1.2}
 P(F_R-m(F_R)\geq 2\beta (R)) 
\leq e \left(1+\frac{1}{C \gamma (R)}\right)^{-1} 
\leq e C \gamma (R),
\end{equation} 
 and from Lemma~\ref{lemma:BHP} with $\tilde \beta(R)=2\beta(R)$, $\tilde{\gamma} (R) = Ce \gamma (R)$ and condition \eqref{*1.1} given by \eqref{eq:*1.2}, we get: 
$$ 
m(F_R)-m(F)\leq 2\beta (R), 
$$ 
 i.e. using \eqref{eq:*1.2}: 
$$
P\big(F_R-m(F)\geq 4 \beta (R)\big) 
 \leq P\big(F_R-m(F_R)\geq 2\beta (R)\big)\leq Ce\: \gamma (R) 
, 
$$
 i.e. for $x \geq 4 \beta (R)$, we have 
\begin{equation} 
\label{esst} 
 P(F_R-m(F)\geq x)\leq 
 Ce\: \gamma \circ \beta ^{-1}\left(\frac x{4}\right) 
, 
\end{equation} 
 under condition (\ref{ccond}) which can be rewritten in terms of $x$ as 
$$ 
x\geq \max \left( 
 2 \beta ( 
 \gamma ^{-1} (\delta) ) , 
 2 \beta \left(\gamma ^{-1} \left( 
 \frac{1}{Ce} 
 \left( 
 \frac{1}{2} - \delta \right)\right)\right) 
 \right) 
,
$$
 i.e 
\begin{equation*}
\gamma\circ \beta^{-1}(x/2)\leq \min\left( 
 \delta,\; \frac{1}{Ce} 
 \left(\frac 12- \delta\right)\right).
\end{equation*} 
 The optimal bound with $\delta \in (0,1/2)$, 
 being obtained for $\delta_0=\disp\frac{1}{2( 1 + Ce )}\in (0,1/2)$, 
 i.e. the condition on $x$ becomes 
\begin{equation} 
\label{x} 
x\geq 2\beta\circ\gamma^{-1}\left(\frac 1{2( 1 + Ce )}\right).
\end{equation} 
 The estimate \eqref{esst}, together with 
$$P(\{ \omega \in \Omega^X \ : \  
 \omega^c_R\not =\emptyset \} )\leq \gamma (R) 
 \leq 
 \gamma \circ \beta ^{-1}\left(\frac x {4}\right), 
$$ 
 gives 
$$ 
P(F-m(F)\geq x)\leq \left( C 
 e+ 1 \right)\gamma \circ \beta ^{-1}\left(\frac x{4}\right) 
, 
$$ 
 using \eqref{eq:tech2}, under the condition \eqref{x}. 
\end{Proof}
 Note that in the hypotheses of Theorem~\ref{theo:BHP} it is sufficient to assume 
$$ 
\sup_{y\in B_X(0,R)} 
 \vert 
 D_yF (\omega_R) 
 \vert  
 \leq \beta (R),\quad P(d\omega)\mbox{-a.s.} 
, 
$$ 
 and 
$$\Vert DF_R \Vert_{L^\infty(\Omega^X , L^2(\nu_R))}^2 
 \leq C \beta^2 (R)\: \gamma (R) 
, 
$$ 
 instead of $(i)$ and $(ii)$. 
 The next corollary presents a particular and more tractable case 
 of Theorem~\ref{theo:BHP}. 
\begin{corollary} 
\label{c} 
 Let $F : \Omega^X \to \real$, and let 
$$
\gamma (R) =  1-e^{-\nu ( \{ y\in X \ : \ d_X(0,y) > R\} )}, \quad \quad R>0 
, 
$$ 
 and assume that 
$$
\sup_{y\in B_X(0,R) \leq R} 
 \vert D_y F \vert \leq C' R \quad \mbox{and} \quad 
\Vert DF \Vert_{L^\infty(\Omega^X, L^2(\nu_R))}^2 
 \leq C 
 R^2 \gamma (R), 
$$ 
 for all $R \geq R_0 >0$. 
 Then 
\begin{equation*}
P(F-m(F)\geq x)\leq \left(1 + \frac{Ce}{(C')^2} \right) 
 \gamma \left( 
 \frac{x}{4C'}\right), 
 \quad \quad 
x\geq 2\: C' \gamma^{-1}\left(\frac 1{2(1+eC/(C')^2)}\right).
\end{equation*}
\end{corollary} 
\bigskip 
 On $\real^n$ equipped with the Euclidean norm 
 $\vert \cdot \vert_2$, 
 consider an 
 $\ell^2$-Lipschitz$(c)$ function $f:\real^n \lto \real$ 
 and a $n$-dimensional infinitely divisible random 
 vector $F=(F_1,\dots, F_n)$ without Gaussian component 
 and with L\'evy measure $\nu$. 
 Let us apply Corollary~\ref{c} to the random functional $G=f(F)$, 
 where $F$ is given as in \eqref{eq:vaX2} by: 
$$ 
F=\left(\int_{ \{\vert y \vert_2 \leq 1\}} y_k\:  (\om(dy)-\nu (dy))+\int_{\{\vert y \vert_2 > 1\}} y_k  \: \om(dy )+b_k\right)_{1\leq k\leq n} 
.
$$ 
 For the gradient, we have if $y\notin \omega$:
\begin{eqnarray*} 
|D_y G(\om)|&=&|G(\om\cup\{y\})-G(\om)|\\
&=&\disp \left|f\left(\int_{\{\vert u\vert_2 \leq 1\}} u (\om(du)-\nu (du))+ 
 \int_{\{\vert u\vert_2 >1\}} \: \om(du) +y+b \right)-f(F)\right|\\
&\leq & c\vert y\vert_2 , 
\end{eqnarray*} 
 since $f$ is $\ell^2$-Lipschitz$(c)$, and we obtain $|D_y G(\om)|\leq cR$, 
 for $\vert y \vert_2 \leq R$. 
 In this case, for $G=f(F)$ the conclusion of Corollary~\ref{c} reads 
\begin{equation} 
\label{eq:but3}
P(G-m(G)\geq x)\leq (1 + Ce ) 
 \left( 
 1-\exp 
 \left( 
 -\nu 
 \left( 
 \left\{ u \in \real^n \ : \ \vert u \vert_2 > \frac{x}{4c} \right\} 
 \right) 
 \right)\right) 
. 
\end{equation} 
 When $f (x) = \vert x \vert_2$ is the Euclidean norm on $\real^n$, 
 Lemma~\ref{lemme:ineg-Levy} below also yields a lower bound on 
 $P(\vert F -m\vert_2 \geq x)$ 
 which has the same order as the upper bound \eqref{eq:but3}. 
\begin{lemma} 
\label{lemme:ineg-Levy} 
 Let $F$ be an infinitely divisible random vector $ID(b,0,\nu )$ 
 in $\real^n$, with median $m\in\real^n$. Then 
\begin{equation} 
\label{eq:ineg-Levy}
P(\Vert F -m\Vert \geq x)\geq \frac 1 4 \left(1-\exp ( 
 -\nu (\{u\in \real^n \ : \ \Vert u\Vert \geq 2x\})) 
 \right), \qquad x>0, 
\end{equation}
 where $\Vert \cdot \Vert$ denotes any norm on $\real^n$. 
\end{lemma}
\begin{Proof} 
 We start by assuming that $F$ is symmetric with median $0$. Then, since $F$ can be 
 taken to be the value $F(1)$ at time $1$ of 
 a L\'evy process $(F(t))_{0\leq t\leq 1}$ starting from $F(0)=0$, we have from 
 L\'evy's inequality: 
\begin{eqnarray*} 
\nonumber
P(\Vert F \Vert \geq x)&=&P\left(\left\Vert \sum_{k=1}^n F \left(\frac kn\right)-F \left(\frac{k-1}n\right)\right\Vert \geq x\right)\\
&\geq&\frac 1 2 P\left(\max_{1\leq j\leq n} \left\Vert F \left(\frac jn\right)-F \left(\frac{j-1}n\right)\right\Vert \geq x\right) 
. 
\end{eqnarray*} 
 Hence 
\begin{eqnarray} 
\nonumber
P(\Vert F \Vert \geq x)&\geq&\liminf_{n\to\infty} \frac 1 2 P\left(\max_{1\leq j\leq n} \left\Vert F \left(\frac jn\right)-F \left(\frac{j-1}n\right)\right\Vert \geq x\right)\\
\nonumber
&\geq&\frac 1 2 P\left(\liminf_{n\to\infty} \max_{1\leq j\leq n} \left\Vert F \left(\frac jn\right)-F \left(\frac{j-1}n\right)\right\Vert > x\right)\\
\nonumber 
&\geq&\frac 1 2 P\left(\max_{s\in[0,1]} \Vert F (s)-F (s^-)\Vert > x\right)\\
\label{eq:ineg-Sato-ex}
&\geq& \frac 1 2 (1-\exp\left(-\nu \{u\in \real^n \ : \ \Vert u\Vert \geq x\}\right)), 
\end{eqnarray} 
 where \eqref{eq:ineg-Sato-ex} is a $n$-dimensional 
 extension of Ex. 22.1 in \cite{Sato}, which relies on the fact that 
 if $\omega$ on $\real_+\times \real^n$ has a jump of $\Vert \cdot \Vert$-norm greater than 
 $x$, then 
 $\max_{s\in[0,1]} \Vert F (s)-F (s^-)\Vert > x$. 
 In the general case where $F$ is not necessarily symmetric we apply the above to 
 $F-G$, where $G$ denotes an independent copy of $F$, and use the inequality 
$$ 
P(\Vert F -m\Vert \geq x)  
 = \frac{1}{2} 
P(\Vert F -m\Vert \geq x) 
+ 
\frac{1}{2}  
P(\Vert G -m\Vert \geq x) 
 \geq 
\frac 1 2 P(\Vert F -G \Vert \geq 2x)
. 
$$ 
\end{Proof} 
 We now present several examples of L\'evy measures 
 $\nu$ for which 
 the function $\gamma$ can be explicitly computed, 
 and where $F$ has infinite variance, 
 i.e. $\int_{\real^n} \Vert y \Vert^2\nu (dy ) = \infty$, 
 but where $f(F)$ satisfies the above hypothesis for $f$ an 
 $\ell^2$-Lipschitz$(c)$ function. 
\begin{enumerate} 
\item Let $\real^n\setminus\{0\}$ be equipped with the measure given for $B\in{\cal B}(\real^n\setminus\{0\})$ by 
\begin{equation}
\label{eq:Levy-ex1}
\nu (B) =\int_{S^{n-1}}\sigma(d\xi)\int_0^\infty\ind_B(r\xi)\frac{|\log r|}{r^2} dr 
, 
\end{equation} 
 where $\sigma$ is again a spherical finite measure. 
 Since 
$$\int_{\{\vert y \vert_2 \leq 1\}} \vert y \vert_2 ^2 \nu (dy ) 
 = \sigma(S^{n-1}) < \infty$$ 
 and 
$$
\nu ( \{\vert x\vert_2 \geq 1\} ) = 
 \int_{S^{n-1}}\sigma(d\xi)\int_1^\infty \frac{|\log r|}{r^2}dr<\infty,
$$ 
 $\nu$ is a L\'evy measure. 
 Moreover $\int_{\{\vert y \vert_2 \geq 1\}}\vert y \vert_2 ^2 \nu (dy ) =\infty$, 
 hence $F$ has infinite variance. 
 As before: 
\begin{eqnarray*}
P(\{\omega \in \Omega^X \ : \  
 \om^c_R\not =\emptyset \})& = 
 & 
 1 -\exp\left(\int_{S^{n-1}}\sigma(d\xi)\int_{\{\vert r\xi\vert_2 \geq R\}}\frac{|\log r|}{r^2}dr\right)\\
&=&1-\exp\left(-\sigma(S^{n-1})\int_R^{\infty}\frac{\log r}{r^2} \: dr\right)\\
&\leq&\sigma(S^{n-1})\frac{1+\log R}{R}, \quad R>1.
\end{eqnarray*}
 Thus, choose $\gamma(R)=2\sigma(S^{n-1})\frac{\log R}{R}$. 
 On the other hand, 
\begin{eqnarray*}
\lefteqn{ 
\hspace{-1cm} 
\|D f(F) \|_{L^\infty(\Omega^X, L^2(\nu_R))}^2 
 \leq  
 \int_{\{\vert y \vert_2 \leq R\}} c^2\vert y \vert_2 ^2\: \nu (dy ) 
} 
\\
&=&-c^2\int_{S^{n-1}}\sigma(d\xi)\int_0^1 \log r\: dr+c^2\int_{S^{n-1}}\sigma(d\xi)\int_1^R\log r \:dr\\
&=& c^2\sigma(S^{n-1})(R\log R-R+2)\\
&\leq &c^2R^2\gamma (R)/2 
. 
\end{eqnarray*} 
\item Let $X=\real^n$, with the finite measure $\nu$ given for $B\in {\cal B}(\real^n)$ by 
\begin{equation} 
\label{eq:Levy-ex2}
\nu (B)
 =\int_{S^{n-1}}\sigma(d\xi) \int_0^\infty \ind_B(r\xi)\frac{e^{-1/(2r^2)}}{r^2\sqrt{2\pi}} dr.
\end{equation} 
 We have 
$$ 
\int_{\{\vert y \vert_2 \leq 1\}}\vert y \vert_2 ^2 \nu (dy ) =\sigma(S^{n-1}) \int_0^1 \frac{e^{-1/(2r^2)}}{\sqrt{2\pi}} dr<\infty, 
$$
 so that $\nu$ is a L\'evy measure. 
 The infinitely divisible random variable given by the Poisson stochastic 
 integral in \eqref{eq:vaX2} is thus another example of  a random variable without finite variance since 
$$ 
\int_{\{\vert y \vert_2 \geq 1\}} \vert y \vert_2 ^2 \nu (dy ) = \infty.
$$ 
 Once more: 
\begin{eqnarray*}
 P (\om^c_R\not =\emptyset)&=&1-\exp\left(-\int_{S^{n-1}}\sigma(d\xi)\int_{\{\vert r\xi\vert_2 \geq R\}} \frac{e^{-1/(2r^2)}}{r^2\sqrt{2\pi}}dr\right)\\
&\leq& \sigma(S^{n-1}) \int_R^\infty \frac{e^{-1/(2r^2)}}{r^2\sqrt{2\pi}}dr
\\ 
 &= &\sigma(S^{n-1})\int_0^{1/R} \frac{e^{-u^2/2}}{\sqrt{2\pi}} du  
. 
\end{eqnarray*} 
 Choose $\disp \gamma(R)=\frac{\sigma(S^{n-1})}{\sqrt{2\pi}R}$. 
 Moreover, 
\begin{eqnarray*}
\|D f(F) \|_{L^\infty(\Omega^X, L^2(\nu_R))}^2&\leq &
 c^2 \int_{B(0,R)} \vert y \vert_2 ^2\: \nu_R(dy ) 
\\ 
 & = & 
 c^2\int_{S^{n-1}}\sigma(d\xi)\int_0^R r^2 \frac{e^{-1/(2r^2)}}{r^2\sqrt{2\pi}}\: dr\\
&\leq&c^2\sigma(S^{n-1})\int_0^R \frac{e^{-1/(2r^2)}}{\sqrt{2\pi}}\: dr 
\\ 
 & =& 
 c^2\sigma(S^{n-1})\int_{1/R}^\infty \frac{e^{-u^2/2}}{\sqrt{2\pi}}\: \frac{du}{u^2}\\
&\leq&c^2\sigma(S^{n-1})\left( \frac{R\: e^{-1/(2R^2)}}{\sqrt{2\pi}}-
 \int_{1/R}^\infty\frac{e^{-u^2/2}}{\sqrt{2\pi}}\: du\right)\\
&\leq& c^2R^2\gamma(R) 
. 
\end{eqnarray*}
\item The above deviation results 
 for $f(F)$ with $F$ as in (\ref{eq:vaX2}) a stable or an infinitely 
 divisible random variable, and with L\'evy measure either given 
 by (\ref{eq:Levy-stable}) or (\ref{eq:Levy-ex1}) or (\ref{eq:Levy-ex2}), 
 continue to hold after minor changes for H\"older continuous functions of order $0<h<1$. 
 Indeed, for such a function $f$ we have 
$$ 
\vert 
 D_y f(F) \vert_2  \leq c\vert y\vert_2 ^h \leq c R^h,  \quad \vert y\vert_2 \leq R 
. 
$$
 For instance in the case of the L\'evy measure 
 \eqref{eq:Levy-ex1} we have 
\begin{eqnarray*}
\lefteqn{ 
\! \! \! \! \! \! \! \! \! \! \! \! \! 
\|D f(F) \|_{L^\infty(\Omega^X, L^2(\nu_R))}^2 \leq 
\int_{S^{n-1}}\sigma(d\xi) \int_0^R c^2r^{2h}\: \nu_R(dr) 
} 
\\ 
 & = & 
 c^2\sigma(S^{n-1})\int_0^R r^{2h-2} |\log r| dr
\\
&=& -\frac{c^2\sigma(S^{n-1})}{2} 
 \int_0^1 r^{2h-2}\log r \: dr+c^2\int_1^R r^{2h-2}\log r \:dr\\
&\leq&c^2\sigma(S^{n-1})\left(\frac{R^{2h-1}\log R}{2h-1}-\frac{R^{2h-1}}{(2h-1)^2}+\frac{1}{(2h-1)^2}+\frac 1 {2h+1} \right)
. 
\end{eqnarray*}
 We can thus apply Theorem~\ref{theo:BHP} to $G=f(F)$ with 
 (up to multiplicative constants) 
 the functions: 
$$ 
\beta(R)=R^h, 
 \quad \mbox{and} 
 \quad \gamma(R)=\frac{\log R} R.
$$
 A similar computation yields in the case of the L\'evy measure (\ref{eq:Levy-ex2}):
\begin{eqnarray*}
\|DF\|_{L^\infty(\Omega^X, L^2(\nu_R))}^2&\leq &
\int_{\{\vert x\vert_2 \leq R\}} c^2\vert x\vert_2 ^{2h} \nu_R (dr) 
 \\ 
 & \leq & 
 \frac{c^2}{\sqrt{2\pi}} \int_{S^{n-1}}\sigma(d\xi)
\int_0^R r^{2h-2}e^{-1/(2r^2)}\: dr 
\\
& \leq & c^2\sigma(S^{n-1})\frac{R^{2h-1}e^{-1/(2R^2)}}{(2h-1)\sqrt{2\pi}}.
\end{eqnarray*} 
 Once more, Theorem~\ref{theo:BHP} applies here, in the H\"older continuous case, with 
 up to multiplicative constants, 
$$ 
\beta(R)=R^h, \quad \mbox{and} 
 \quad \gamma(R)=\frac{e^{-1/(2R^2)}}{R\sqrt{2\pi}}.
$$
\end{enumerate} 
 Before turning to the case of stable intensity measures in the next 
 section, 
 we prove the following lemma for a general intensity measure $\nu$, 
 which is a generalization of Lemma~2 in \cite{HM}. 
\begin{lemma} 
\label{lemme:bis2}
 Let $F:\Omega^X \longrightarrow \real$ and $\alpha_2,\alpha_3,\alpha_4, K>0$, 
 such that 
\begin{description} 
\item{(i)} $\sup_{y\in X} \vert D_yF (\omega ) \vert\leq K < \infty$, $P(d\omega ) $-a.s.
\item{(ii)} $\Vert DF \Vert_{L^\infty(\Omega^X, L^k(\nu))}^k 
 \leq \alpha_k < \infty$, $k=2,3,4$. 
\end{description}
Assume moreover $\alpha_3\leq 2\alpha_4/K$ and $K^2\alpha_2/\alpha_4\geq 2$. Let $s_0$ be the (unique) positive solution of 
\begin{equation} 
\label{eq:s0}
\disp s\left(\alpha_2 -\frac{\alpha_4}{K^2}\right)=\frac{\alpha_4}{K^3} (e^{sK}-1).
\end{equation}
Let $x_0=3s_0(\alpha_2-\alpha_4/K^2)$. Then for all $x\leq x_0$, 
\begin{equation}
\label{eq:bis2-1}
P(F-E[F] \geq x)\leq \exp\left(-\frac{x^2}{6(\alpha_2-\alpha_4/K^2)}\right), 
\end{equation}
while for $x\geq x_0$,  
\begin{equation}
\label{eq:bis2-2}
P(F-E[F] \geq x)\leq K_0 \exp 
 \left(
\frac{x}{K}
-\left(\frac{x}{K}+\frac{3\alpha_4}{K^4}\right)\log\left(1+\frac{K^3x}{3\alpha_4}\right)\right), 
\end{equation}
 with 
\begin{equation}
\label{eq:K0-2}
K_0=
\exp 
 \left(
-\frac{x_0}{K}
+\left(\frac{x_0}{K}+\frac{3\alpha_4}{K^4}\right)\log\left(1+\frac{K^3x_0}{3\alpha_4}\right)-\frac {x_0^2}{6(\alpha_2-\alpha_4/K^2)}\right). 
\end{equation} 
\end{lemma}
\begin{Proof}
 From  Proposition~\ref{thd3} we have 
\begin{equation}
\label{eq:dev0}
P(F-E[F]\geq x)\leq \exp\left(-\int_0^x h^{-1}(s) ds \right), \med 0<x<h(t_0^-)
\end{equation}
 with $h$ given in \eqref{eq:h1}. Using the bounds 
 $|D_yF |\leq K$ and 
\begin{equation}
\label{eq:h-bound2}
e^{su}-1\leq su+\frac{s^2}{2} u^2 + \frac{e^{sK}-1-sK- s^2 K^2/2}{K^3} u^3, 
 \qquad 0 \leq u \leq K, \ s\geq 0, 
\end{equation} 
 we have 
\begin{eqnarray} 
\nonumber
\lefteqn{ 
 h(s) \leq \sup_{\omega, \omega'\in \Omega_X} 
\int_X \Big( 
 s|D_yF(\omega)|  |D_yF(\omega')| +\frac {s^2}2|D_yF(\omega)|^2 |D_yF(\omega')| 
} 
\\
\nonumber
& & +\frac{e^{sK}-1-sK-s^2 K^2/2}{K^3}|D_yF(\omega)|^3
|D_yF(\omega')| \Big) \nu(dy) 
\\
\nonumber
&\leq & s\sup_{\omega, \omega'\in \Omega_X} 
\int_X |D_yF(\omega)|  |D_yF(\omega')| \nu(dy) +\frac{s^2}2\sup_{\omega, \omega'\in \Omega_X} \int_X |D_yF(\omega)|^2  |D_yF(\omega')| \nu(dy) 
\\ 
\label {eq:alpha4}
& & +\frac{e^{sK}-1-sK- s^2K^2/2}{K^3}\sup_{\omega, \omega'\in \Omega_X} 
\int_X |D_yF(\omega)|^3
|D_yF(\omega')| \nu(dy).
\end{eqnarray} 
 Using the inequality $xy\leq x^p/p+y^q/q$ for $p^{-1}+q^{-1}=1$ and $x,y\geq 0$, we have 
 for $p=q=2$: 
\begin{eqnarray*}
\lefteqn{ 
\sup_{\omega, \omega'\in \Omega_X} \int_X |D_yF(\omega)|  |D_yF(\omega')| \nu(dy)
} 
\\
& \leq & \frac 12\sup_{\omega, \omega'\in \Omega_X} 
\int_X |D_yF(\omega)|^2\nu(dy) +\frac 12\sup_{\omega, \omega'\in \Omega_X} 
\int_X |D_yF(\omega')|^2 \nu(dy) 
 \leq \alpha_2
, 
\end{eqnarray*} 
 for $q=3$: 
\begin{eqnarray*} 
\lefteqn{ 
\sup_{\omega, \omega'\in \Omega_X} 
\int_X |D_yF(\omega)|^2|D_yF(\omega')| \nu(dy) 
} 
\\
&\leq & \frac{2}{3} \sup_{\omega, \omega'\in \Omega_X} 
\int_X |D_yF(\omega)|^3 \nu(dy)+ 
 \frac{1}{3} \sup_{\omega, \omega'\in \Omega_X} 
\int_X |D_yF(\omega')|^3\nu(dy) 
 \leq \alpha_3, 
\end{eqnarray*}
 and similarly for $q=4$: 
\begin{eqnarray*}
\lefteqn{ 
\sup_{\omega, \omega'\in \Omega_X} 
\int_X |D_yF(\omega)|^3|D_yF(\omega')| \nu(dy) 
} 
\\
&\leq & \frac{3}{4} \sup_{\omega, \omega'\in \Omega_X} 
\int_X |D_yF(\omega)|^4 \nu(dy)+\frac{1}{4} \sup_{\omega, \omega'\in \Omega_X} 
\int_X |D_yF(\omega')|^4 \nu(dy) 
 \leq \alpha_4. 
\end{eqnarray*}
 From \eqref{eq:alpha4} we get 
\begin{eqnarray}
\nonumber
h(s) & \leq &s\alpha_2 +\frac{s^2}2\alpha_3+\frac{e^{sK}-1-sK-\frac {s^2}2K^2}{K^3} \alpha_4 
\\ 
\label{eq:h123}
&=& 
 s\left(\alpha_2 -\frac{\alpha_4}{K^2}\right)+\frac{s^2}2\left(\alpha_3-\frac{\alpha_4}{K}\right)+\frac{\alpha_4}{K^3} (e^{sK}-1). 
\end{eqnarray}
Since we assume $\alpha_3\leq 2\alpha_4/K$, the second summand in the right-hand side of \eqref{eq:h123} is bounded by the third one for all $s\geq 0$. We may now end the proof as in Lemma 2 of \cite{HM}: 
\begin{eqnarray*} 
h(s)&\leq& 3\max 
 \left( s\left(\alpha_2 -\frac{\alpha_4}{K^2}\right),\frac{s^2}2\left(\alpha_3-\frac{\alpha_4}{K}\right),\frac{\alpha_4}{K^3} (e^{sK}-1)\right) 
\\
&=& 3\max 
 \left( s\left(\alpha_2 -\frac{\alpha_4}{K^2}\right), \frac{\alpha_4}{K^3} (e^{sK}-1)\right) 
\\
&\leq& \left\{
\begin{array}{ll}
\disp 3s\left(\alpha_2-\frac{\alpha_4}{K^2}\right), & 0 \leq s\leq s_0, 
\\
\disp 3\frac {\alpha_4}{K^3} (e^{sK}-1),& s\geq s_0, 
\end{array}
\right .
\end{eqnarray*}
where $s_0$ is the unique positive solution of \eqref{eq:s0} which is well defined since $K^2\alpha_2/\alpha_4\geq~2$. 
 Hence, for $x_0=3s_0\left(\alpha_2-\alpha_4/K^2\right)$, 
$$
h^{-1}(t)=\left\{
\begin{array}{ll}
\disp\frac t{3(\alpha_2-\alpha_4/K^2)}&\mbox{ for } 0 \leq t\leq x_0, 
\\ 
\disp \frac 1{K} \log\left(1+\frac{K^3}{3\alpha_4}t \right)&\mbox{ for } t\geq x_0, 
\end{array}
\right .
$$ 
 which yields \eqref{eq:bis2-1} and \eqref{eq:bis2-2} from \eqref{eq:dev0}. 
\end{Proof} 

\vskip 5pt 
\noindent
 Lemma~\ref{lemme:bis2} will be used in the proof of Theorem~\ref{theo:BHP-bis2} 
 below to obtain a deviation result under $\alpha$-stable L\'evy measures 
 for all value of $\alpha\in(0,2)$. 
 The following lemma applies only for $\alpha\geq 1$, but will yield a slightly 
 better range condition in Theorem~\ref{theo:BHP-bis}, and is stated 
 without boundedness assumption on $4$th the order moment. 
\begin{lemma} 
\label{lemme:bis}
 Let $F:\Omega^X \longrightarrow \real$ and $\alpha_2,\alpha_3,K>0$, 
 such that $K\alpha_2 \geq 2 \alpha_3$ and 
\begin{description} 
\item{(i)} $\Vert DF \Vert_{L^\infty(\Omega^X, L^2(\nu))}^2 
 \leq \alpha_2 < \infty$, 
\item{(ii)} $\Vert DF \Vert_{L^\infty(\Omega^X, L^3(\nu))}^3
 \leq \alpha_3 < \infty$,
\item{(iii)} $\sup_{y\in X} \vert D_yF (\omega ) \vert\leq K <\infty$, 
 $P(d\omega ) $-a.s. 
\end{description} 
 Denote by $s_0$ the unique solution of 
$$\disp \frac{e^{sK}-1}{sK} = 
 K \frac{\alpha_2}{\alpha_3}-1 
. 
$$ 
 Let also $x_0= 2s_0(\alpha_2-\alpha_3/K)$. 
 Then 
\begin{equation}
\label{eq:bis1}
P(F-E[F] \geq x)\leq \exp 
 \left(-\frac{x^2}{4(\alpha_2-\alpha_3/K)}\right)
, \qquad 0 \leq x\leq x_0, 
\end{equation}
 and 
\begin{equation}
\label{eq:bis2}
P(F-E[F] \geq x)\leq K_0 \exp 
 \left(
\frac{x}{K}
-\left(\frac{x}{K}+2\frac{\alpha_3}{K^3}\right)\log\left(1+\frac{K^2x}{2\alpha_3}\right)\right), 
 \qquad x\geq x_0 , 
\end{equation}
 with 
\begin{equation}
\label{eq:K0}
K_0=
\exp\left(-x_0/K+\left(\frac{x_0}K+\frac{2\alpha_3}{K^3}\right)\log\left(1+\frac{K^2x_0}{2\alpha_3}\right)-\frac{x_0^2}{4(\alpha_2-\alpha_3/K)}\right). 
\end{equation} 
\end{lemma}
\begin{Proof}
As in the proof of Lemma \ref{lemme:bis2}, apply Proposition~\ref{thd3} with $h$ given in \eqref{eq:h1} and bounded by 
$$ 
h(s) \leq s\alpha_2 +\frac{e^{sK}-1-sK}{K^2} \alpha_3. 
$$ 
using  
$$ 
e^{su}-1\leq su+\frac{e^{sK}-1-sK}{K^2} u^2, 
 \qquad u\in [0,K], 
$$
instead of \eqref{eq:h-bound2}. 
 We get 
$$
h(s) \leq2\max\left( 
 s \left(\alpha_2+\frac{\alpha_3}{K}\right) ,
\frac{e^{sK}-1}{K^2}\alpha_3\right)\leq \left\{ 
\begin{array}{ccc}
2s (\alpha_2 - \alpha_3 / K ), &\mbox{  }& s\leq s_0\\
2(e^{sK}-1)\alpha_3 / K^2, &\mbox{  }& s\geq s_0
\end{array} 
\right. 
$$
 which allows to conclude as in the proof of Lemma \ref{lemme:bis2}.
\end{Proof}

\section{The case of stable L\'evy measures}
\label{sec:stable} 
 Let $0<\alpha<2$, $X=\real^n$ and the stable L\'evy measure given by 
\begin{equation} 
\label{eq:Levy-stable}
\nu(B)= \int_{S^{n-1}}\sigma(d\xi)\int_0^\infty \ind_B(r\xi) r^{-1-\alpha} dr, 
 \qquad 
 B\in{\cal B}(\real^n) 
, 
\end{equation} 
 where $\sigma$ is a finite positive measure on $S^{n-1}$, 
 the unit sphere of $\real^n$, called the spherical component of $\nu$. 
 We have 
\begin{eqnarray}
\nonumber
P(\{ \omega \in \Omega^X \ : \  
 \om^c_R\not =\emptyset \} 
 )& 
 = & 
 1 - P(\{ \omega \in \Omega^X \ : \  
 \om^c_R =\emptyset \} 
 )
\\ 
\nonumber
 & = & 
 1-\exp\left(-\int_{\{\vert y \vert_2 >R\}}\nu (dy )\right)\\
\nonumber
&=&1 -\exp\left(\int_{S^{n-1}}\sigma(d\xi)\int_{\{\vert r\xi\vert_2 \geq R\}}\frac{dr}{r^{1+\alpha}}\right)\\
\nonumber
&=&1-\exp\left(-\frac{\sigma(S^{n-1})}\alpha R^{-\alpha}\right)\\
\label{eq:gamma-stable}
&\leq&\frac{\sigma(S^{n-1})}{\alpha R^{\alpha}}. 
\end{eqnarray}
 Thus we can take 
$$ 
\gamma(R)=\frac{\sigma(S^{n-1})}{\alpha R^\alpha}, 
 \qquad R>0. 
$$ 
 in Theorem~\ref{theo:BHP}. 
 Let $f:\real^n\to \real$ be $\ell^2$-Lipschitz$(\constante)$. 
 In case $F$ is a stable random variable represented by a single Poisson stochastic 
 integral of the form \eqref{eq:vaX2}, we have from \eqref{hjl}: 
\begin{eqnarray*}
\|D f(F) \|_{L^\infty(\Omega^X, L^2(\nu_R))}^2& \leq &
\int_{\{\vert y \vert_2 \leq R\}} \constante^2\vert y \vert_2 ^2\: \nu (dy )
\\ 
 &=& \constante^2\int_{S^{n-1}}\sigma(d\xi)\int_{\{\vert r\xi\vert_2 \geq R\}} 
 r^{1-\alpha} dr
\\ 
&\leq& 
 \frac{\constante^2\sigma(S^{n-1})}{2-\alpha} R^{2-\alpha} 
\\ 
& \leq & \frac{2\constante^2}{2-\alpha} R^2 \gamma (R), 
\end{eqnarray*} 
 hence Theorem~1 of \cite{HM} is recovered 
 taking $\beta (r) = \constante r$ and $C= 2/(2-\alpha)$ in Theorem \ref{theo:BHP}, i.e. 
\begin{equation} 
\label{cstexpl} 
 P(f(F)-m(f(F))\geq x)\leq 
 \left(1 + \frac{2e}{2-\alpha} \right) 
 \frac{\sigma(S^{n-1})}\alpha \left(  \frac{x}{4\constante}\right)^{-\alpha}, 
\end{equation} 
 for all $x$ such that 
$$ 
x\geq 2 \constante \gamma^{-1}(2-\alpha ) 
 \geq 2 \constante \gamma^{-1}\left(\frac 1{2(1+2e/(2-\alpha) )}\right), 
$$ 
 where $F$ is a stable random variable with parameter $\alpha$. 
 The constant in front of $x^{-\alpha}$ in \eqref{cstexpl} explodes as 
 $\alpha$ goes to $0$ or to $2$. 
 In fact, 
 as noted in \cite{HM}, the dependency in $\alpha^{-1}$ of the constant 
 is sharp as $\alpha$ goes to $0$ (take for example 
 a symmetric $\alpha$-stable random variable). 
 This explosion does not occur however when $\alpha$ goes to $2$, 
 and the aim of the next result is to provide a deviation bound with 
 such a non-exploding constant, for general random variables 
 on Poisson space under $\alpha$-stable intensity measures. 
 The proof relies on Lemma~\ref{lemme:bis2}, 
 and in the particular case of stable random variables, this result also 
 improves Theorem~2 of \cite{HM} by allowing $\alpha$ to be arbitrary. 
\begin{theorem} 
\label{theo:BHP-bis2}
 Let $\alpha \in (0,2)$ and $F:\Omega^X \to \real$ such that 
$$ 
 \vert 
 D_yF (\omega ) 
 \vert  
 \leq \constante \vert y \vert_X, \quad P(d\omega )\otimes \nu(dy)\mbox{-a.e.},  
$$ 
 with $\constante >0$. 
 Then we have 
\begin{equation}
\label{eq:but2-2}
P(F-m(F)\geq x) 
 \leq \sigma(S^{n-1}) 
 \left(\frac 32 e^2+\frac{1}{\alpha}\right) 
 \frac{(4\constante)^{\alpha}}{x^\alpha} 
, 
\end{equation}
 for all 
\begin{equation} 
\label{eq:range-f-2}
x\geq 4\constante 
 \sigma(S^{n-1})^{1/\alpha} 
 \left(
 \left( 
 \frac 32  \left(1+\frac{4}{2-\alpha} \log \frac 2{2-\alpha}\right)\log\left(1+\frac{8}{2-\alpha}\log \frac 2{2-\alpha}\right)
 \right) 
 \vee \frac 4\alpha 
 \vee (6e^2) 
 \right)^{1/\alpha}. 
\end{equation}
\end{theorem}

\begin{Proof} 
 Using the notation of the proof of Theorem~\ref{theo:BHP} we have 
$$ 
D_yg(F_R) (\omega ) \leq \vert D_y F(\omega_R ) \vert \leq \constante \vert y
\vert_X, \quad P(d\omega)\otimes \nu(dy) \mbox{ a.e.} 
, 
$$ 
 where $g(x) = (x-m(F_R))^+\wedge r$. Thus 
\begin{align*}
&\sup_{y \in B_X(0,R)} 
 D_yg(F_R)\leq \constante R, \quad P-a.s. 
, 
\\ 
& \|Dg(F_R)\|_{L^\infty(\Omega^X , L^2(\nu_R))}^2\leq 
 \frac{\constante^2 \sigma(S^{n-1})}{2-\alpha} R^{2-\alpha} 
,
\\
&
\|Dg(F_R)\|_{L^\infty(\Omega^X , L^3(\nu_R))}^3 \leq 
 \frac{\constante^3 \sigma(S^{n-1})}{3-\alpha} R^{3-\alpha}
,
\end{align*} 
 and 
$$ 
\|Dg(F_R)\|_{L^\infty(\Omega^X , L^4(\nu_R))}^4 \leq 
 \frac{\constante^4 \sigma(S^{n-1})}{4-\alpha} R^{4-\alpha}. 
$$ 
 We now apply Lemma~\ref{lemme:bis2} to $\nu_R$ and $F_R$ with 
$$
K=\constante R, \quad 
 \alpha_2 = \frac{\constante^2 \sigma(S^{n-1})}{2-\alpha} R^{2-\alpha}, 
 \quad 
 \alpha_3 = \frac{\constante^3 \sigma(S^{n-1})}{3-\alpha} R^{3-\alpha}, 
\quad 
 \alpha_4 = \frac{\constante^4 \sigma(S^{n-1})}{4-\alpha} R^{4-\alpha} 
. 
$$ 
 Using \eqref{and}, equation \eqref{eq:s0} reads 
$$
\varphi(scR):=e^{s\constante R} 
 - \frac{2s\constante R}{2-\alpha} 
 -1 
 = 0. 
$$
Since for all $\alpha\in (0,2)$, $\varphi(\log \frac 2{2-\alpha})\leq 0$ and $\varphi(2\log \frac 2{2-\alpha})\geq 0$ we have 
$$
\log\frac 2{2-\alpha}\leq s_0\constante R\leq 2\log \frac 2{2-\alpha}, 
$$
so that for $\disp x_0 =3s_0\left(\alpha_2 -\frac{\alpha_4}{K^2}\right)=\frac{6\constante^2\sigma(S^{n-1})R^{2-\alpha}}{(2-\alpha)(4-\alpha)}s_0$, we have  
\begin{eqnarray} 
\label{eq:x0-2}
\constante\sigma(S^{n-1})R^{1-\alpha}\frac{3}{2(2-\alpha)} \log\frac 2{2-\alpha}&\leq x_0 \leq& \constante\sigma(S^{n-1})R^{1-\alpha}\frac{6}{2-\alpha} \log\frac 2{2-\alpha}.
\end{eqnarray} 
For 
\begin{eqnarray*}
r&\geq& 2x_0 =\frac{12\constante^2\sigma(S^{n-1})R^{2-\alpha}}{(2-\alpha)(4-\alpha)}s_0,
\end{eqnarray*} 
 we get from Lemma~\ref{lemme:bis2}: 
\begin{eqnarray} 
\label{eq:tech1.0-2} 
\lefteqn{ 
P(F_R - m(F_R)\geq r ) 
\leq P(g(F_R)-E[g(F_R)]\geq r/2) 
} 
\\ 
\nonumber 
& \leq & 
 K_0 \exp 
 \left(
\frac{r}{2\constante R}
-\left(\frac{r}{2\constante R}+3\frac{\sigma(S^{n-1})}{(4-\alpha)R^\alpha}\right)\log\left(1+\frac{(4-\alpha)r}{6\sigma(S^{n-1})\constante R^{1-\alpha}}\right)\right)
\end{eqnarray} 
 with from \eqref{eq:K0-2} and \eqref{eq:x0-2}: 
\begin{align*}
&K_0\\
&=\exp\left(-\frac{x_0}{\constante R}+\left(\frac{x_0}{\constante R}+3\frac{\sigma(S^{n-1})}{(4-\alpha) R^\alpha}\right)\log\left(1+\frac{(4-\alpha)x_0}{3\sigma(S^{n-1})\constante R^{1-\alpha}}\right)-\frac {(2-\alpha)(4-\alpha)x_0^2}{12\sigma(S^{n-1})\constante ^2 R^{2-\alpha}}\right)\\
&\leq\exp\left(\frac{3\sigma(S^{n-1})}{2R^\alpha}\left(1+\frac{4}{2-\alpha} \log \frac 2{2-\alpha}\right)\log\left(1+\frac{8}{2-\alpha}\log \frac 2{2-\alpha}\right)\right).
\end{align*}
Hence under the condition 
\begin{equation}
\label{eq:cond0-2}
 \sigma(S^{n-1})R^{-\alpha} 
 \leq \frac{2}{3\left(1+\frac{4}{2-\alpha} \log \frac 2{2-\alpha}\right)\log\left(1+\frac{8}{2-\alpha}\log \frac 2{2-\alpha}\right)}
\end{equation} 
 we get $K_0\leq e$ and 
\begin{eqnarray*}
2x_0 & \leq & 12\constante\sigma(S^{n-1})R^{1-\alpha}\frac{1}{2-\alpha} \log\frac 2{2-\alpha}\\ 
&\leq& 
 4\constante R 
 \frac{\frac{2}{2-\alpha} \log\frac 2{2-\alpha}}{\left(1+\frac{4}{2-\alpha} \log \frac 2{2-\alpha}\right)\log\left(1+\frac{8}{2-\alpha}\log \frac 2{2-\alpha}\right)}
\\
&\leq&\constante R 
, 
\end{eqnarray*}
 i.e. $r\geq 2x_0$ with $r= 2 \constante R$. 
 Then from \eqref{eq:tech1.0-2} and $K_0\leq e$ we get 
\begin{eqnarray} 
\nonumber 
P\big(F_R-m(F_R)\geq 2 \constante R \big) 
&\leq &\exp 
 \left(2 -\left(1+\frac{3\sigma(S^{n-1})}{(4-\alpha)R^\alpha}\right)\log\left(1+\frac{(4-\alpha)R^{\alpha}}{3\sigma(S^{n-1})}\right)\right)
\\
\nonumber 
&\leq & e^2 \left(1+\frac{(4-\alpha)R^{\alpha}}{3\sigma(S^{n-1})}\right)^{-1} 
\\
\nonumber 
&\leq & \frac{3 e^2 \sigma(S^{n-1})}{(4-\alpha)R^\alpha}
\\ 
\label{eq:tech1-1-2} 
&\leq  & \frac{3e^2\sigma (S^{n-1})}{2R^\alpha} 
\\ 
\label{*1.3-2}
& = & \frac 32 e^2\alpha\gamma (R), 
\end{eqnarray} 
 as long as \eqref{eq:cond0-2} holds. 
 In order to control $P(F_R-m(F)\geq x)$ from (\ref{eq:tech1-1-2}), 
 we need to control $m(F_R)-m(F)$.
 For this we apply Lemma~\ref{lemma:BHP} 
 with $\tilde \beta(R) =2\constante R$, $\tilde{\gamma} (R) = \frac 32 e^2\alpha \gamma (R)$, 
$$
R_0=\left(\frac 32 \sigma(S^{n-1})
 \left(1+\frac{4}{2-\alpha} \log \frac 2{2-\alpha}\right)\log\left(1+\frac{8}{2-\alpha}\log \frac 2{2-\alpha}\right)\right)^{1/\alpha}
$$
 and \eqref{*1.3-2}. This yields, 
 with $x = 4 \constante R$: 
$$ 
m(F_R)-m(F)\leq x/2, 
$$ 
 and 
$$ 
 P\big(F_R-m(F)\geq x \big) 
 \leq P\big(F_R-m(F_R)\geq x/2 \big) 
 \leq 
\frac 32 e^2 \sigma(S^{n-1}) 
 \left(\frac{4\constante}{x}\right)^\alpha, 
$$ 
 provided 
\begin{equation} 
\label{eq:cond0.1-2}
R\geq \max \left(R_0,  
 \gamma^{-1}(\delta ) 
 , 
 \gamma^{-1} \left( 
 \frac{2}{3\alpha e^2} 
 \left( \frac{1}{2} - \delta 
 \right) \right) 
 \right) 
, 
\end{equation} 
 for any given $\delta \in (0,1/2)$. 
 When $x=4\constante R$, this estimate together with 
$$
P(\{ \omega \in \Omega^X \ : \  
 \omega^c_R\not =\emptyset \} )\leq 
 \gamma (R) = \frac{\sigma(S^{n-1}) R^{-\alpha}}{\alpha} 
, 
$$ 
 gives, using \eqref{eq:tech2}: 
\begin{eqnarray}
\nonumber
P(F-m(F)\geq x)&=&P(F-m(F)\geq x,\; \omega^c_R=\emptyset)+P(F-m(F)\geq x,\; \omega^c_R\not =\emptyset) 
\\ 
\label{eq:fin1-2} 
 & \leq & \sigma(S^{n-1}) 
 \left(\frac 1\alpha+\frac 32 e^2\right) 
 \left(\frac{x}{4\constante}\right)^{-\alpha}  ,
\end{eqnarray} 
 as long as (\ref{eq:cond0-2}) and \eqref{eq:cond0.1-2} hold. 
 Now, conditions \eqref{eq:cond0-2} and \eqref{eq:cond0.1-2} 
 can be rewritten in terms of $x$ as 
$$
 x\geq 4\constante\left(
\frac 32 \sigma(S^{n-1}) \left(1+\frac{4}{2-\alpha} \log \frac 2{2-\alpha}\right)\log\left(1+\frac{8}{2-\alpha}\log \frac 2{2-\alpha}\right)\right)^{1/\alpha} 
$$
 and 
$$ 
x\geq 4 \constante 
 \max 
 \left( \left( 
 \frac{\sigma(S^{n-1})}{\alpha \delta} 
 \right)^{1/\alpha}, 
 \left(\frac{3\sigma(S^{n-1})e^2}{2(1/2-\delta)}\right)^{1/\alpha} 
 \right),
$$ 
 When e.g. $\delta = 1/4$, the range of \eqref{eq:fin1-2} can be written  
$$ 
x\geq 4\constante 
 \sigma(S^{n-1})^{1/\alpha} 
 \left(
 \left( 
 \frac 32  \left(1+\frac{4}{2-\alpha} \log \frac 2{2-\alpha}\right)\log\left(1+\frac{8}{2-\alpha}\log \frac 2{2-\alpha}\right)
 \right) 
 \vee \frac 4\alpha 
 \vee (6e^2 ) 
 \right)^{1/\alpha}. 
$$
\end{Proof} 
 Using Lemma \ref{lemme:bis} instead of Lemma \ref{lemme:bis2}, we can state a 
 similar deviation result 
 under a slight better range condition on $x$, in case $\alpha\in[1,2)$. 
\begin{theorem}
\label{theo:BHP-bis}
 Assume that $\alpha \geq 1$ and let $F:\Omega^X \to \real$ such that 
$$ 
 \vert 
 D_yF (\omega ) 
 \vert  
 \leq \constante \vert y \vert_X, \quad P(d\omega )\otimes \nu(dy)\mbox{-a.e.},  
$$ 
 with $\constante >0$. 
 Then we have 
\begin{equation}
\label{eq:but2}
P(F-m(F)\geq x) 
 \leq \sigma(S^{n-1}) 
 \left(1+\frac{e^2}{2}\right) 
 \frac{(4\constante)^{\alpha}}{x^\alpha} 
, 
\end{equation}
 for all 
\begin{equation} 
\label{eq:range-f}
 x\geq 4 \constante 
 \sigma(S^{n-1})^{1/\alpha} 
 \left( 
 \left( 
 \left( 
 1+\frac 2{2-\alpha} \log\frac 1{2-\alpha}\right) 
 \log \left(1+\frac 4{2-\alpha}\log \frac 1{2-\alpha} \right) 
 \right)  
 \vee ( 4e^2 ) 
 \right)^{1/\alpha}. 
\end{equation}
\end{theorem} 

\begin{Proof} 
 We sketch the modifications of the proof, 
 following the argument of Theorem \ref{theo:BHP-bis2} and 
 applying Lemma~\ref{lemme:bis} instead of Lemma \ref{lemme:bis2} to $\nu_R$ and $F_R$, with 
$$
K=\constante R, \qquad 
 \alpha_2 = \frac{\constante^2 \sigma(S^{n-1})}{2-\alpha} R^{2-\alpha}, 
 \qquad 
 \alpha_3 = \frac{\constante^3 \sigma(S^{n-1})}{3-\alpha} R^{3-\alpha}.
$$ 
Under the condition 
\begin{equation}
\label{eq:cond0}
 \sigma(S^{n-1})R^{-\alpha} 
 \leq \frac{1}{2(1+\frac 2{2-\alpha} \log \frac 1{2-\alpha} ) \log \big(1+\frac 4{2-\alpha}\log \frac 1{2-\alpha}\big)}
\end{equation} 
and since 
$$
\frac{2\alpha_3}{\constante^2R^2} \frac 1{2-\alpha}\log \frac 1{2-\alpha}\leq x_0
\leq \frac{4\alpha_3}{\constante^2R^2} \frac 1{2-\alpha}\log \frac 1{2-\alpha}, 
$$ 
we have $K_0\leq e$ and also $x_0\leq \constante R$.
Using \eqref{and} with $r=2\constante R\geq 2x_0$, we get applying Lemma \ref{lemme:bis}: 
\begin{eqnarray} 
\nonumber 
P\big(F_R-m(F_R)\geq 2 \constante R \big) 
&\leq & \exp\left( 2 -\left( 
 1 +2\frac{\alpha_3}{\constante^3R^3}\right)\log 
 \left( 
 1+\frac{2 \constante^3R^3}{\alpha_3} 
 \right) 
 \right) 
\\
\nonumber 
&\leq & e^2 \left( 
 1+\frac{2 \constante^3R^3}{\alpha_3}\right)^{-1} 
\\
\nonumber 
&\leq & \frac{e^2 \alpha_3}{2 \constante^3R^3}
\\ 
\label{eq:tech1-1} 
&\leq  & \frac{e^2\sigma (S^{n-1})}{2R^\alpha} 
\\ 
\label{*1.3}
& = & 
 \frac{e^2\alpha}{2} \gamma (R), 
\end{eqnarray} 
 as long as \eqref{eq:cond0} holds. 
 Finally, applying Lemma~\ref{lemma:BHP} 
 with $\tilde \beta(R) =2\constante R$, $\tilde{\gamma} (R) = e^2\alpha \gamma (R) /2$ and condition \eqref{*1.1} given by \eqref{*1.3}, with $x = 4 \constante R$, derive $m(F_R)-m(F)\leq x/2$, and 
$$ 
 P\big(F_R-m(F)\geq x \big) 
 \leq P\big(F_R-m(F_R)\geq x/2 \big) 
 \leq 
 \frac{e^2 \sigma(S^{n-1}) }{2} 
 \left(\frac{4\constante}{x}\right)^\alpha, 
$$ 
 provided moreover for any, $0<\delta<1/2$,
\begin{equation} 
\label{eq:cond0.1} 
R\geq \max \left( 
 \gamma^{-1}(\delta ) 
 , 
 \gamma^{-1} \left( 
 \frac{2}{\alpha e^2} 
 \left( \frac{1}{2} - \delta 
 \right) \right) 
 \right).
\end{equation} 
 With $x=4\constante R$, this estimate together with \eqref{eq:gamma-stable} 
 gives, using $\alpha\geq 1$ and \eqref{eq:tech2}: 
\begin{eqnarray}
\nonumber
P(F-m(F)\geq x)&=&P(F-m(F)\geq x,\; \omega^c_R=\emptyset)+P(F-m(F)\geq x,\; \omega^c_R\not =\emptyset) 
\\ 
\label{eq:fin1}
 & \leq & \sigma(S^{n-1}) 
 \left(1+\frac{e^2}{2}\right) 
 \left(\frac{x}{4\constante}\right)^{-\alpha},
\end{eqnarray} 
 as long as (\ref{eq:cond0}) and \eqref{eq:cond0.1} hold. 
 Now, conditions \eqref{eq:cond0} and \eqref{eq:cond0.1} can be rewritten in terms of $x$ as \eqref{eq:range-f} with e.g. $\delta=1/4$. 
\end{Proof} 
 Finally, we extend a recent result of \cite{marchal} to Poisson functionals 
 under stable intensity measures. 
\begin{theorem} 
\label{theo:BHP-ter}
Let $F:\Omega^X \to \real$ such that for some $\constante>0$, 
$$ 
 \vert 
 D_yF (\omega ) 
 \vert  
 \leq \constante \vert y \vert_X, \quad P(d\omega )\otimes \nu(dy)\mbox{-a.e.}  
$$
1)  Let $\varepsilon >0$, then if $\alpha$ is sufficiently close to $2$, 
\begin{equation}
\label{eq:i}
P(F-m(F)\geq x)\leq 
 (\varepsilon + \sqrt e ) \exp \left( -\frac{(2-\alpha) x^\alpha}{2(4\constante)^\alpha\sigma(S^{n-1})}\right), 
\end{equation}
 provided 
\begin{equation}
\label{eq:C0}
\frac{2\sigma(S^{n-1}) (4\constante)^\alpha}{2-\alpha} \log(4(1 + \sqrt e)) 
 \leq x^\alpha \leq 
 \frac{\sigma(S^{n-1}) (4\constante)^\alpha}{2(2-\alpha)} 
 \frac{\log(\frac 1{2-\alpha})}{3-\alpha}. 
\end{equation} 
2) 
 Let $b>3$, $\varepsilon > 0$, 
 and $x=4b\constante \sigma(S^{n-1})\frac 1{2-\alpha} 
 \log \frac1{2-\alpha}$. 
 For $\alpha$ close enough to $2$ we have 
\begin{equation}
\label{eq:ii}
P(F-m(F)\geq x)\leq \frac{(4\constante)^\alpha\sigma(S^{n-1})}{x^\alpha}\left(\frac 1\alpha+(2+\varepsilon) \exp\left( 
 \frac{(2+\varepsilon)(4\constante)^\alpha\sigma(S^{n-1})g(2-\alpha)}{x^\alpha} 
 \right) 
 \right) 
\end{equation}
where $\disp g(x)=\left(\frac 1x \log \frac 1x\right) \log \left(\frac 1x\log \frac 1x\right)$.
\end{theorem} 
\begin{Proof} 
 We follow \cite{marchal} as in the proofs of Theorem
 \ref{theo:BHP} and Proposition \ref{theo:BHP-bis} above. 
First, using the same notation as before, we have: 
\begin{eqnarray}
\label{eq:F_R-g}
P(F_R-m(F_R)\geq r)&\leq& P(g(F_R)-E[g(F_R)]\geq r/2)\\
\label{eq:mar-maj1}
&\leq&  \exp\left(-\int_0^{r/2} h_R^{-1}(s) ds \right), \med 0<x<h_R(t_0^-), 
\end{eqnarray}
 with 
\begin{eqnarray*} 
h_R(s)&\leq& \left(\alpha_2-\frac{\alpha_3} K\right)s 
 +\frac{\alpha_3}{K^2} \left(e^{sK}-1\right) 
\\
&=&\frac{\constante^2\sigma(S^{n-1}) R^{2-\alpha}}{(2-\alpha)(3-\alpha)}s+(e^{s\constante R}-1)\frac{\constante \sigma(S^{n-1}) R^{1-\alpha}}{3-\alpha}
\end{eqnarray*}
since again
$$
\alpha_2 = \frac{\constante^2 \sigma(S^{n-1})}{2-\alpha} R^{2-\alpha}, 
 \qquad 
 \alpha_3 = \frac{\constante^3 \sigma(S^{n-1}) }{3-\alpha}R^{3-\alpha}
\qquad K=\constante R,
$$ 
where \eqref{eq:F_R-g} above comes as in \eqref{and} in the proof of Theorem
\ref{theo:BHP} and \eqref{eq:mar-maj1} comes from the proofs of
Lemma~\ref{lemme:bis} and Proposition~\ref{thd3}. 
 Following \cite{marchal}, for  $\delta, s,R$ satisfying 
\begin{equation} 
\label{eq:marchal1}
\frac{e^{s\constante R}-1}{\constante sR}\leq \frac\delta{2-\alpha}
\end{equation}
we have 
$$
h_R(s)\leq (1+\delta) \frac{\constante^2\sigma(S^{n-1}) R^{2-\alpha}}{(2-\alpha)(3-\alpha)} s
$$
and 
\begin{equation}
\label{eq:mar2}
\int_0^yh_R^{-1}(t) dt\geq \frac{(3-\alpha)(2-\alpha)y^2}{2(1+\delta)\constante^2\sigma(S^{n-1}) R^{2-\alpha}}
\end{equation}
 for all $y$ such that 
$$\frac{(3-\alpha)(2-\alpha)y}{(1+\delta)\constante^2\sigma(S^{n-1}) R^{2-\alpha}} \leq s, 
$$ 
 where $s$ satisfies \eqref{eq:marchal1}.
 Taking for some $A>0$, $\disp R^\alpha=
\frac{A\sig}{(2-\alpha)(3-\alpha)}$ and $y=R\constante$, since 
 $s\constante R=A/(1+\delta)$, \eqref{eq:marchal1} can be rewritten 
 as 
$$
\disp (1+\delta)\frac{e^{\frac A{1+\delta}}-1}{A}\leq \frac{\delta}{2-\alpha}
$$
which is satisfied whenever 
$$
(1+\delta) \frac{e^A}A\leq \frac\delta{2-\alpha}. 
$$
Choosing $\disp\delta=\frac{e^A(2-\alpha)}{A-e^A(2-\alpha)}$ which is positive for 
 $0<a<A< - \log (2-\alpha )$ when $\alpha$ is close enough to $2$, we derive 
 from \eqref{eq:mar2} for $a<A< - \log ( 2-\alpha )$
\begin{equation} 
\label{eq:mar3}
\exp\left( 
 - \int_0^{\constante R} h_R^{-1}(t) dt\right) 
 \leq e^{-\frac A 2}\exp\left(\frac{e^A(2-\alpha)}{2}\right).
\end{equation}
But since 
$$
\lim_{\alpha\to 2^-} \sup_{a<A<-\log (2-\alpha)}e^{-\frac A 2}\exp\left(\frac{e^A(2-\alpha)}{2}\right) e^{\frac A{2(3-\alpha)}}=\sqrt e,
$$
for any $\varepsilon\in (0,1)$ and $\alpha$ close to $2$, from \eqref{eq:mar-maj1} with $r=2\constante R$, 
\begin{eqnarray}
\nonumber
P(F_R-m(F_R)\geq 2\constante R)&\leq& \exp\left( 
 - \int_0^{\constante R} h_R^{-1}(t) dt\right) 
 \leq \left(\sqrt e+\frac \varepsilon 2\right)e^{-\frac A {2(3-\alpha)}}\\
\label{eq:mar4} 
&\leq& \left(\sqrt e+\frac\varepsilon 2\right) 
 \exp\left(-\frac{(2-\alpha)}{2\sig}R^\alpha\right) 
\end{eqnarray}
 for 
\begin{equation}
\label{eq:range1}
\frac{a\sigma(S^{n-1})}{(2-\alpha)(3-\alpha)}<R^\alpha<\frac{\sigma(S^{n-1})\log \frac 1{2-\alpha}}{2(2-\alpha)(3-\alpha)}.
\end{equation}
 Next, control $m(F)-m(F_R)$ using Lemma~\ref{lemma:BHP} 
 with $\tilde \beta(R)=2\constante R$, 
$$
\tilde{\gamma} (R) = \left( 1 + \sqrt e \right) 
 \exp\left(-\frac{(2-\alpha)}{2\sig}R^\alpha\right) 
$$ 
 and condition \eqref{*1.1} given by \eqref{eq:mar4} (with $\varepsilon\leq 2$). This yields 
\begin{equation}
\label{eq:MFR}
 m(F_R)-m(F) \leq 2\constante R 
, 
\end{equation} 
 provided \eqref{ccond}, rewritten as 
\begin{equation}
\label{eq:range2}
R^\alpha\geq \max\left(\frac{\sig}{\alpha\delta}, -\frac{2\sig}{2-\alpha}\log\frac{1/2-\delta}{ 1 + \sqrt e}\right)
, 
\end{equation} 
 and \eqref{eq:range1} above still hold. 
 Equations \eqref{eq:mar4} and \eqref{eq:MFR} yield
\begin{eqnarray}
\nonumber 
 P\big(F_R-m(F)\geq 4\constante R )\big) & \leq& P\big(F_R-m(F_R)\geq 2\constante R \big)\\
\label{eq:mar5}
&\leq&  \left( 
 \sqrt e+\frac\varepsilon 2 
 \right) 
 \exp\left(-\frac{(2-\alpha)}{2\sig} R^\alpha\right)
\end{eqnarray}
 provided \eqref{eq:range1} and \eqref{eq:range2} hold. 
 Next, when \eqref{eq:range1} holds, \eqref{eq:gamma-stable} gives for $\alpha$ close 
 enough to $2$: 
\begin{equation}
\label{eq:mar6}
P(\omega_R^c\not=\emptyset)\leq \frac{\sig}{\alpha R^\alpha}\leq \frac
\varepsilon 2 \exp\left( 
 -\frac{(2-\alpha)R^\alpha}{2\sig}\right) 
. 
\end{equation}
 Finally, \eqref{eq:tech2} together with \eqref{eq:mar5} and \eqref{eq:mar6} yields with $x=4\constante R$,  
$$
P(F-E[F_R]\geq x)\leq 
(\sqrt e+\varepsilon)\exp\left( 
 -\frac{(2-\alpha)x^\alpha}{2(4\constante)^\alpha\sig}\right) 
$$
 as long as 
\begin{equation}
\label{eq:C1}
a\frac{\sigma(S^{n-1}) (4\constante)^\alpha}{(2-\alpha)(3-\alpha)}\leq x^\alpha \leq 
 (4\constante)^\alpha \sig\frac{\log (1/(2-\alpha))}{2(2-\alpha)(3-\alpha)}
, 
\end{equation}
 and 
\begin{equation}
\label{eq:C2}
x^\alpha\geq (4\constante)^\alpha\max\left(\frac{\sig}{\alpha\delta}, -\frac{2\sig}{2-\alpha}\log\frac{1/2 -\delta}{1 + \sqrt e}\right) 
, 
\end{equation}
 for any $0 < \delta < 1/2$. 
 Taking $\delta=1/4$, conditions \eqref{eq:C1} and \eqref{eq:C2} can be rewritten  
 as 
$$ 
\frac{2\sigma(S^{n-1}) (4\constante)^\alpha}{2-\alpha} \log(4(1 + \sqrt e )) 
 \leq x^\alpha \leq 
 \frac{\sigma(S^{n-1}) (4\constante)^\alpha}{2(2-\alpha)} 
 \frac{\log(1/(2-\alpha))}{3-\alpha} 
, 
$$
 which yields \eqref{eq:ii}. 

\vskip 10pt
 
We now deal with the second part of Theorem \ref{theo:BHP-ter}, still following \cite{marchal}. Take for some $b>0$, $R^\alpha=\disp\frac{b\sigma(S^{n-1})\log (1/(2-\alpha))}{2-\alpha}$ 
 and let $A>0$. 
 For $\alpha$ close to $2$ and $s\constante R\geq \log (1/(2-\alpha)) 
 +\log \log (1/(2-\alpha)) +A$ we have
$$
\frac{e^{s\constante R}-1}{s\constante R}\geq \frac{1}{(2-\alpha)(e^{-A}+\varepsilon)}
, 
$$ 
hence 
\begin{equation}
\label{eq:mar8}
h_R ^{-1}(u)\geq \frac 1{\constante R}\log\left(1+\frac{(3-\alpha)u}{(e^{-A}+1+\varepsilon)\constante\sigma(S^{n-1})R^{1-\alpha}}\right)
\end{equation}
whenever 
$$u>u_1=\disp \frac{(1+e^{-A}+\varepsilon)\constante\sigma(S^{n-1})e^A 
 \log (1/(2-\alpha)) 
 }{(2-\alpha)(3-\alpha) R^{\alpha-1}}. 
$$  
 For $A>0$ small enough and $b>2$, we have $\constante R>u_1$, and integrating \eqref{eq:mar8} over $[u_1, \constante R]$, we obtain 
\begin{equation}
\label{eq:mar9} 
\int_{u_1}^{\constante R} h_R^{-1}(t) dt \geq \left( 
 \left(1+\frac 1{\theta \constante R}\right)\log (1+\theta \constante R)-1\right) 
 -\left( 
 \left(\frac{u_1}{\constante R}+\frac 1{\theta \constante R}\right)\log (1+\theta u_1)-\frac{u_1}{\constante R}\right) 
\end{equation}
 with 
 $\theta =\disp \frac{(3-\alpha)R^{\alpha-1}}{(1+e^{-A}+\varepsilon)\constante\sigma(S^{n-1})}$.
 For $\alpha$ close to $2$ and $A,\varepsilon >0$ small enough, using  \eqref{eq:F_R-g}, \eqref{eq:mar-maj1}, \eqref{eq:mar9}, we derive 
\begin{align}
\nonumber 
&P(F_R-m(F_R)\geq 2\constante R)\\ 
\nonumber 
&\leq \exp\left(-\int_0^{\constante R} h_R^{-1}(t) dt\right)
\leq \exp\left(-\int_{u_1}^{\constante R} h_R^{-1}(t) dt\right)\\
\nonumber 
&\leq
\frac{(2+\varepsilon)\sigma(S^{n-1})}{R^\alpha} \exp\left(\frac{(2+\varepsilon)\constante\sigma(S^{n-1})\left(\frac 1{2-\alpha} 
 \log (1/(2-\alpha)) 
 \right)\log \left(\frac 1{2-\alpha} 
 \log (1/(2-\alpha)) 
 \right)}{R^\alpha}\right)\\
\label{eq:mar10}
&\leq
\frac{(2+\varepsilon)\sigma(S^{n-1})}{R^\alpha} \exp\left(\frac{(2+\varepsilon)\sigma(S^{n-1})g(2-\alpha)}{R^\alpha}\right).
\end{align}
with $g(x)=\left(\frac 1x \log \frac 1x\right) \log \left(\frac 1x\log \frac 1x\right)$ and $\varepsilon$ some (new) positive constant.
 It is easy now to control $m(F)-m(F_R)$ using 
 once more Lemma~\ref{lemma:BHP} with $\tilde \beta(R)=2\constante R$, $\tilde\gamma (x) = \alpha^{-1}
 \sigma(S^{n-1}) x^{-\alpha}$, and condition \eqref{*1.1} given by \eqref{eq:mar10}. This yields 
\begin{equation}
\label{eq:MFR-bis}
m(F_R)-m(F) \leq 2\constante R
\end{equation}
as long as 
\begin{equation}
\label{eq:range2-bis}
R^\alpha\geq \frac{\sig}{\alpha\delta},\med \mbox{ and }\med \frac{3\sigma(S^{n-1})}{R^\alpha} \exp\left(\frac{3\sigma(S^{n-1})g(2-\alpha)}{R^\alpha}\right)\leq 1/2-\delta. 
\end{equation}
Then with $x=4\constante R$, \eqref{eq:mar10}, \eqref{eq:MFR-bis} yield 
$$
P(F_R-m(F)\geq x)\leq  \frac{(2+\varepsilon)\sigma(S^{n-1})}{R^\alpha} \exp\left(\frac{(2+\varepsilon)\sigma(S^{n-1})g(2-\alpha)}{R^\alpha}\right)
$$
as long as \eqref{eq:range2-bis} holds. Together with \eqref{eq:tech2} and \eqref{eq:gamma-stable}, this  gives 
$$ 
P(F-m(F)\geq x)\leq \frac{\sigma(S^{n-1})}{R^\alpha}\left(\frac 1\alpha+(2+\varepsilon) \exp\left(\frac{(2+\varepsilon)\sigma(S^{n-1})g(2-\alpha)}{R^\alpha}\right)\right) 
$$
for $x=4\constante R$, that is \eqref{eq:ii} as long as \eqref{eq:range2-bis} holds. This 
 latter condition can be rewritten for $\delta\in (0,1/2)$ and $b>3$:
$$
\frac 1{2-\alpha} 
 \log \frac 1{2-\alpha} 
 \geq \max 
 \left(\frac 1{\alpha \delta b}, \left(\frac{3}{b(1/2-\delta)} 
 \right)^{\frac b{b-3}}\right)
$$ 
 which is obviously true for $\alpha$ close enough to $2$ since $\alpha, b$ are bounded below and $\disp\frac b{b-3}$ is bounded above. 
\end{Proof}


\footnotesize



\def\cprime{$'$}

\bigskip 

\sc
\noindent
J.-C. B.: Laboratoire de Math\'ematiques,
 Universit\'e de La Rochelle, Avenue Michel Cr\'epeau, 
 17042 La Rochelle, France. 
\\ 
{\tt jcbreton@univ-lr.fr}

\bigskip
\noindent C.H.: Laboratoire d'Analyse et de Math\'ematiques Appliqu\'ees, 
CNRS UMR 8050, Universit\'e Paris XII, 94010 Cr\'eteil Cedex, France, and 
School of Mathematics,
Georgia Institute of Technology,
Atlanta, Ga 30332 USA. 
\\ 
{\tt houdre@math.gatech.edu}

\bigskip
\noindent
N.P.: Laboratoire de Math\'ematiques, 
 Universit\'e de La Rochelle, Avenue Michel Cr\'epeau, 
 17042 La Rochelle, France. 
\\ 
{\tt nprivaul@univ-lr.fr}

\end{document}